\documentclass[11pt]{article}
\usepackage{amsfonts}
\usepackage{amsfonts,mathrsfs,amssymb,amsthm,mathptm}
\usepackage{amsmath,amscd}
\usepackage{mathptm,pslatex}
\usepackage{color}
\usepackage[dvips]{graphicx}
\usepackage[all]{xy}
\usepackage{graphicx}
\usepackage{fancyhdr}

\begin{document}
\newtheorem{Def}{Definition}[section]
\newtheorem{Bsp}[Def]{Example}
\newtheorem{Prop}[Def]{Proposition}
\newtheorem{Theo}[Def]{Theorem}
\newtheorem{Lem}[Def]{Lemma}
\newtheorem{Koro}[Def]{Corollary}
\theoremstyle{definition}
\newtheorem{Rem}[Def]{Remark}

\newcommand{\add}{{\rm add}}
\newcommand{\gd}{{\rm gl.dim}}
\newcommand{\dm}{{\rm dom.dim}}
\newcommand{\E}{{\rm E}}
\newcommand{\Mor}{{\rm Morph}}
\newcommand{\End}{{\rm End}}
\newcommand{\ind}{{\rm ind}}
\newcommand{\rsd}{{\rm res.dim}}
\newcommand{\rd} {{\rm rep.dim}}
\newcommand{\ol}{\overline}
\newcommand{\overpr}{$\hfill\square$}
\newcommand{\rad}{{\rm rad}}
\newcommand{\soc}{{\rm soc}}
\renewcommand{\top}{{\rm top}}
\newcommand{\pd}{{\rm proj.dim}}
\newcommand{\id}{{\rm inj.dim}}
\newcommand{\fld}{{\rm flat.dim}}
\newcommand{\Fac}{{\rm Fac}}
\newcommand{\Gen}{{\rm Gen}}
\newcommand{\fd} {{\rm fin.dim}}
\newcommand{\DTr}{{\rm DTr}}
\newcommand{\cpx}[1]{#1^{\bullet}}
\newcommand{\D}[1]{{\mathscr D}(#1)}
\newcommand{\DF}[1]{{\mathscr D}_F(#1)}
\newcommand{\Dz}[1]{{\mathscr D}^+(#1)}
\newcommand{\Df}[1]{{\mathscr D}^-(#1)}
\newcommand{\Db}[1]{{\mathscr D}^b(#1)}
\newcommand{\C}[1]{{\mathscr C}(#1)}
\newcommand{\Cz}[1]{{\mathscr C}^+(#1)}
\newcommand{\Cf}[1]{{\mathscr C}^-(#1)}
\newcommand{\Cb}[1]{{\mathscr C}^b(#1)}
\newcommand{\Dc}[1]{{\mathscr D}^c(#1)}
\newcommand{\K}[1]{{\mathscr K}(#1)}
\newcommand{\Kz}[1]{{\mathscr K}^+(#1)}
\newcommand{\Kf}[1]{{\mathscr  K}^-(#1)}
\newcommand{\Kb}[1]{{\mathscr K}^b(#1)}
\newcommand{\modcat}{\ensuremath{\mbox{{\rm -mod}}}}
\newcommand{\Modcat}{\ensuremath{\mbox{{\rm -Mod}}}}

\newcommand{\stmodcat}[1]{#1\mbox{{\rm -{\underline{mod}}}}}
\newcommand{\pmodcat}[1]{#1\mbox{{\rm -proj}}}
\newcommand{\imodcat}[1]{#1\mbox{{\rm -inj}}}
\newcommand{\Pmodcat}[1]{#1\mbox{{\rm -Proj}}}
\newcommand{\Imodcat}[1]{#1\mbox{{\rm -Inj}}}
\newcommand{\opp}{^{\rm op}}
\newcommand{\otimesL}{\otimes^{\rm\mathbf {L}}}
\newcommand{\rHom}{{\rm\mathbf {R}}{\rm Hom}\,}
\newcommand{\projdim}{\pd}
\newcommand{\Hom}{{\rm Hom}}
\newcommand{\Coker}{{\rm Coker}}
\newcommand{\Ker  }{{\rm Ker}}
\newcommand{\Cone }{{\rm Con}}
\newcommand{\Img  }{{\rm Im}}
\newcommand{\Ext}{{\rm Ext}}
\newcommand{\StHom}{{\rm \underline{Hom}}}

\newcommand{\gm}{{\rm _{\Gamma_M}}}
\newcommand{\gmr}{{\rm _{\Gamma_M^R}}}
\newcommand{\Pf}[1]{{\mathscr P}^{\infty}(#1)}

\def\vez{\varepsilon}\def\bz{\bigoplus}  \def\sz {\oplus}
\def\epa{\xrightarrow} \def\inja{\hookrightarrow}

\newcommand{\lra}{\longrightarrow}
\newcommand{\lraf}[1]{\stackrel{#1}{\lra}}
\newcommand{\llaf}[1]{\stackrel{#1}{\longleftarrow}}
\newcommand{\ra}{\rightarrow}
\newcommand{\dk}{{\rm dim_{_{k}}}}
\newcommand{\Add}{{\rm Add }}
\newcommand{\colim}{{\rm colim\, }}
\newcommand{\limt}{{\rm lim\, }}
\newcommand{\Tor}{{\rm Tor}}
\newcommand{\Cogen}{{\rm Cogen}}
\newcommand{\Tria}{{\rm Tria}}
\newcommand{\tria}{{\rm thick}}

{\Large \bf
\begin{center}
Recollements of derived categories II:  Algebraic
$K$-theory
\end{center}}

\medskip
\centerline{\textbf{Hongxing Chen} and \textbf{Changchang Xi}$^*$}

\renewcommand{\thefootnote}{\alph{footnote}}
\setcounter{footnote}{-1} \footnote{ $^*$ Corresponding author.
Email: xicc@cnu.edu.cn; Fax: 0086 10 58808202.}
\renewcommand{\thefootnote}{\alph{footnote}}
\setcounter{footnote}{-1} 
\renewcommand{\thefootnote}{\alph{footnote}}
\setcounter{footnote}{-1} 

\begin{abstract}

For a recollement of derived module categories of rings, we provide sufficient conditions to guarantee the additivity formula of higher algebraic $K$-groups of the rings involved, and establish a long Mayer-Vietoris exact sequence of higher algebraic $K$-groups for homological exact contexts introduced in the first paper of this series. Our results are then applied to recollements induced from homological ring epimorphisms and noncommutative localizations. Consequently, we get an infinitely long Mayer-Vietoris exact sequence of $K$-theory for Milnor squares, re-obtain a result of Karoubi (Corollary \ref{Karoubi}) on localizations and a result on generalized free products pioneered by Waldhausen and developed by Neeman and Ranicki. In particular, we describe algebraic $K$-groups of the free product of two groups over a regular coherent ring as the ones of the noncommutative tensor product of an exact context. This yields a new description of algebraic $K$-theory of infinite dihedral group.
\end{abstract}

\tableofcontents

\section{Introduction}\label{Section 1}
Algebraic $K$-theory of rings and algebras in the sense of Quillen (see \cite{Quillen}) collects elaborate invariants for
rings, groups and algebras. One of the most fundamental and important questions in this
theory is to understand and calculate these invariants: algebraic
$K$-groups $K_n$ of rings, which are closely connected with
Hochschild homologies $HH_n$ and with cyclic homologies $HC_n^-$ of
rings by Chern characters on higher $K$-theory (see \cite[Chapter
6]{rosenberg} for a survey). In computation of higher algebraic
$K$-groups of rings, Quillen, Suslin and many others have made important contributions
in the cases of finite fields, algebraically closed fields and certain integral domains (see \cite{qff} and the references in \cite{rosenberg}). For arbitrary rings, however,
the question is too hard, and remains little to be known, though general,
abstract algebraic $K$-theories have been explosively developed in
the last a few decades. In order to understand these algebraic $K$-groups $K_n(R)$ for arbitrary rings $R$, it is reasonable to investigate relationship between these $K$-groups of different rings which are linked in certain nice ways.

Along this direction are there some interesting and remarkable investigations in the literature. 
For example, if two rings are derived equivalent, then they have isomorphic $K_n$-groups by a result of  Dugger and Shipley (see \cite{DS}). For a homological noncommutative localization $\lambda: R\ra S$ of rings, Neeman and Ranicki established a long exact sequence of algebraic $K$-groups of $R$ and $S$ (see \cite{nr}). As an application of this result, Ranicki gave a new interpretation of Waldhausen's result on algebraic $K$-theory of generalized free products from the viewpoint of noncommutative localizations
(see \cite{wald1,Ranicki}). Later, Krause relaxed noncommutative locations to homological ring epimorphisms with the property that the chain map lifting problem has a positive answer, and established the same long exact sequence of $K$-groups (see \cite{kr}). Recently, we show in \cite{xc4} that, for a homological ring epimorphism $\lambda: R\ra S$, if the left $R$-module $S$ has a finite-type resolution, then $K_n(R)$ is the direct sum of $K_n(S)$ and $K_n(\mathbf R)$ where $\mathbf{R}$ is a Waldhausen category determined by $\lambda$. This result is then applied to study algebraic $K$-groups of endomorphism rings, matrix subrings and rings with idempotent ideals (see \cite{xc4} for detail).

Another useful type of natural linkages among rings is recollements of derived module categories, which were introduced by Beilinson, Bernstein and Deligne in \cite{BBD} for triangulated categories. Roughly speaking, a recollement consists of three derived (or triangulated) categories linked by two triangle functors both of which have left and right adjoint functors. The notion of recollements is an analogue of exact sequences for derived (or triangulated) categories, which generalizes derived equivalences and is closely related to homological ring epimorphisms. Here, a natural question is whether and when the additivity formula still holds true for algebraic $K$-groups of rings involved in an \emph{arbitrary} recollement of derived module categories. Namely, we consider the following question:

\medskip
{\bf Question.} Let $R$, $S$ and $T$ be rings with identity. Suppose
that there is a recollement among the derived module categories
$\D{T}$, $\D{R}$ and $\D{S}$ of the rings $T$, $R$ and $S$:

$$\xymatrix@C=1.2cm{\D{S}\ar[r]&\D{R}\ar[r]
\ar@/^1.2pc/[l]\ar@/_1.2pc/[l]
&\D{T}\ar@/^1.2pc/[l]\ar@/_1.2pc/[l]}\vspace{0.3cm}.$$
When does the following additivity isomorphism hold true:
$$K_n(R)\simeq K_n(S)\oplus K_n(T)\quad \mbox{for each}\;\; n\in
\mathbb{N}?$$
where $K_n(R)$ always means the $n$-th algebraic $K$-group of the ring $R$.

\medskip
This question may trace back to the work of Berrick and Keating on $K$-theory of the matrix rings of $2$ by $2$ triangular matrices (see \cite{bk}), where they described the $K_n$-group of a triangular matrix ring as the direct sum of the ones of rings in the diagonal. Recently,  we show in \cite{xc4} that, if the ideal $ReR$ of a ring $R$ generated by an idempotent element $e$ has a special finite-type resolution, then the $K_n$-group of $R$ is the direct sum of the ones of $R/ReR$ and $eRe$. In both cases, we do have recollements of derive module categories of rings and additivity formula of $K$-groups. However, the isomorphism $K_0(R)\simeq K_0(S)\oplus K_0(T)$ does not have to hold for arbitrary recollements. This can been seen by an example in \cite[Section 8, Remark (2)]{xc1}. So, answers to the above question seem to be mysterious.

In this paper we shall apply representation-theoretic methods to investigate the above question in detail and establish an additivity formula for higher algebraic $K$-groups of rings involved in recollements with a compactness condition. Thus we provide a general answer to the above question. Further, we apply our result to homological ring epimorphisms, exact contents, extensions and free products of groups.

In dealing with the above-mentioned question, a number of technical obstacles occur: To understand $K_n(R)$, one has to choose certain models of algebraic $K$-theory space $K(R)$ of $R$, up to homotopy equivalence. For example, the usual favorite models are the category of finitely generated projective $R$-modules (see \cite{Quillen, rosenberg, xc4}), the category of bounded complexes of finitely generated projective $R$-modules (see \cite{wald, kr, ne2, nr, Sch2}), and certain full subcategories of the category of complexes over $R$ with countable direct sums (see \cite{nr, DS}). When comparing algebraic $K$-theory of different rings, one has first to fix a suitable model to define $K$-theory, and then to find exact functors compatible with the chosen model. Unfortunately, given an arbitrary recollement of derived module categories, nothing is known about the concrete forms of the six triangle functors. This means that it would be quite difficult to find a suitable model for all three rings in the recollement such that the given six functors can induce compatible functors on the model for all rings and  connect $K$-theory space $K(R)$ with $K$-spaces $K(S)$ and $K(T)$ in a reasonable way. Hence the methods used in \cite{xc4, nr, ne2, kr} actually does not work any more for the present case, and therefore some new ideas are necessary for attacking the above question.

To overcome these obstacles, we pass to differential graded (dg) algebras and introduce a new definition of algebraic $K$-theory spaces for dg algebras, which captures the usual definition of algebraic $K$-theory spaces of ordinary rings up to homotopy equivalence. Our definition of $K$-theory spaces is a modification of Schlichting's definition in \cite{Sch1}, and excludes the potential set-theoretic difficulties in the corresponding definition given by Dugger and Shipley in \cite{DS}. Also, this new definition gives much freedom for choices of compatible functors among models that define $K$-theory. Under a compactness assumption, we can identify $K(S)$ and $K(T)$ with algebraic $K$-theory spaces of dg endomorphism algebras $\mathbb{S}$ and $\mathbb{T}$ of perfect complexes over $R$, respectively. After a systematical study on homotopy equivalences of $K$-theory spaces related to perfect dg modules, we establish decomposition formulas for algebraic $K$-groups of dg algebras. Particularly, this leads to the following main result in this paper.

\begin{Theo}\label{new-theorem}
Let $R$, $S$ and $T$ be rings with identity. Suppose
that there is a recollement among the derived module categories
$\D{T}$, $\D{R}$ and $\D{S}$ of the rings $T$, $R$ and $S$:

$$\xymatrix@C=1.2cm{\D{S}\ar[r]^-{i_*}&\D{R}\ar[r]
\ar@/^1.2pc/[l]\ar@/_1.2pc/[l]
&\D{T}.\ar@/^1.2pc/[l]\ar@/_1.2pc/[l]}\vspace{0.3cm}$$
If $i_*(S)$ is quasi-isomorphic to a bounded complex of finitely
generated projective $R$-modules, that is, $i_*(S)$ is compact in $\D R$, then
$$K_n(R)\simeq K_n(S)\oplus K_n(T)\quad \mbox{for all}\;\; n\in
\mathbb{N}.$$
\end{Theo}

\medskip
We remark that, under the compactness condition in Theorem \ref{new-theorem}, it is not difficult to prove that $K_0(S)$ is a direct summand of $K_0(R)$. However, the key point here, which seems to be highly non-trivial, is to prove that an additive complement to $K_n(S)$ is just $K_n(T)$ for all $n\ge 0$. Also, we note that Theorem \ref{new-theorem} cannot be extended to dg algebras because derived equivalences of dg algebras do not preserve algebraic $K$-groups, as pointed out by an example in \cite{DS}.

\medskip
First, we apply Theorem \ref{new-theorem} to recollements of derived module categories arising from homological ring epimorphisms.

Recall that a ring epimorphism $\lambda: R\ra S$ is said to be \emph{homological} if $\Tor_j^R(S,S)=0$ for all $j>0$. An $R$-module $M$ has a \emph{finite-type resolution} provided that there is a finite
projective resolution by finitely generated projective $R$-modules, that is, there is an exact sequence $0\ra P_m\ra \cdots\ra P_1\ra P_0\ra M\ra 0$ for some $m\in \mathbb{N}$ such that all $R$-modules $P_j$ are projective and finitely generated.


\begin{Koro} \label{corollary1}
Suppose that $\lambda: R\to S$ is a homological ring epimorphism which induces a recollement of derived module categories of rings $T,R,S$:
$$
\xymatrix@C=1.2cm{\D{S}\ar[r]^-{i_*}&\D{R}\ar[r]
\ar@/^1.2pc/[l]\ar@/_1.2pc/[l]
&\D{T}\ar@/^1.2pc/[l]\ar@/_1.2pc/[l]}\vspace{0.3cm}$$
where $i_*$ is the restriction functor induced from $\lambda$. If $_RS$ or $S_R$ has a finite-type resolution, then
$$K_n(R)\simeq K_n(S)\oplus K_n(T)\quad \mbox{for all}\;\; n\in
\mathbb{N}.$$
\end{Koro}

We should note that not every homological ring epimorphisms $R\ra S$ can induce a recollement of derived module categories of rings because the Verdier quotient of $\D{R}$ by $\D{S}$ may not be realized as the derived category of a usual ring. This can be seen by the counterexample given by Bernhard Keller to the Telescope conjecture. Comparing Corollary \ref{corollary1} with \cite[Theorem 1.1]{xc4}, we see that the conclusion of Corollary \ref{corollary1}, under the assumption of existence of a recollement, is quite strong. In fact, by \cite[Theorem 1.1]{xc4}, we have  $K_n(R)\simeq K_n(S)\oplus K_n(R,\lambda)$ for all $n\in\mathbb{N}$, where $K_n(R,\lambda)$ is the $n$-th algebraic $K$-group of the category $\mathbf{R}$ mentioned before, while Corollary \ref{corollary1} describes $K_n(R,\lambda)$ explicitly as the $K_n$-group of a ring $T$ if such a ring $T$ exists. Moreover, since stratifying ideals give rise to recollements of derived module categories, Corollary \ref{corollary1} also generalizes \cite[Corollary 1.3]{xc4}.

\smallskip
Next, we consider $K$-theory of reollements arising from exact contexts introduced in the first paper of this series (see \cite{xc3}). This kind of recollements involves noncommutative localizations in ring theory, which occur often in algebraic topology and representation theory (see \cite{schofield, Ranicki}).

Let $R$, $S$ and $T$ be associative rings with identity, and let $\lambda:R\to S$ and $\mu:R\to T$ be ring homomorphisms. Suppose that $M$ is an $S$-$T$-bimodule together with an element $m\in M$. The quadruple $(\lambda, \mu, M, m)$ is called an \emph{exact context} if the following sequence
$$
0\lra R\lraf{(\lambda,\,\mu)}S\oplus T
\lraf{\left({\cdot\,m\,\atop{-m\,\cdot}}\right)}M\lra 0
$$
is an exact sequence of abelian groups, where $\cdot m$ and $m \cdot$ denote the right and left multiplication by $m$ maps, respectively.
An exact context $(\lambda, \mu, M, m)$ is called an \emph{exact pair} if $M=S\otimes_RT$ and $m=1\otimes 1$. In this case we simply say that $(\lambda,\mu)$ is an exact pair. The exact context $(\lambda,\mu, M, m)$ is said to be \emph{homological} if $\Tor^R_i(T,S)=0$ for all $i\geq 1$.

For each exact context $(\lambda,\mu, M, m)$, we associate it with a new ring $T\boxtimes_RS$, called the \emph{noncommutative tensor product} of $(\lambda. \mu, M, m)$ in \cite[Section 4.1]{xc3}, which is a generalization of the usual tensor products over commutative rings, and captures coproducts of rings and dual extensions.

For a homological exact context $(\lambda,\mu, M, m)$, we have the following long Mayer-Vietoris sequence which links algebraic $K$-groups of the rings $R,\,S,\,T$ and $T\boxtimes_RS$ together.

\begin{Theo}\label{K-theory}
Let $(\lambda, \mu, M, m)$ be a homological exact context. Then the following statements hold true:

$(1)$ There exists a long exact sequence of algebraic $K$-groups:
$$
\xymatrix{\cdots\lra K_{n+1}(T\boxtimes_RS)\lra
K_n(R)\ar[rr]^-{\big(-K_n(\lambda),\, K_n(\mu)\big)} && K_n(S)\oplus
K_n(T) \ar[rr]^-{\left({K_n(\rho)\; \atop{K_n(\phi)}}\right)} &&
K_n(T\boxtimes_RS)\lra K_{n-1}(R)\lra }
$$
$$\cdots\lra K_0(R)\lra K_0(S)\oplus K_0(T)\lra
K_0(T\boxtimes_RS)$$ for all $n\in \mathbb{N}.$

$(2)$ If the left $R$-module $S$ or the right $R$-module $T$ has a finite-type resolution, then
$K_n(R)\oplus K_n(T\boxtimes_RS)\simeq K_n(S)\oplus K_n(T)$ for all $n\in
\mathbb{N}.$
\end{Theo}

Since a Milnor square of rings provides a typical exact pair (see \cite[Example (3), Section3; Corollary 4.3]{xc3}, we have the following long Mayer-Vietoris exact sequence which extends and amplifies the $K$-theory sequence in \cite{land}.

 \begin{Koro}
 Given a pullback square of rings and surjective homomorphisms
 $$
\xymatrix{
R \ar[d]_-{i_2} \ar[r]^-{i_1} & R_1 \ar[d]^-{j_1}\\
R_2\ar[r]^-{j_2} &  \; \, R',}$$ if $\Tor_j^{R}(R_2,R_1)=0$ for all $j>0$, then there is a long Mayer-Vietoris exact sequence:
$$\cdots\lra K_{n+1}(R')\lra
K_n(R)\lra  K_n(R_1)\oplus
K_n(R_2) \lra K_n(R')\lra K_{n-1}(R)\lra $$
$$\cdots\lra K_0(R)\lra K_0(R_1)\oplus K_0(R_2)\lra
K_0(R')$$ for all $n\in \mathbb{N}.$
\label{milnorsquare}
 \end{Koro}

As another consequence of Theorem \ref{K-theory}, we obtain the following result on ring extensions.

\begin{Koro}\label{ring extension}
Suppose that $R\subseteq S$ is an extension of rings, that is, $R$ is a subring of the ring $S$ with the same identity. Let $S'$ be the endomorphism ring of the left $R$-module $S/R$. If the left $R$-module $S$ is projective and finitely generated, then
$$K_n(R)\oplus K_n(S'\boxtimes_RS)\simeq K_n(S)\oplus K_n(S')\quad \mbox{for all}\;\; n\in
\mathbb{N},$$ where $S'\boxtimes_RS$ is the noncommutative tensor product of an exact context defined by the extension.
\end{Koro}

A rather striking application of Theorem \ref{K-theory} is that algebraic $K$-groups of the free products of finite groups can be characterized by noncommutative tensor products which have finite ranks over ground rings, while the group rings of free products usually have infinite ranks.

Let $H$ and $G$ be two groups, and let $RH$ and $RG$ be the group rings of $H$ and $G$ over a ring $R$, respectively. Then the canonical maps from $R$ to $RH$ and $RG$ can be completed into an exact context (see Section \ref{groupring} for details) and the associated noncommutative tensor product $RH\boxtimes_RRG$ can be described explicitly as follows:

As an abelian group, $RH\boxtimes_RRG$ coincides with the group ring $R(H\times G)$ of the direct product $H\times G$ over $R$. Thus $RH\boxtimes_RRG$ is a finitely generated free $R$-module if $G$ and $H$ are finite. As an associative ring, it admits the following multiplication:
$$r(h,g)=(h,g)r \quad \mbox{and}\quad (h,g)(h',g')=(h,gg')+(hh',g')-(h,g'),$$
where  $r\in R$, $h,h'\in H$ and $g,g'\in G$,

\smallskip
Recall that the free product of $H$ and $G$, denoted by $H*G$, is the coproduct of $H$ and $G$ in the category of groups. In general, the free product of finite groups may be infinite. For example, the free product of two cyclic groups of order $2$ is the infinite dihedral group $D_\infty$.

We say that a ring $R$ is \emph{regular coherent} if any finitely presented left $R$-module has a finite-type resolution. A typical example of regular coherent rings is the ring of integers.

\smallskip
The following corollary follows from Theorem \ref{K-theory} together with \cite[Theorems 1 and 4]{wald1}, which reduces surprisingly $K$-theory of group rings of infinite $R$-rank to the one of rings of finite $R$-rank.

\begin{Koro}\label{GP}
Let $R$ be a regular coherent ring and let $H$ and $G$ be two groups. Then
$$ K_n\big(R(H*G)\big)\simeq K_n\big(RH\boxtimes_RRG\big)\simeq K_n\big(RG\boxtimes_RRH\big) \quad \mbox{for all}\;\; n\geq 1.$$
\end{Koro}

As a consequence of our methods, we get a new description of algebraic $K$-theory for infinite dihedral group $D_{\infty}$: For an arbitrary ring $R$, $K_n\big(R(D_\infty)\big)\simeq K_n(R\,\mathbb{Z}_2\boxtimes_R R\,\mathbb{Z}_2)\oplus
\widetilde{\mathfrak{Nil}}_{n-1}(R)$ for $n\ge 1$, where $\widetilde{\mathfrak{Nil}}_{n}(R)$ is the $n$-th reduced Nilgroup of $R$. This decomposition is different from the result in \cite{DKR}.

\medskip
This paper is organized as follows: In Section \ref{sect2}, we briefly recall some definitions and basic facts on triangulated categories, recollements and homological ring epimorphisms. In Section \ref{sect3}, we first recall the algebraic $K$-theories developed by Waldhausen for Waldhausen categories and Schlichting for Frobenius pairs, and then mention several fundamental theorems in algebraic $K$-theory of Frobenius pairs. In Section \ref{3.5}, we first introduce our definition of algebraic $K$-theory spaces for differential graded algebras, and then discuss homotopy equivalences of $K$-theory spaces constructed from perfect dg modules in detail. As a result, we establish a reduction  in Proposition \ref{dg-5} for calculation of algebraic $K$-groups of dg algebras. At the end of this section, we prove Theorem \ref{new-theorem} as well as Corollary \ref{corollary1}. In Section \ref{Exact context}, we apply our results to homological exact contexts, and prove Theorem \ref{K-theory} and Corollaries \ref{milnorsquare}, \ref{ring extension} and \ref{GP}.

In the third paper of this series, we shall discuss finitistic dimension theory for recollements of derived module categories of rings (see \cite{xc6}).

\section{Preliminaries\label{sect2}}
In this section, we shall fix notation employed
throughout the paper, and provide some basic facts for later proofs.

\subsection{General terminology and notation on categories} \label{subsection2.1}

Let $\mathcal C$ be an additive category.

We always assume that a full subcategory $\mathcal B$ of $\mathcal
C$ is closed under isomorphisms, that is, if $X\in {\mathcal B}$ and
$Y\in\cal C$ with $Y\simeq X$, then $Y\in{\mathcal B}$.

Given two morphisms $f: X\to Y$ and $g: Y\to Z$ in $\mathcal C$, we
denote the composite of $f$ and $g$ by $fg$ which is a morphism from
$X$ to $Z$, while given two functors $F:\mathcal {C}\to \mathcal{D}$
and $G: \mathcal{D}\to \mathcal{E}$ among three categories
$\mathcal{C}$, $\mathcal{D}$ and $\mathcal{E}$, we denote the
composite of $F$ and $G$ by $GF$ which is a functor from $\mathcal
C$ to $\mathcal E$.

Let $\Ker(F)$ and $\Img(F)$ be the kernel and image of the functor
$F$, respectively. That is, $\Ker(F):=\{X\in \mathcal{C}\mid
FX\simeq 0\}$ and $\Img(F):=\{Y\in \mathcal{D}\mid \exists\, X\in
\mathcal{C}, FX\simeq Y\}$. In particular, $\Ker(F)$ and $\Img(F)$
are closed under isomorphisms in $\mathcal{C}$ and $\mathcal{D}$,
respectively.

An additive functor $F:\mathcal{A}\to\mathcal{B}$ between two
additive categories $\mathcal{A}$ and $\mathcal{B}$ is called an
\emph{equivalence up to factors} if $F$ is fully faithful and each
object of $\mathcal{B}$ is isomorphic to a direct summand of the
image of an object of $\mathcal{A}$ under $F$.

Let $\mathcal{A}$ be a triangulated category and $\mathcal{X}$ a
full triangulated subcategory of $\mathcal{A}$. Then,
due essentially to Verdier, there exists a triangulated category
$\mathcal{A}/\mathcal{X}$, and a triangle functor $q:\mathcal{A}\to
\mathcal{A}/\mathcal{X}$ with $\mathcal{X}\subseteq\Ker(q)$ such
that $q$ has the following universal property: If $q':
\mathcal{A}\to\mathcal{T}$ is a triangle functor with
$\mathcal{X}\subseteq\Ker(q')$, then $q'$ factorizes uniquely
through $\mathcal{A}\lraf{q} \mathcal{A}/\mathcal{X}$ by
\cite[Theorem 2.18]{ne3}. The category $\mathcal{A}/\mathcal{X}$ is
called the \emph{Verdier quotient} of $\mathcal{A}$ by
$\mathcal{X}$, and the functor $q$ is called the \emph{Verdier
localization functor}. In this case, $\Ker(q)$ is the full
subcategory of $\mathcal{A}$ consisting of direct summands of all
objects in $\mathcal{X}$ (see \cite[Chapter 2]{ne3} for details). We remark that the objects of the category
$\mathcal{A}/\mathcal{X}$ are the same as the objects of
$\mathcal{A}$.

A sequence $\mathcal{A}\lraf{F}\mathcal{B}\lraf{G}\mathcal{C}$ of
triangle functors $F$ and $G$ between triangulated categories is
said to be \emph{exact} if the following four conditions are
satisfied:

$(i)$ The functor $F$ is fully faithful.

$(ii)$ The composite $GF:\mathcal{A}\to\mathcal{C}$ of $F$ and $G$
is zero.

$(iii)$ The image $\Img(F)$ of $F$ is equal to the kernel of $G$.

$(iv)$ The functor $G$ induces an equivalence from the Verdier
quotient of $\mathcal{B}$ by $\Img(F)$ to $\mathcal{C}$.

\medskip
Clearly, if $\mathcal{X}$ is closed under direct summands in
$\mathcal{A}$, then we have an exact sequence of triangulated
categories:
$$
\xymatrix{\mathcal{X}\ar@{^{(}->}[r] & \mathcal{A}\ar[r]^-{q}&
\mathcal{A}/\mathcal{X}}.$$

Let $\mathcal{T}$ be a triangulated category with small coproducts
(that is, coproducts indexed over sets exist in ${\mathcal T}$).

An object $U\in \mathcal{T}$ is said to be \emph{compact} if
$\Hom_\mathcal{T}(U,-)$ commutes with small coproducts in
$\mathcal{T}$. The full subcategory of $\mathcal{T}$ consisting of
all compact objects is denoted by $\mathcal{T}^c$.

For any non-empty class $\mathscr{S}$ of objects in $\mathcal{T}$,
we denote by $\Tria(\mathscr S)$ (respectively, $\tria(\mathscr S)$)
the smallest full triangulated subcategory of $\mathcal{T}$
containing $\mathscr S$ and being closed under small coproducts
(respectively, direct summands). If $\mathscr{S}$ consists of a single object $U$, then we simply write $\Tria(U)$ and $\tria(U)$ for
$\Tria(\{U\})$ and $\tria(\{U\}),$ respectively. The notation
$\Tria(\mathscr S)$ without referring to $\mathcal T$ will not cause
any confusions because this notation can be clarified from the
contexts of our considerations.

\smallskip
The following facts are in the literature (see \cite[Proposition
1.6.8]{ne3} and \cite[Section 2.1]{xc3}).

\begin{Lem}\label{literature}
$(1)$ If $\mathcal{T}_0$ is a full triangulated subcategory of
$\mathcal{T}$ such that $\mathcal{T}_0$ is closed under countable
coproducts, then $\mathcal{T}_0$ is closed under direct summands in
$\mathcal{T}$.

$(2)$ Let $\mathcal{T}'$ be a triangulated category with small
coproducts, and let $F: \mathcal{T}\ra \mathcal{T}'$ be a triangle
functor. If $F$ preserves small coproducts, then
$F(\Tria(U))\subseteq \Tria(F(U))$ for any $U\in\mathcal{T}$.
\end{Lem}

Special examples of triangulated categories are the derived module categories of (associative) rings with identity, which are of our particular interest in this paper. Now, let us fix some notation for rings.

Let $R$ be a ring with identity. We denote by $R\Modcat$ the category of all left $R$-modules. The
complex, homotopy and derived categories of $R\Modcat$ are usually denoted
by $\C{R}, \K{R}$ and $\D{R}$, respectively. It is well-known that both $\K{R}$ and $\D{R}$ are triangulated categories, and that $\D{R}=\Tria({_R}R)$. As usual we write $\Dc{R}$ for $\D{R}^c$, which is equal to the full subcategory of $\D{R}$ consisting all those complexes that are quasi-isomorphic to bounded complexes of finitely generated projective $R$-modules.

\subsection{Recollements and homological ring epimorphisms}\label{subsection2.2}

In this subsection, we recall the notion of recollements
introduced by Beilinson, Bernstein and Deligne (see \cite{BBD}), which is widely used in algebraic geometry and representation theory.
Some prominent examples of recollements can be constructed from certain homological
ring epimorphisms.

Let $\cal D$, $\mathcal{D'}$ and $\mathcal{D''}$ be triangulated
categories with shift functors denoted universally by [1].

We say that $\mathcal{D}$ is a \emph{recollement} of $\mathcal{D'}$
and $\mathcal{D''}$ if there are six triangle functors indicated in the
following diagram
$$\xymatrix{\mathcal{D''}\ar^-{i_*=i_!}[r]&\mathcal{D}\ar^-{j^!=j^*}[r]
\ar^-{i^!}@/^1.2pc/[l]\ar_-{i^*}@/_1.6pc/[l]
&\mathcal{D'}\ar^-{j_*}@/^1.2pc/[l]\ar_-{j_!}@/_1.6pc/[l]}$$ such
that:

$(1)$ The 4 pairs $(i^*,i_*),(i_!,i^!),(j_!,j^!)$ and $(j^*,j_*)$
are adjoint pairs of functors.

$(2)$ The 3 functors $i_*,j_*$ and $j_!$ are fully faithful.

$(3)$ The composite of two functors in each row is zero, that is,
$i^!j_*=0$ (and thus also $j^! i_!=0$ and $i^*j_!=0$).

$(4)$ There are 2 canonical triangles in $\mathcal D$ for each
object $X\in\mathcal{D}$:
$$
j_!j^!(X)\lra X\lra i_*i^*(X)\lra j_!j^!(X)[1],\qquad
i_!i^!(X)\lra X\lra j_*j^*(X)\lra i_!i^!(X)[1],
$$ where
$j_!j^!(X)\ra X$ and $i_!i^!(X)\ra X$ are counit adjunction maps,
and where $ X\ra i_*i^*(X)$ and $X\ra j_*j^*(X)$  are unit
adjunction maps.

It is known that, up to equivalence of categories, recollements of
triangulated categories are the same as torsion torsion-free triples
(TTF-triples) of triangulated categories (see, for example,
\cite{BI} and \cite[Section 2.3]{xc1} for details). In the following lemma we
mention some facts about recollements for later proofs.

\begin{Lem}\label{prop-recell} Suppose that the above recollement is
given. Then the following hold:

 $(a)$ The images of the three fully faithful functors $j_!$, $i_*$ and $j_*$ are closed
under direct summands in $\mathcal{D}$.

$(b)$ The Verdier quotients of $\mathcal{D}$ by the images of the
triangle functors $j_!$ and $i_*$ are equivalent to $\mathcal{D''}$
and $\mathcal{D'}$, respectively.

$(c)$ Assume that $\mathcal{D}$, $\mathcal{D}'$ and $\mathcal{D}''$
admit small coproducts. Then both $j_!$ and $i^*$ preserve compact
objects. Suppose further that $\mathcal{D}$ is compactly generated,
that is, there is a set $S$ of compact objects in $\mathcal D$ such
that $\Tria(S)=\mathcal D$, then $i_*$ preserves compact objects if
and only if so is $j^!$. In this case, we can obtain a ``half
recollement" of subcategories of compact objects:
$$\xymatrix@C=1.2cm{(\mathcal{D}'')^c\ar[r]^-{i_*}&\mathcal{D}^c\ar[r]^-{j^!}
\ar@/_1.2pc/[l]_-{i^*}&(\mathcal{D}')^c\ar@/_1.2pc/[l]_-{j_!}}\vspace{0.3cm}$$

\end{Lem}
Note that $(a)$ and $(b)$ follow from \cite[Chapter I, Proposition
2.6]{BI}, while $(c)$ follows from \cite[Chapter III, Lemma 1.2 (1)
and Chapter IV, Proposition 1.11]{BI}.

A typical example of recollements occurs in the following two situations.

(1) Recollements of derived module categories.

Let $R$ be a ring with an idempotent ideal $I=ReR$ for $e^2=e\in R$.
Suppose that $I$ is a \emph{stratifying ideal} of $R$, that is, the
multiplication map $Re\otimes_{eRe}eR\ra ReR$ is an isomorphism and
$\Tor^{eRe}_j(Re,eR)=0$ for $j\ge 1$, then there is a recollement of
derived module categories:

$$\xymatrix{\D{R/I}\ar^-{D(\pi_*)}[r]&\D{R}\ar^-{eR\otimes_R^{\mathbb{L}}-}[r]
\ar^-{\mathbb{R}\Hom_R(R/I,-)\;}@/^1.2pc/[l]\ar_-{(R/I)\otimes_R^{\mathbb
L}-}@/_1.6pc/[l] &\D{eRe}\ar_-{Re\otimes_{eRe}^{\mathbb
L}-}@/_1.6pc/[l]\ar^-{\;\mathbb{R}\Hom_{eRe}(eR,-)}@/^1.2pc/[l] }$$
where $D(\pi_*)$ is the restriction functor induced from the canonical surjection $\pi:R\to R/I$, and where $Re\otimes_{eRe}^{\mathbb{L}}-$ is the total left-derived
functor of $Re\otimes_{eRe}-$ and $\mathbb{R}\Hom_{eRe}(eR,-)$ is
the total right-derived functor of $\Hom_{eRe}(eR,-)$.  For more
details, we refer the reader to \cite{cps}.

In Section \ref{Exact context} and \cite{xc1, xc5, xc3}£¬ one may find more examples of recollements of derived module categories, which have not to be induced from idempotent elements.

\medskip
(2) Recollements of triangulated categories induced from ring epimorphisms.

Recall that a ring epimorphism $\lambda:R\to S$ is said to be \emph{homological} if
$\Tor^R_n(S, S)=0$ for all $n>0$ (see \cite{GL, nr}). This is also equivalent to that the
restriction functor $D(\lambda_*):\D{S}\to\D{R}$ is fully faithful.

According to \cite[Section 4]{NS}, for an arbitrary homological ring epimorphism, we obtain the following recollement of triangulated categories, of which the right-hand term is not necessarily the derived category of an ordinary ring.

\begin{Lem}\label{homological}
Let $\lambda:R\to S$ be a homological ring epimorphism. Then there
is a recollement of triangulated categories:
$$\xymatrix@C=1.2cm{\D{S}\ar[r]^-{i_*}&\D{R}\ar[r]^-{j^!}
\ar@/^1.2pc/[l]\ar@/_1.2pc/[l]_-{i^*}
&{\rm{Tria}}({_R}\cpx{Q})\ar@/^1.2pc/[l]\ar@/_1.2pc/[l]_-{j_!}}\vspace{0.3cm}$$
where $\cpx{Q}$ is the two-term complex $0\ra R\lraf{\lambda} S\ra
0$ with $R$ and $S$ in degrees $0$ and $1$, respectively, and where
$j_!$ is the canonical embedding and
$j^!=\cpx{Q}\otimesL_R-,\,
i^*=S\otimesL_R-, \,i_*=D(\lambda_*).$
\end{Lem}

Thus, if we define $\mathscr{Y}:=\{Y\in\D{R}\mid \Hom_{\D{R}}(X,
Y)=0 \mbox{\, for any \;} X\in\Tria({_R}\cpx{Q})\},$ then it follows
from Lemma \ref{homological} that
$$\mathscr{Y}=\{Y\in \D{R}\mid \Hom_{\D{R}}(\cpx{Q},
Y[n])=0\mbox{\,for \;} n\in\mathbb{Z}\}=\{Y\in\D{R}\mid
\cpx{Q}\otimesL_RY=0\},$$ and that $i_*$ induces an equivalence
$\D{S}\lraf{\simeq} \mathscr{Y}$.

\smallskip
Finally, we point out that if a homological ring epimorphism induces a recollement of derived module categories of rings, then it also gives a recollement of derived module categories of opposite rings, though the categories $\D R$ and $\D
{R^{\opp}}$ for a ring $R$ may not be triangle equivalent. This fact will be used in the proof of Corollary \ref{corollary1}.

\begin{Lem}\label{dual}
Let $\lambda:R\to S$ be a homological ring epimorphism. Then the
following are equivalent for a ring $T$:

$(1)$ There is a recollement of derived categories:

$$\xymatrix@C=1.2cm{\D{S}\ar[r]^-{D(\lambda_*)}&\D{R}\ar[r]
\ar@/^1.2pc/[l]\ar@/_1.2pc/[l]
&\D{T}\ar@/^1.2pc/[l]\ar@/_1.2pc/[l]}\vspace{0.3cm}$$

$(2)$ There is a recollement of derived categories:
$$\xymatrix@C=1.2cm{\D{S\opp}\ar[r]^-{D(\lambda_*)}&\D{R\opp}\ar[r]
\ar@/^1.2pc/[l]\ar@/_1.2pc/[l]
&\D{T\opp}\ar@/^1.2pc/[l]\ar@/_1.2pc/[l]}\vspace{0.3cm}$$
\end{Lem}

{\it Proof.} Observe that if $\lambda:R\ra S$ is a homological ring
epimorphism, then so is the map $\lambda: R\opp\to S\opp$ by
\cite[Theorem 4.4]{GL}. Moreover, it follows from \cite[Corollary
3.4]{NS} that $(1)$ holds if and only if there is a complex
$\cpx{P}\in\Cb{\pmodcat R}$ such that
$\Tria(\cpx{P})=\Tria({_R}\cpx{Q})$, $\End_{\D R}(\cpx{P})\simeq T$
and $\Hom_{\D R}(\cpx{P}, \cpx{P}[n])=0$ for any $n\neq 0$, where
$\cpx{Q}$ is the complex $0\ra R\ra S\ra 0$. However, for such a
complex $\cpx{P}$, we always have
$$\Hom_{\D{R\opp}}(\cpx{P}{^*}, \cpx{P}{^*}[n])\simeq\Hom_{\D
R}(\cpx{P}, \cpx{P}[n])\;\;\mbox{for all}\;\; n\in\mathbb{Z},$$
where $\cpx{P}{^*}:=\Hom_R(\cpx{P}, R)\in\Cb{\pmodcat {R\opp}}.$ So,
to prove that $(1)$ and $(2)$ are equivalent, it is enough to prove
the following statement:

If $\cpx{P}\in\Cb{\pmodcat R}$ such that
$\Tria(\cpx{P})=\Tria({_R}\cpx{Q})$, then
$\Tria(\cpx{P}{^*})=\Tria(\cpx{Q}_R).$

In fact, let $\cpx{P}$ be such a complex and define
$$\mathscr{Y}':=\{Y\in \D{R\opp}\mid \Hom_{\D{R\opp}}(X, Y)=0
\mbox{\, for \;} X\in{\rm{Tria}}(\cpx{P}{^*})\}.$$ Since
$\cpx{P}\in\Cb{\pmodcat R}$, we have $\cpx{P}{^*}\in\Cb{\pmodcat
{R\opp}}$. It follows from \cite[Lemma 2.8]{xc1} that there is a
recollement:
$$\xymatrix@C=1.2cm{\mathscr{Y}'\ar[r]^-{\mu}&\D{R\opp}\ar[r]
\ar@/^1.2pc/[l]\ar@/_1.2pc/[l]
&{\rm{Tria}}(\cpx{P}{^*})\ar@/^1.2pc/[l]\ar@/_1.2pc/[l]}\vspace{0.3cm}$$
where  $\mu$ is the inclusion. This implies that
$$(a)\quad \;\Tria(\cpx{P}{^*})=\{X\in \D{R\opp}\mid \Hom_{\D{R\opp}}(X, Y)=0
\mbox{\, for \;} Y\in\mathscr{Y}'\}.$$ Furthermore, we remark that
$$\mathscr{Y}'=\{Y\in \D{R\opp}\mid \Hom_{\D{R\opp}}(\cpx{P}{^*},
Y[n])=0\mbox{\;for \;} n\in\mathbb{Z}\}=\{Y\in \D{R\opp}\mid
{\mathbb R}\Hom_{R\opp}(\cpx{P}{^*}, Y)=0\},$$ and that $${\mathbb
R}\Hom_{R\opp}(\cpx{P}{^*}, -)\simeq -\otimesL_R\cpx{P}:\;
\D{R\opp}\lra\D{\mathbb Z}$$ by \cite[Section 2.1]{xc3}. Thus
$\mathscr{Y}'=\{Y\in \D{R\opp}\mid Y\otimesL_R\cpx{P}= 0\}$.
However, by Lemma \ref{literature} (2), for a given $Y\in
\D{R\opp}$, the left-derived tensor functor $Y\otimesL_R-:\D{R}\to
\D{\mathbb{Z}}$ sends $\Tria(\cpx{Q})$ (respectively,
$\Tria({_R}\cpx{P})$) to zero if and only if $Y\otimesL_R\cpx{Q}= 0$
(respectively, $Y\otimesL_R\cpx{P}= 0$). Since
$\Tria(\cpx{P})=\Tria({_R}\cpx{Q})$ by assumption, we certainly
obtain $\mathscr{Y}'=\{Y\in \D{R\opp}\mid Y\otimesL_R\cpx{Q}=0\}.$

Since $\lambda: R\opp\to S\opp$ is  a homological ring
epimorphism, we obtain another recollement by Lemma
\ref{homological}:
$$\xymatrix@C=1.2cm{\D{S\opp}\ar[r]^-{D(\lambda_*)}&\D{R\opp}\ar[r]^-{G}
\ar@/^1.2pc/[l]\ar@/_1.2pc/[l]
&{\rm{Tria}}(\cpx{Q}_R)\ar@/^1.2pc/[l]\ar@/_1.2pc/[l]_-{F}}\vspace{0.3cm}$$
where $F$ is the inclusion and $G$ is the tensor functor
$-\otimesL_R\cpx{Q}$. This implies that $\Img\big(D(\lambda_*)\big)=\Ker(G)$
and
$$(b)\quad \;\Tria(\cpx{Q}_R)=\{X\in \D{R\opp}\mid \Hom_{\D{R\opp}}(X, Y)=0
\mbox{\, for \;} Y\in\Ker(G)\}.$$ Since $\mathscr{Y}'=\Ker(G)$, we
conclude from $(a)$ and $(b)$ that
$\Tria(\cpx{P}{^*})=\Tria(\cpx{Q}_R)$. This finishes the proof of
Lemma \ref{dual}.  $\square$

\section{Algebraic $K$-theory \label{sect3}}

In this section, first, we briefly recall some basics on algebraic
$K$-theory of Waldhausen categories and Frobenius pairs developed in
\cite{wald} and \cite{Sch1}, respectively. And then we discuss
algebraic $K$-theory of differential graded algebras and prove a few
facts as preparations for proofs of the main results.

\subsection{$K$-theory spaces of small Waldhausen categories}\label{subsection3.1}

Let us first recall some elementary notion and facts about
$K$-theory of small Waldhausen categories (see \cite{wald, tt,
Quillen}).

Let $\mathcal{C}$ be a small Waldhausen category, that is, a pointed
category (equipped with a zero object) with cofibrations and weak
equivalences. In \cite[Section 1.3]{wald}, Waldhausen has defined a
$K$-theory space $K(\mathcal{C})$ for $\mathcal{C}$, which is a
pointed topological space, and an $n$-th homotopy group
$K_n(\mathcal C)$ of $K(\mathcal C)$ for each $n\in\mathbb{N}$,
which is called the $n$-th $K$-group of $\mathcal{C}$. Clearly, if a
Waldhausen category $\mathcal{C}'$ is \emph{essentially small}, that
is, the isomorphism classes of objects of $\mathcal{C}'$ form a set,
then the definition of Waldhausen $K$-theory still makes sense for
$\mathcal{C}'$ because, in this case, one can choose a small
Waldhausen subcategory $\mathcal{C}$ of $\mathcal{C}'$ such that
$\mathcal{C}$ is equivalent to $\mathcal{C}'$, and define the
$K$-theory of $\mathcal{C}'$ through that of $\mathcal{C}$.

Note that $K(\mathcal C)$ is always homotopy equivalent to a
CW-complex. In fact, this follows from the following observation:
The classifying space of a small category has the structure of a
connected CW-complex and the loop space of a CW-complex is homotopy equivalent
to a CW-complex (see \cite{milnor2}), while $K(\mathcal{C})$ is the
loop space of the classifying space constructed from $\mathcal{C}$.

The $K$-theory space defined by Waldhausen is natural in the
following sense: Each exact functor $F:\mathcal{C}\to \mathcal{D}$
between Waldhausen categories $\mathcal{C}$ and $\mathcal{D}$
induces a continuous map $K(F): K(\mathcal C)\to K(\mathcal D)$ of
(pointed) topological spaces, and a homomorphism $K_n(F):
K_n(\mathcal{C})\to K_n(\mathcal{D})$ of abelian groups for each
$n\in\mathbb{N}$. If $G:\mathcal{D}\ra \mathcal{E}$ is another exact
functor between Waldhausen categories, then $K(GF)=K(F)K(G)$ in our
notation.

Note that the associated point $e_{\mathcal C}$ of $K(\mathcal C)$
corresponds to the image of the map $K(\{0\})\to K(\mathcal C)$
induced from the inclusion $\{0\}\hookrightarrow\mathcal{C}$, where
$0$ denotes the zero object of $\mathcal{C}$.

Finally, we recall some definitions and basic facts in homotopy theory for later proofs. For more details, we refer the reader to \cite[Chapters III and IV]{Wh} and \cite[Chapter 7]{se}.  Those
readers who are familiar with homotopy theory may skip the rest of this subsection.

Let $g: Y\to Z$ be a continuous map of topological spaces. We say that $g$ is a \emph{homotopy equivalence} if there is a continuous map $h:Z\to Y$ such that $gh: Y\to Y$ and $hg:Z\to Z$ are homotopic to the identities of $Y$ and $Z$, respectively. If there is a homotopy equivalence between $Y$ and $Z$, then we say that $Y$ and $Z$ are \emph{homotopy equivalent}, and simply write
$Y\lraf{\sim} Z$.

Assume that $Y$ and $Z$ are pointed topological spaces with the base-points $y_0$ and $z_0$,
respectively, and that the map $g:Y\to Z$ sends $y_0$ to $z_0$.  The {\emph {homotopy fibre}} $F(g)$ of $g$  is defined to be the following pointed topological space
$$
F(g):=\{(\omega, y)\mid \omega:[0,1]\to Z,\,y\in
Y,\,(0)\omega=z_0,\, (1)\omega=(y)g\}
$$
with the base-point $\big(c_{z_0},y_0\big)$, where $c_{z_0}$ is the
constant path $t\mapsto z_0$ for $t\in [0,1]$. Note that homotopy fibres are well defined up to homotopy equivalences.

The homotopy fibre of the map $\{z_0\}\hookrightarrow Z$ is called the \emph{loop space} of $(Z,z_0)$, and denoted by $\Omega(Z,z_0)$. Note that we can identify  $\Omega(Z,z_0)$ with the set $\{\omega:[0,1]\to Z\mid(0)\omega=z_0=(1)\omega\}$, and that there is a canonical map
$$\partial:\;\Omega(Z,z_0)\lra F(g),\; \omega\mapsto (\omega, y_0)\;\;\mbox{for}\;\; \omega\in \Omega(Z,z_0).$$
Let $\pi_n(Z,z_0)$ denote the \emph{$n$-th homotopy group} of $(Z,z_0)$ for each $n\in\mathbb{N}$. Then $\pi_n\big(\Omega(Z,z_0)\big)=\pi_{n+1}(Z,z_0)$.

Further, we define $h:F(g)\to Y$ by $(\omega, y)\mapsto y$ for any $(\omega, y)\in F(g)$. Then the sequence $$\Omega(Z,z_0)\lraf{\partial} F(g)\lraf{h} (Y,y_0)\lraf{g} (Z,z_0) $$
gives rise to a long exact sequence of homotopy groups:
$$\cdots \lra \pi_{n+1}(Z,z_0)\lraf{\pi_n(\partial)} \pi_n\big(F(g),(c_{z_0},y_0)\big)\lraf{\pi_n(h)}
\pi_n(Y,y_0)\lraf{\pi_n(g)} \pi_n(Z,z_0)\lra
\pi_{n-1}\big(F(g),(c_{z_0},y_0)\big)\lra$$
$$\cdots\lra \pi_0\big(F(g),(c_{z_0},y_0)\big)\lra \pi_0(Y,y_0)\lra
\pi_0(Z,z_0).$$ For a proof, we refer the reader to
\cite[Corollary IV. 8.9]{Wh}.

A sequence $(X,x_0)\lraf{f} (Y,y_0)\lraf{g} (Z,z_0)$ of pointed
topological spaces is called \emph{a homotopy fibration} if the
composite of $f$ and $g$ is equal to the constant map which sends
every $x$ in $X$ to the base-point of $Z$, and if the natural map
$$X\lra F(g), \quad \; x\mapsto \big(c_{\,z_0}, (x)f\big)\;\;\mbox{for}\;\; x\in X$$ is a homotopy equivalence. In this case, the loop space $\Omega(Z,z_0)$ is homotopy equivalent to the homotopy fibre of $f$.

The sequence $(X,x_0)\lraf{f} (Y,y_0)\lraf{g} (Z,z_0)$ of pointed
topological spaces is called \emph{a weak homotopy fibration} if
there is a pointed topological space $(Z', z_0')$, and two pointed
maps $g_1: Y\to Z'$ and $g_2:Z'\to Z$ with $g=g_1 g_2$ such that

$(1)$ the sequence $(X,x_0)\lraf{f} (Y,y_0)\lraf{g_1} (Z',z_0')$ is
a homotopy fibration, and that

$(2)$ $g_2$ induces an injection $\pi_0(Z', z_0')\to \pi_0(Z, z_0)$
and a bijection $\pi_n(Z',z_0')\to \pi_n(Z, z_0)$ for $n>0$.

Assume that $(X,x_0)\lraf{f} (Y,y_0)\lraf{g} (Z,z_0)$ is a weak
homotopy fibration. Then there is a long exact sequence of homotopy
groups:
$$\cdots \lra \pi_{n+1}(Z,z_0)\lra \pi_n(X,x_0)\lraf{\pi_n(f)}
\pi_n(Y,y_0)\lraf{\pi_n(g)} \pi_n(Z,z_0)\lra \pi_{n-1}(X,x_0)\lra$$
$$\cdots\lra \pi_0(X,x_0)\lra \pi_0(Y,y_0)\lra
\pi_0(Z,z_0)$$ for all $n\in \mathbb{N},$ and $g_2$ induces a homotopy
equivalence $\Omega(g_2):\Omega(Z',z_0')\lraf{\sim} \Omega(Z,z_0)$.
Thus $\Omega(Z,z_0)$ is homotopy equivalent to the homotopy
fibre of $f$.

\subsection{Frobenius pairs and their $K$-theory spaces\label{fpairs}}
We recall some definitions given in \cite{Sch1}.

By a \emph{Frobenius category} we mean an exact category (see
\cite{Quillen, keller2}) with enough projective and injective
objects such that projectives and injectives coincide. A map between
two Frobenius categories is an exact functor which preserves
projective objects.

Let $ \mathcal{C}$ be a Frobenius category.

We denote by $\mathcal{C}$-{\rm proj} the full subcategory of
$\mathcal{C}$ consisting of all projective objects. It is well known
that the factor category $\underline{\mathcal C}$ of $\mathcal{C}$
modulo $\mathcal{C}$-proj, called the \emph{stable category} of
$\mathcal{C}$, is a triangulated category. Moreover, two objects $X$
and $Y$ of $\mathcal{C}$ are isomorphic in $\underline{\mathcal C}$
if and only if $X\oplus P\simeq Y\oplus Q$ in $\mathcal{C}$ for some
$P, Q\in\mathcal{C}\mbox{-proj}.$ In particular, $X\simeq 0$ in
$\underline{\mathcal C}$ if and only if $X\in
\mathcal{C}\mbox{-proj}$.

A subcategory $\mathcal{X}$ of $\mathcal{C}$ is called a
\emph{Frobenius subcategory} of $\mathcal{C}$  if $\mathcal{X}$ is a
Frobenius category and the inclusion $\mathcal{X}\subseteq
\mathcal{C}$ is a fully faithful map of Frobenius categories. In
this case, $\mathcal{X}\mbox{-proj}\subseteq
\mathcal{C}\mbox{-proj}$, and a morphism in $\mathcal{X}$ factorizes
through $\mathcal{X}\mbox{-proj}$ if and only if it factorizes
through $\mathcal{C}\mbox{-proj}$. This implies that the inclusion
$\mathcal{X}\subseteq \mathcal{C}$ induces a fully faithful
inclusion $\underline{\mathcal{X}}\subseteq \underline{\mathcal{C}}$
of triangulated categories. In general, $\underline{\mathcal X}$
does not have to be a triangulated subcategory of
$\underline{\mathcal{C}}$ since $\underline{\mathcal X}$ is not
necessarily closed under isomorphisms in $\underline{\mathcal{C}}$.
However, by our convention, the image of the inclusion
$\underline{\mathcal{X}}\subseteq \underline{\mathcal{C}}$ is indeed
a triangulated subcategory of $\underline{\mathcal{C}}$.

A pair $\mathbf{C}:=(\mathcal{C}, \mathcal{C}_0)$ of Frobenius
categories is called a \emph{Frobenius pair} if $\mathcal{C}$ is a
small category and $\mathcal{C}_0$ is a Frobenius subcategory of
$\mathcal{C}$. A map from a Frobenius pair $(\mathcal{C},
\mathcal{C}_0)$ to another Frobenius pair $(\mathcal{C}',
\mathcal{C}_0')$ is a map of Frobenius categories
$\mathcal{C}\to\mathcal{C}'$ such that it restricts to a map from
$\mathcal{C}_0$ to $\mathcal{C}_0'$ (see \cite[Section 4.3]{Sch1}).

Let $\mathbf{C}:=(\mathcal{C}, \mathcal{C}_0)$ be a Frobenius pair.
Then the image of the inclusion $\underline{{\mathcal
C}_0}\subseteq\underline{\mathcal C}$ is a triangulated subcategory
of $\underline{\mathcal{C}}$. So we can form the Verdier quotient of
$\underline{\mathcal C}$ by this image, denoted by
$$\DF{\bf \mathbf C}:=\underline{\mathcal C}\,/\underline{{\mathcal C}_0}$$
which is called the $\emph{derived category}$ of the Frobenius pair
$\mathbf{C}$. Here, we use the same notation $\underline{\mathcal
C}\,/\underline{{\mathcal C}_0}$ as in \cite{Sch1} to denote the
derived category of $\mathbf{C}$, but the meaning of $
\underline{\mathcal C}\,/\underline{{\mathcal C}_0}$ in our paper is
slightly different from the one in \cite{Sch1} because we require
that the image of an inclusion functor is closed under isomorphisms.
Nevertheless, all results in \cite{Sch1} work with this modified
definition of derived categories.

Clearly, if $\mathcal{C}_0=\mathcal{C}\mbox{-\rm proj}$, then
$\DF{\mathbf{C}}=\underline{\mathcal C}$. In this case, we shall
often write $\mathcal{C}$ for the Frobenius pair $(\mathcal{C},
\mathcal{C}\mbox{-\rm proj})$.

The category $\mathcal{C}$ of a Frobenius pair
$\mathbf{C}:=(\mathcal{C}, \mathcal{C}_0)$ can be regarded as a
small Waldhausen category (for definition, see \cite{wald} or
\cite{xc4}): The inflations in $\mathcal{C}$ form the cofibrations
of $\mathcal{C}$, and the morphisms in $\mathcal{C}$ which are
isomorphisms in $\DF{\mathbf C}$ form the weak equivalences of
$\mathcal{C}$. In this note, we shall write $\mathbf{C}$ for the
Waldhausen category $\mathcal{C}$ to emphasize the role of
$\mathcal{C}_0$. According to our foregoing notation, we denote by
$\mathcal{C}$ the Waldhausen category defined by the Frobenius pair
$(\mathcal{C}, \mathcal{C}\mbox{-\rm proj})$. For the Waldhausen
category $\mathbf{C}$, we denote  the $K$-theory space of
$\mathbf{C}$ in the sense of Waldhausen by $K(\mathbf C)$ which is a
pointed topological space, and  the $n$-th $K$-group of $K(\mathbf
C)$ by $K_n(\mathbf C)$ for each $n\in\mathbb{N}$.

It is known that $K_0(\mathbf C)$ is naturally isomorphic to the
Grothendieck group $K_0(\DF{\mathbf C})$ of the small triangulated
category $\DF{\mathbf C}$ (see \cite[Section 1.5.6]{tt},
\cite[Chapter IV, Proposition 8.4]{weibelbook} and \cite[Proposition
3.2.22]{Sch2}).

Let $G:\mathbf C\to\mathbf {C\,'} $ be a map of Frobenius pairs. On
the one hand, $G$ automatically induces a triangle functor
$\DF{G}:\DF{\mathbf C}\to\DF{\mathbf {C\,'}}$, which sends
$X\in\mathcal{C}$ to $G(X)\in\mathcal{C}'$. On the other hand,
$G:\mathcal{C}\to \mathcal{C\,'}$ is an exact functor of associated
Waldhausen categories, which induces a continuous map $K(G):
K(\mathbf C)\to K(\mathbf{C\,'})$.

In this paper, we assume that all Waldhausen categories considered
arise from Frobenious pairs. Two typical examples of Frobenius pairs are of our particular interest.

(a) The first typical example of Frobenius pairs is provided by
 the categories of bounded complexes over exact categories.

Let $\mathscr {E}$ be a small exact category (for definition, see
\cite{Quillen} and \cite{keller2}). We denote by $\Cb{\mathscr E}$
the category of bounded chain complexes over $\mathscr E$. Then
$\Cb{\mathscr E}$ is a small, exact category with degreewise split
conflations, that is, a sequence $\cpx{X}\to\cpx{Y}\to \cpx{Z}$ is a
conflation in $\Cb{\mathscr E}$ if $X^i\to Y^i\to Z^i$ is isomorphic
to the split conflation $X^i\to X^i\oplus Z^i\to Z^i$ for each
$i\in\mathbb{Z}$. Actually, $\Cb{\mathscr E}$ is even a Frobenius
category in which projective objects are exactly bounded
contractible chain complexes over $\mathscr{E}$. Recall that a chain
complex $\cpx{X}$ is called \emph{contractible} when the identity on
$\cpx{X}$ is null-homotopic. Moreover, the stable category of
$\Cb{\mathscr E}$ is the usual bounded homotopy category
$\Kb{\mathscr E}$, that is,
$\DF{\mathscr{C}^b(\mathscr{E})}=\Kb{\mathscr{E}}$.

Recall that a complex $\cpx{X}=(X^i,d^i)_{i\in \mathbb{Z}}$ over
$\mathscr{E}$ is called \emph{acyclic} if $d^i$ is a composite of a
deflation $\pi^i$ with an inflation $\lambda^i$ such that
$(\lambda^i,\pi^{i+1})$ is a conflation for all $i$. Let
$\mathscr{C}^b_{ac}{(\mathscr E)}\subseteq \Cb{\mathscr E}$ be the
full subcategory of objects which are homotopy equivalent to acyclic
chain complexes over $\mathscr{E}$. Then
$\mathscr{C}^b_{ac}{(\mathscr E)}$ contains all projective objects
of the Frobenius category $\Cb{\mathscr{E}}$, and is closed under
extensions, kernels of deflations as well as cokernels of inflations
in $\Cb{\mathscr{E}}$. Thus $\mathscr{C}^b_{ac}{(\mathscr E)}$
inherits a Frobenius structure from $\Cb{\mathscr E}$ and
$$\mathbf{C}:=\big(\Cb{\mathscr E}, \mathscr{C}^b_{ac}{(\mathscr E)}\big)$$
is a Frobenius pair. In particular, the pair $\mathbf{C}$ (or the
associated category $\Cb{\mathscr E}$) can be regarded as a
Waldhausen category: A chain map $\cpx{f}:\cpx{X}\to \cpx{Y}$ in
$\Cb{\mathscr E}$ is called a cofibration if $f^i: X^i\to Y^i$ is  a
split inflation in $\mathscr{E}$ for each $i\in\mathbb{Z};$ a weak
equivalence if the mapping cone of $\cpx{f}$ belongs to
$\mathscr{C}^b_{ac}{(\mathscr E)}.$ Moreover, $\DF{\mathbf{C}}$
coincides with the bounded derived category $\Db{\mathscr E}$ of
$\mathscr{C}^b(\mathscr{E})$, which is defined as follows:

Let $\mathscr E'$ be an arbitrary exact category. The objects of
$\Db{\mathscr E'}$ are the objects of $\mathscr{C}^b(\mathscr{E}')$.
The morphisms of $\Db{\mathscr E'}$ are obtained from the chain maps
by formally inverting the maps whose mapping cones are acyclic (as
complexes of objects in $\mathscr{E}'$). For example, if
$\mathscr{E}'$ is the usual exact category $R\Modcat$ with $R$ a
ring, then $\Db{\mathscr E'}$ is the usual derived category
$\Db{R}$. For more details, see \cite{keller2}.

Assume that the exact structure of $\mathscr{E}$ is induced from an
abelian category $\mathscr{A}$. That is, $\mathscr{E}\subseteq
\mathscr{A}$ is a full subcategory such that it is closed under
extensions, and that a sequence $X\lraf{f} Y\lraf{g} Z$ with all
terms in $\mathscr{E}$ is a conflation in $\mathscr{E}$ if and only
if $0\to X\lraf{f} Y\lraf{g} Z\to 0$ is an exact sequence in
$\mathscr{A}$. Furthermore, assume that $\mathscr{E}$ is closed
under kernels of epimorphisms in the abelian category. In this case,
the chain map $\cpx{f}:\cpx{X}\to \cpx{Y}$ is a weak equivalence in
$\mathbf{C}$ if and only if $\cpx{f}$ is a quasi-isomorphism in
$\mathscr{C}(\mathscr{A})$, that is, $H^i(\cpx{f}): H^i(\cpx{X})\to
H^i(\cpx{Y})$ is an isomorphism in $\mathscr{A}$ for each
$i\in\mathbb{Z}$.

Note that an exact category $\mathscr{E}$ itself can also be
understood as a Waldhausen category with cofibrations being
inflations, and weak equivalences being isomorphisms. Up to now,
there are at least three algebraic $K$-theory spaces associated with
a small exact category $\mathscr E$: The Quillen $K$-theory space of
the exact category $\mathscr{E}$, the Waldhausen $K$-theory space
with respect to the Waldhausen category $\mathscr{E}$, and the
Waldhausen $K$-theory space of the Waldhausen category defined by
the Frobenius pair $\big(\Cb{\mathscr
E},\mathscr{C}^b_{ac}({\mathscr E})\big)$. However, these spaces are
the same up to homotopy equivalence (see \cite[Section 1.9]{wald})
and \cite[Theorem 1.11.7]{tt}). So, in this paper, we always
identify these spaces.

(b) The next example of Frobenius pairs is constructed from
categories of finitely generated projective modules.

Let $R$ be a ring. Then the category $\pmodcat{R}$ of finitely
generated projective $R$-modules is a small exact category with
split, short exact sequences as its conflations. Clearly, this exact
structure on $\pmodcat R$ is induced from the usual exact structure
of the abelian category $R\Modcat$. Following Quillen
\cite{Quillen}, the \emph{algebraic $K$-theory space} $K(R)$ of $R$
is defined to be the space $K(\pmodcat R)$ of $\pmodcat R$, and the
\emph{$n$-th algebraic $K$-group} $K_n(R)$ of $R$ to be the $n$-th
homotopy group of $K(R)$.

\medskip
We know from $(a)$ that the pair $\big(\Cb{\pmodcat R}, \mathscr{C}^b_{ac}{(\pmodcat R})\big)$ is a Frobenius pair. In this way, $\Cb{\pmodcat R}$ can be regarded as a small Waldhausen
category. Moreover, $\mathscr{C}^b_{ac}{(\pmodcat R)}$ consists of all bounded contractible chain complexes over $\pmodcat R$, which are exactly projective objects in the Frobenius category
$\Cb{\pmodcat R}$. In other words, we have $\mathscr{C}^b(\pmodcat R)\mbox{-proj} =
\mathscr{C}^b_{ac}{(\pmodcat R)}$. Thus $\DF{\Cb{\pmodcat R}}$ is
the bounded homotopy category $\Kb{\pmodcat R}$. Since each compact
object of $\D{R}$ is quasi-isomorphic to an object of $\Cb{\pmodcat
R}$, the Verdier localization functor $\K{R}\to\D{R}$ restricts to a triangle equivalence $\Kb{\pmodcat R}\lraf{\simeq} \mathscr{D}^c(R)$.

Hence, we see that $K(R)$, $K(\Cb{\pmodcat R})$ and $K(\mathbf{C})$
with $\mathbf{C}:=\big(\Cb{\pmodcat R}, \mathscr{C}^b_{ac}{(\pmodcat
R)}\big)$ are homotopy equivalent, and therefore their algebraic
$K_n$-groups are all isomorphic.

Let $S$ be another ring and $\cpx{N}$ a bounded complex of $S$-$R$-bimodules. If ${_S}\cpx{N}\in\Cb{\pmodcat S}$, then the tensor functor $\cpx{N}\cpx{\otimes}_R-:\Cb{\pmodcat R}\to\Cb{\pmodcat S}$ is a map of Frobenious pairs. So, we obtain a map $K(\cpx{N}\cpx{\otimes}_R-):K(R)\to K(S)$ of $K$-theory spaces.

In case $\lambda: R\to S$ is a ring homomorphism, we choose $\cpx{N}=S$ and denote simply by $K(\lambda)$ the map $K(S\otimes_R-):K(R)\to K(S)$. Since the homomorphism  $K_n(S[1]\cpx{\otimes}_R-):K_n(R)\to K_n(S)$, induced from the map $K(S[1]\cpx{\otimes}_R-):K(R)\to K(S)$, is equal to the minus of $K_n(\lambda)$, we shall denote the map $K(S[1]\cpx{\otimes}_R-)$ by $-K(\lambda)$.

Note that the shift functor $[1]:\Cb{\pmodcat{R}}\to \Cb{\pmodcat{R}}$ is also a map of Frobenius pairs. Now, let $\Delta$ be the diagonal map $x\mapsto (x,x)$ for
$x\in K(R)$ and let $\sqcup:\pmodcat{R}\times \pmodcat{R}\to \pmodcat{R}$ be
the coproduct functor. Then the induced map $K([1]): K(R)\to K(R)$ is a homotopy equivalence and a homotopy inverse of $K(R)$ in the sense that the composite of the following
maps:
$$\xymatrix{
K(R)\ar[r]^-{\Delta} & K(R)\times K(R)\ar[rr]^-{K([1])\times
Id} && K(R)\times K(R)\ar[r]^-{K(\sqcup)} & K(R)}
$$
is homotopic to the constant map which sends $x$ to the base-point of $K(R)$.

\subsection{Fundamental theorems in algebraic $K$-theory of Frobenius pairs}

Now, we recall some basic results on algebraic $K$-theory of
Frobenious pairs in terms of derived categories. Our main reference
in this section is the paper \cite{Sch1} by Schlichting.

The following localization theorem may trace back to the
localization theorem in \cite[Section 5, Theorem 5]{Quillen} for
exact categories, the fibration theorem in \cite[Theorem
1.6.4]{wald} for Waldhausen categories, and the localization theorem
in \cite[Theorem 1.8.2]{tt} for complicial biWaldhausen categories.
For a proof of the present form, we refer the reader to
\cite[Propositions 3 and 5, p.126 and p.128]{Sch1}.  Also, the
approximation and cofinality theorems are taken from
\cite[Propositions 3 and 4]{Sch1}.

\begin{Lem}\label{loc}

$(1)$ \emph{Localization Theorem:}

Let $\mathbf{A} \lraf{F} \mathbf{B} \lraf{G} \mathbf{C}$ be a
sequence of Frobenius pairs. If the sequence $\DF{\mathbf A}
\lraf{\DF{F}}\DF{\mathbf B} \lraf{\DF{G}}\DF{\mathbf C}$ of derived
categories is exact, then the induced sequence $K(\mathbf {A})
\lraf{K(F)} K(\mathbf{B}) \lraf{K(G)} K(\mathbf{C})$ of $K$-theory
spaces is a homotopy fibration, and therefore there is a long exact
sequence of $K$-groups $$\cdots \lra K_{n+1}(\mathbf{C})\lra
K_n(\mathbf{A})\lraf{K_n(F)} K_n(\mathbf{B})\lraf{K_n(G)}
K_n(\mathbf{C})\lra K_{n-1}(\mathbf{A})\lra$$
$$\cdots\lra K_0(\mathbf{A})\lra K_0(\mathbf{B})\lra
K_0(\mathbf{C})\lra 0$$ for all $n\in \mathbb{N}.$

$(2)$ \emph{Approximation Theorem:}

Let $G:\mathbf{B}\to \mathbf C$ be a map of Frobenius pairs. If the
associated functor $\DF{G}: \DF{\mathbf B}\to \DF{\mathbf C}$ of
derived categories is an equivalence, then the induced map $K(G):
K(\mathbf {B})\to K(\mathbf{C})$ of $K$-theory spaces is a homotopy
equivalence. In particular, $K_n(G): K_n(\mathbf{B}) \lraf{\simeq}
K_n(\mathbf{C})$ for all $n\in \mathbb{N}.$

$(3)$ \emph{Cofinality Theorem:}

Let $G:\mathbf{B} \to \mathbf{C}$ be a map of Frobenius pairs. If
the associated functor $\DF{G}: \DF{\mathbf{B}}\to \DF{\mathbf{C}}$
of derived categories is an equivalence up to factors, then the
induced map $K(G): K(\mathbf{B})\to K(\mathbf{C})$ of $K$-theory
spaces gives rise to an injection $K_0(G): K_0(\mathbf{B})\ra
K_0(\mathbf{C})$ and an isomorphism: $K_n(G):
K_n(\mathbf{B})\lraf{\simeq} K_n(\mathbf{C})$ for all $n>0$.
\end{Lem}

Note that the surjectivity of the last map in the long exact
sequence in Lemma \ref{loc} (1) follows from the fact that
$K_0(\mathbf{C})$ is isomorphic to the Grothendieck group
$K_0\big(\DF{\mathbf{C}}\big)$ of $\DF{\mathbf{C}}$.

The localization theorem is useful, but when we deal with $K$-theory of recollements, the obstacle for us to use it is that, in a given recollement of derived module categories, we do not know whether the given functors between derived categories are induced from exact functors between Frobenius pairs.

For our purpose of later proofs, we mention the following result which is a slight variation of \cite[Section
6.1]{Sch1} and has been mentioned there without proof. For the
convenience of the reader, we include here a proof (see also
\cite[Lemma 2.5]{nr} for a special case).

\begin{Lem} \emph{Thickness Theorem:}\label{thickness}

Let $\mathbf{C}:=(\mathcal{C}, \mathcal{C}_0)$ be a Frobenius pair.
Suppose that there is a triangulated category $\mathscr{C}$ together
with a triangle equivalence $G:\DF{\mathbf C}\to \mathscr{C}$. Let
$\mathscr{X}$ be a full triangulated subcategory of $\mathscr{C}$ .
Define $\mathcal{X}$ to be the full subcategory of $\mathcal{C}$
consisting of objects $X$ such that $G(X)\in\mathscr{X}$. Then the
following statements are true:

$(1)$ The category $\mathcal{X}$ contains $\mathcal{C}_0$ and is
closed under extensions in $\mathcal{C}$. Moreover, $\mathcal{X}$
naturally inherits a Frobenius structure from $\mathcal{C}$, and
becomes a Frobenius subcategory of $\mathcal{C}$ such that
$\mathcal{X}$-{\rm proj} = $\mathcal{C}$-{\rm proj}.

$(2)$ Both $\mathbf{X}:=(\mathcal{X}, \mathcal{C}_0)$ and
$\mathbf{C}_{\mathscr{X}}:=(\mathcal{C}, \mathcal{X})$ are Frobenius
pairs, and the inclusion functor $\mathcal{X}\to\mathcal{C}$ and the
identity functor $\mathcal{C}\to\mathcal{C}$ induce the following
commutative diagram of triangulated categories:
$$
\xymatrix{\DF{\mathbf {X}}\ar@{^{(}->}[r]\ar[d]^-{\simeq}
&\DF{\mathbf
C}\ar[r]\ar[d]^-{\simeq}_-{G} &\DF{\mathbf{C}_{\mathscr{X}}}\ar[d]^-{\simeq}\\
\mathscr{X}\ar@{^{(}->}[r] & \mathscr{C}\ar[r]&
\mathscr{C}/\mathscr{X}}
$$

$(3)$ If $\mathscr{X}$ is closed under direct summands in
$\mathscr{C}$, then both rows in the diagram of $(2)$ are exact
sequences of triangulated categories.
\end{Lem}

{\it Proof.} $(1)$ By definition of
$\DF{\mathbf{C}}:=\underline{\mathcal C}\,/\underline{{\mathcal
C}_0}$, the objects of $\DF{\mathbf{C}}$ are the same as the objects
of $\mathcal C$. Thus, if $M\in\mathcal{C}_0$ or
$M\in\mathcal{C}\mbox{-proj}$, then $M\simeq 0$ in
$\DF{\mathbf{C}}$. This implies that $\mathcal{X}$ contains both
$\mathcal{C}_0$ and $\mathcal{C}\mbox{-proj}$. Since $G$ is a
triangle functor and $\mathscr{X}$ is a full triangulated
subcategory of $\mathscr{C}$, it is easy to see that $\mathcal{X}$
is closed under extensions in $\mathcal{C}$.

Since $\mathcal{X}$ is closed under extensions in $\mathcal{C}$, we
can endow $\mathcal{X}$ with an  exact structure induced from the
one of $\mathcal{C}$, namely, a sequence $X\to Y\to Z$ with all
terms in $\mathcal{X}$ is called a conflation in $\mathcal{X}$ if it
is a conflation in $\mathcal{C}$. Then one can check that, with this
exact structure, $\mathcal{X}$ becomes an exact category. Now, we
claim that $\mathcal{X}$ is even a Frobenius category such that
$\mathcal{X}$-{\rm proj} = $\mathcal{C}$-{\rm proj}. Indeed, it
suffices to show that if $L\to P\to N$ is a conflation in
$\mathcal{C}$ with $P\in\mathcal{C}\mbox{-proj}$, then
$L\in\mathcal{X}$ if and only if $N\in\mathcal{X}$. Actually, such a
conflation can be extended to a distinguished triangle $L\to P\to
N\to L[1]$ in $\underline{\mathcal C}$, and further, to a
distinguished triangle in $\DF{\mathbf{C}}$. Since $P\simeq 0$ in
$\DF{\mathbf{C}}$, we have $N\simeq L[1]$ in $\DF{\mathbf{C}}$. As
$\mathscr{X}$ is closed under shifts in $\mathscr{C}$ and $G$ is a
triangle functor, we know that $G(L)\in\mathscr{X}$ if and only if
$G(N)\in\mathscr{X}$. In other words, $L\in\mathcal{X}$ if and only
if $N\in\mathcal{X}$. This verifies the claim.

$(2)$ Note that $\mathcal{C}_0\subseteq
\mathcal{X}\subseteq\mathcal{C}$ and
$\mathcal{C}_0\mbox{-proj}\subseteq
\mathcal{X}\mbox{-proj}=\mathcal{C}\mbox{-proj}.$ Thus
$\mathbf{X}:=(\mathcal{X},\mathcal{C}_0)$ and
$\mathbf{C}_{\mathscr{X}}:=(\mathcal{C},\mathcal{X}) $ are Frobenius
pairs.

Recall that $\DF{\mathbf{X}}:=\underline{\mathcal
X}\,/\underline{{\mathcal C}_0}$ and
$\DF{\mathbf{C}_\mathscr{X}}:=\underline{\mathcal
C}\,/\underline{{\mathcal X}}.$ Clearly, the inclusion functor
$\lambda: \mathcal{X}\to\mathcal{C}$ and the identity functor
$Id_\mathcal{C}: \mathcal{C}\to\mathcal{C}$ are maps from the
Frobenius pair $\mathbf{X}$ to the Frobenius pairs $\mathbf{C}$, and
from $\mathbf{C}$ to $\mathbf{C}_{\mathscr{X}}$, respectively. So we
have two triangle functors $\DF{\lambda}:\underline{\mathcal
X}\,/\underline{{\mathcal C}_0}\to \underline{\mathcal
C}\,/\underline{{\mathcal C}_0}$ and
$\DF{Id_\mathcal{C}}:\underline{\mathcal C}\,/\underline{{\mathcal
C}_0}\to \underline{\mathcal C}\,/\underline{\mathcal X}$, which are
induced from the inclusion $\underline{\mathcal
X}\subseteq\underline{\mathcal{C}}$ and the identity functor of
$\underline{\mathcal{C}}$, respectively.

Clearly, $\underline{\mathcal{X}}$ contains
$\underline{\mathcal{C}_0}$, that is, the objects of
$\underline{\mathcal{C}_0}$ is a subclass of the objects of
$\underline{\mathcal{X}}$ with the morphism set
$\Hom_{\underline{\mathcal{C}_0}}(X,Y)=\Hom_{\underline{\mathcal{X}}}(X,Y)$
for all objects $X,Y$ in $\mathcal{C}_0$. Since the inclusion
$\underline{\mathcal X}\subseteq\underline{\mathcal{C}}$ is fully
faithful, the functor $\DF{\lambda}$ is also a fully faithful
inclusion which gives rise to the following commutative diagram:
$$
(\ast)\quad \xymatrix{\underline{\mathcal X}\,/\underline{{\mathcal
C}_0}\ar@{^{(}->}[r]^-{\DF{\lambda}} \ar[d]^-{\simeq} &
\underline{\mathcal C}\,/\underline{{\mathcal
C}_0}\ar[d]^-{\simeq}_-{G} \\
\mathscr{X}\ar@{^{(}->}[r] & \mathscr{C}.}
$$
Consequently, $G$ induces a triangle equivalence
$$G_1: (\underline{\mathcal
C}\,/\underline{{\mathcal C}_0})/(\underline{\mathcal
X}\,/\underline{{\mathcal
C}_0})\lraf{\simeq}\mathscr{C}/\mathscr{X}.$$

By the universal property of the Verdier localization functor
$q_1:\underline{\mathcal C}\to \underline{\mathcal
C}\,/\underline{\mathcal X}$ (respectively,
$q_2:\underline{\mathcal{C}}/\underline{\mathcal{C}}_0\to
(\underline{\mathcal C}\,/\underline{{\mathcal
C}_0})/(\underline{\mathcal X}\,/\underline{{\mathcal C}_0})$),
there is a triangle functor $\phi: \underline{\mathcal
C}\,/\underline{{\mathcal X}}\ra (\underline{\mathcal
C}\,/\underline{{\mathcal C}_0})/(\underline{\mathcal
X}\,/\underline{{\mathcal C}_0})$ (respectively, $\psi:
(\underline{\mathcal C}\,/\underline{{\mathcal
C}_0})/(\underline{\mathcal X}\,/\underline{{\mathcal
C}_0})\to\underline{\mathcal C}\,/\underline{{\mathcal X}}\,$) such
that $q_2q_0 = \phi q_1$ (respectively, $\DF{Id_\mathcal{C}}=\psi
q_2$), where $q_0:\underline{\mathcal C}\ra
\underline{\mathcal{C}}/\underline{\mathcal{C}}_0$ is the Verdier
localization functor. Since $q_1=\DF{Id_\mathcal{C}}q_0$, we have
$$\psi\phi q_1=\psi
q_2q_0=\DF{Id_\mathcal{C}}q_0=q_1\mbox{\; and \;}
\phi\DF{Id_\mathcal{C}} q_0=\phi q_1=q_2q_0.$$ It follows that
$\psi\phi=Id$ and $\phi\DF{Id_\mathcal{C}}=q_2$. As $\phi\psi
q_2=\phi\psi\phi\DF{Id_\mathcal{C}}=\phi\DF{Id_{\mathcal C}}=q_2$,
we obtain $\phi\psi=Id$. Thus $\phi$ is a triangle isomorphism.

Now, we define $\ol{G}:=G_1\phi:\underline{\mathcal
C}\,/\underline{{\mathcal X}}\to \mathscr{C}/\mathscr{X}$. Then the
following diagram of triangulated categories
$$ (\ast\ast)\quad
\xymatrix{\underline{\mathcal C}\,/\underline{{\mathcal
C}_0}\ar[r]^-{\DF{Id_\mathcal{C}}}\ar[d]^-{\simeq}_-{G}
&\underline{\mathcal C}\,/\underline{\mathcal X}\ar[d]^-{\simeq}_-{\ol{G}}\\
\mathscr{C}\ar[r]^-{q}& \mathscr{C}/\mathscr{X}.}
$$
is commutative, where $q$ is the Verdier localization functor. Now,
$(2)$ follows from ($\ast$) and ($\ast\ast$).

$(3)$ In this case, $\mathscr{X}$ is the kernel of the localization
functor $q:\mathscr{C}\to \mathscr{C}/\mathscr{X}$. Thus $(3)$
follows. $\square$.

\section{Algebraic $K$-theory of differential graded algebras \label{3.5}}

\subsection{Definitions of $K$-theory spaces of dg algebras}

In this subsection, we shall give a definition of $K$-theory spaces
of differential graded algebras, which generalizes the one of
$K$-theory spaces of usual rings and modifies slightly the
definition in \cite{Sch1}. But, at the level of homotopy groups, the
two definitions give the isomorphic algebraic $K_n$-groups for $n\in
\mathbb N$.

Throughout this subsection, $k$ stands for an arbitrary but fixed
commutative ring, and all rings considered here are $k$-algebras. Note that each
ring with identity can be regarded as a
$\mathbb{Z}$-algebra.

Let $\mathbb{A}$ be a differential graded (dg) associative and
unitary $k$-algebra, that is, $\mathbb{A}=\oplus_{n\in
\mathbb{Z}}A^n$ is a $\mathbb{Z}$-graded $k$-algebra with a
differential $d^n:A^n\ra A^{n+1}$ such that $(A^n, d^n)_{n\in\mathbb
Z}$ is a chain complex of $k$-modules and
$$(xy)d^{m+n}=x(yd^n)+(-1)^n (xd^m)y$$ for $m\in\mathbb{Z}$, $x\in
A^m$ and $y\in A^n$. Thus the map
$\mathbb{A}\cpx{\otimes}_k\mathbb{A}\ra \mathbb{A}$,
$a\otimes_kb\mapsto ba$ for $a,b\in \mathbb{A}$, is a chain map.

A left dg $\mathbb{A}$-module $\cpx{M}$ is a $\mathbb{Z}$-graded
left module $\cpx{M}=\oplus_{n\in \mathbb{Z}}M^n$ over the
$\mathbb{Z}$-graded $k$-algebra $\mathbb{A}$, with a differential
$d$ such that $(M^n,d)_{n\in{\mathbb{Z}}}$ is a complex of
$k$-modules, and for any $a\in A^m, x\in M^n$, the following holds:
$$(ax)d^{m+n}=a(xd^n)+(-1)^n (ad^m)x.$$
In particular, each dg $\mathbb A$-module is a $\mathbb{Z}$-graded
$\mathbb{A}$-module (forgetting the differential).

For a dg $\mathbb{A}$-module $\cpx{M}$, we denote by $\cpx{M}[1]$
the shift of $\cpx{M}$ by degree $1$.

We should observe that the dg algebra $(\mathbb{A},d)$ and left dg
$\mathbb{A}$-module $\cpx{M}$ defined in this paper are actually the
dg algebra $(\mathbb{A}^{\opp},d)$ and right dg
$\mathbb{A}^{\opp}$-module in the sense of \cite[Summary]{keller1},
respectively.

In the following, we give a typical way to obtain dg algebras by taking Hom-complexes of dg modules.

Let $(\cpx{X},d_{\cpx{X}})$ and $(\cpx{Y},d_{\cpx{Y}})$ be two dg $\mathbb{A}$-modules. The \emph{Hom-complex} of $\cpx{X}$ and $\cpx{Y}$ over $\mathbb{A}$ is defined to be the following complex $\cpx{\Hom}_{\mathbb{A}}(\cpx{X},\cpx{Y}):=\big(\Hom_\mathbb{A}^n(\cpx{X},\cpx{Y}),
d^{\,n}_{\cpx{X},\cpx{Y}}\big)_{n\in\mathbb{Z}}$ over $k$:

As a $k$-module, the $n$-th component $\Hom_\mathbb{A}^n(\cpx{X},\cpx{Y})$ is formed by the morphisms $h:\cpx{X}\to \cpx{Y}$ of graded $\mathbb{A}$-modules, homogeneous of degree $n$. In other words, $h$ is a homomorphism of $\mathbb{A}$-modules such that $h=(h^p)_{p\in\mathbb{Z}}$ with $h^p\in\Hom_k(X^p,Y^{p+n})$. Further, the differential $d^{\,n}_{\cpx{X},\cpx{Y}}:\Hom_\mathbb{A}^n(\cpx{X},\cpx{Y})\to \Hom_\mathbb{A}^{n+1}(\cpx{X},\cpx{Y})$ of degree $n$ is given by
$$(h^p)_{p\in\mathbb{Z}}\mapsto \big(h^p
d_{\cpx{Y}}^{p+n}-(-1)^{n}d_{\cpx{X}}^p
h^{p+1}\big)_{p\in\mathbb{Z}}.$$

\smallskip
Furthermore, we take another dg $\mathbb{A}$-module $\cpx{Z}$, and define
$$
\circ: \;\cpx{\Hom}_\mathbb{A}(\cpx{X},\cpx{Y})\times
\cpx{\Hom}_\mathbb{A}(\cpx{Y},\cpx{Z})\lra
\cpx{\Hom}_\mathbb{A}(\cpx{X},\cpx{Z}),\;\;(f, g)\mapsto
(f^pg^{p+m})_{p\in\mathbb Z}
$$
for $f:=(f^p)_{p\in\mathbb{Z}}\in\Hom_\mathbb{A}^m(\cpx{X},\cpx{Y})$ and
$g:=(g^p)_{p\in\mathbb{Z}}\in\Hom_\mathbb{A}^n(\cpx{Y},\cpx{Z})$ with $m,
n\in\mathbb{Z}$. Thus the operation $\circ$ is associative and
distributive. In particular, under this operation, $\cpx{\Hom}_\mathbb{A}(\cpx{X},\cpx{X})$ is a
$\mathbb Z$-graded ring. Moreover, the above-defined operation $\circ$ satisfies the
following identity:
$$
(f\circ g)\,d^{\,m+n}_{\cpx{X},\cpx{Z}}=f\circ
(g)d^{\,n}_{\cpx{Y},\cpx{Z}}+
(-1)^n(f)d^{\,m}_{\cpx{X},\cpx{Y}}\circ g.
$$
This implies that $\cpx{\Hom}_\mathbb{A}(\cpx{X},\cpx{X})$, together with the differential of
itself as a complex over $k$, is a dg algebra. In this sense, $\cpx{\End}_\mathbb{A}(\cpx{X})$ will be called the \emph{dg endomorphism ring} of $\cpx{X}$, and denoted simply by $\cpx{\End}_\mathbb{A}(\cpx{X})$. Also, due to the above identity,
the complex $\cpx{\Hom}_{\mathbb{A}}(\cpx{X},\cpx{Y})$ is actually a left dg
$\cpx{\End}_\mathbb{A}(\cpx{X})$- and right dg $\cpx{\End}_\mathbb{A}(\cpx{Y})$-
bimodule.

Now, we recall the definition of the category $\C{\mathbb A}$of left dg $\mathbb A$-modules.
Actually, this category has all dg $\mathbb{A}$-modules as objects, and
a homomorphism $\cpx{f}: \cpx{X}\to \cpx{Y}$ between dg $\mathbb{A}$-modules
$\cpx{X}$ and $\cpx{Y}$ is a chain map of complexes over $k$, which commutes with the $\mathbb{A}$-actions on $\cpx{X}$ and $\cpx{Y}$. This means that $\Hom_{\C{\mathbb A}}(\cpx{X},\cpx{Y})$ is exactly the $0$-th cocycle of the complex $\cpx{\Hom}_{\mathbb{A}}(\cpx{X},\cpx{Y})$. It is known that $\C{\mathbb A}$ is a Frobenius category (see \cite[Section 2]{keller1}) by declaring a conflation to be a
short sequence of dg $\mathbb{A}$-modules such that the underlying
sequence of graded $\mathbb{A}$-modules (forgetting differentials)
is split exact. The stable category of $\C{\mathbb A}$ is the dg
homotopy category $\K{\mathbb A}$ in which the objects are the dg
$\mathbb A$-modules and the morphisms are the homotopy classes of
homomorphisms of dg $\mathbb{A}$-modules. In other words, $\Hom_{\K{\mathbb A}}(\cpx{X},\cpx{Y})$
is equal to the $0$-th cohomology $H^0\big(\cpx{\Hom}_{\mathbb{A}}(\cpx{X},\cpx{Y})\big)$ of the complex $\cpx{\Hom}_{\mathbb{A}}(\cpx{X},\cpx{Y})$.

We say that $\cpx{f}$ is a \emph{quasi-isomorphism} if it is a quasi-isomorphism as a chain map of complexes over $k$, that is, $H^i(\cpx{f}): H^i(\cpx{X})\to H^i(\cpx{Y})$ is an isomorphism for
every $i\in\mathbb{Z}$. By inverting all quasi-isomorphisms of dg $\mathbb{A}$-modules, we obtain the \emph{dg derived category} $\D{\mathbb A}$ of $\mathbb A$. This is a
triangulated category and generated by the dg module $\mathbb A$,
that is, $\D{\mathbb A}=\Tria(\mathbb A).$

Observe that an ordinary $k$-algebra $A$ can be regarded as a dg algebra concentrated in degree $0$, and that the above-mentioned
categories $\C{A}$, $\K{A}$ and $\D{A}$ coincide with the usual
complex, homotopy and derived categories of the category of left $A$-modules,
respectively.

To give a description of $\D{\mathbb{A}}$ by a triangulated subcategory of $\K{\mathbb A}$ up to equivalence, we shall recall some more definitions in \cite{keller1}.

The dg $\mathbb{A}$-module $\cpx{X}$ is said to be \emph{acyclic} if it is acyclic as a complex of $k$-modules, that is, $H^i(\cpx{X})=0$ for all $i\in\mathbb{Z}$; is said to have the \emph{property $(P)$} if $\Hom_{\K{\mathbb A}}(\cpx{X}, \cpx{Y})=0$ for any acyclic dg
$\mathbb A$-module $\cpx{Y}$, or equivalently, $\cpx{\Hom}_{\mathbb{A}}(\cpx{X},\cpx{Y})$ is acyclic as a complex over $k$. Note that the class of dg $\mathbb
A$-modules with the property $(P)$ is closed under extensions,
shifts, direct summands and direct sums in $\C{\mathbb A}$. We
denote by $\K{\mathbb A}_p$ the full subcategory of $\K{\mathbb A}$
consisting of all modules with the property $(P)$. Then $\K{\mathbb
A}_p\subseteq \K{\mathbb A}$ is a triangulated subcategory
containing $\mathbb A$ and being closed under direct sums. More
important, by \cite[Section 3.1]{keller1}, the Verdier localization functor
$q:\K{\mathbb A}\to\D{\mathbb A}$ restricts to a triangle equivalence
$$\widetilde{q}:\,\K{\mathbb A}_p\lraf{\simeq} \D{\mathbb A}.$$ This implies that any
quasi-isomorphism between two dg $\mathbb{A}$-modules with the
property $(P)$ is an isomorphism in $\K{\mathbb A}$ and that, for
each dg $\mathbb{A}$-module $\cpx{M}$, there is a (functorial)
quasi-isomorphism ${_p}\cpx{M}\to \cpx{M}$ of dg
$\mathbb{A}$-modules such that ${_p}\cpx{M}$ has the property $(P)$.

With the help of the above triangle equivalence, we can define the total left-derived functors of
tensor functors. This procedure is similar to the one for usual complexes over ordinary rings.

Let $\cpx{W}$ be a right dg $\mathbb{A}$-module and $\cpx{X}$ a (left) dg $\mathbb{A}$-module. The \emph{tensor complex}  of $\cpx{W}$ and $\cpx{X}$ over $\mathbb{A}$ is defined to be the following complex $\cpx{W}\cpx{\otimes}_\mathbb{A}\cpx{X}:=\big(\cpx{W}{\otimes}^n_\mathbb{A}\cpx{X},
\partial^n_{\cpx{W},\cpx{X}}\big)_{n\in\mathbb{Z}}$ over $k$:

As a $k$-module, the $n$-th component $\cpx{W}{\otimes}^n_\mathbb{A}\cpx{X}$ is the quotient module of $\bigoplus_{p\in\mathbb{Z}}W^p\otimes_kX^{n-p}$ modulo the
$k$-submodule generated by all elements $ua\otimes v-u\otimes av$ for $u\in W^r$,
$a\in\mathbb{A}^s$ and $v\in X^t$ with $r, s, t\in\mathbb{Z}$ and $n=r+s+t$. Further, the differential $\partial_{\cpx{W},\cpx{X}}$ of degree $n$ is
given by
$$w\otimes x \mapsto(w)d^p_{\cpx{W}}\otimes x
+(-1)^p w\otimes(x)d_{\cpx{X}}^{n-p}$$ for  $w\in W^p $ and $x\in
X^{n-p}$.

Assume further that $\mathbb{B}$ is another dg algebra and that $\cpx{W}$ is a $\mathbb{B}$-$\mathbb{A}$-bimodule. Then $\cpx{W}\cpx{\otimes}_\mathbb{A}\cpx{X}$ is indeed a
dg $\mathbb{B}$-module. This gives rise to the following tensor
functor
$$
{_\mathbb{B}}\cpx{W}\cpx{\otimes}_\mathbb{A}-:\C{\mathbb A}\lra \C{\mathbb
B},\; \cpx{X}\mapsto \cpx{W}\cpx{\otimes}_\mathbb{A}\cpx{X}.
$$
Now, the total left-derived functor
$\cpx{W}\otimesL_\mathbb{A}-: \D{\mathbb A}\to\D{\mathbb B}$ of this
functor is defined by $\cpx{X}\mapsto
\cpx{W}\cpx{\otimes}_\mathbb{A}{(_p\cpx{X})}$. Particularly, if $\cpx{X}$ has the property $(P)$, then $\cpx{W}\otimesL_\mathbb{A}\cpx{X}=\cpx{W}\cpx{\otimes}_\mathbb{A}\cpx{X}$
in $\D{\mathbb B}$.

A dg $\mathbb A$-module $M$ is called \emph{relatively countable
projective} (respectively, \emph{countable projective}) if there is
a dg $\mathbb A$-module $N$ such that $M \oplus N$ is isomorphic to
$\bigoplus_{ i\in I} \mathbb A[n_i]$ as dg $\mathbb A$-modules
(respectively, as $\mathbb{Z}$-graded $\mathbb A$-modules), where
$I$ is a countable set and $n_i \in \mathbb{Z}$. Clearly,
relatively countable projective modules are countable projective
modules, and have the property $(P)$ since
$\Hom_{\K{\mathbb{A}}}(\mathbb{A}[i],M)\simeq H^{-i}(M)$ for all
$i$.

Let $\mathcal{X}(\mathbb A)$ be the full subcategory of $\C{\mathbb
A}$ consisting of countable projective $\mathbb A$-modules. Then
$\mathcal{X}(\mathbb A)$ is an essentially small category. This is
due to the following observation: Let $\mathcal{G}(\mathbb A)$ be
the category of $\mathbb{Z}$-graded $\mathbb{A}$-modules. For every
$X:=\bigoplus_{i\in\mathbb Z}X^i\in\mathcal{G}(\mathbb A)$, we have
the following: $(a)$ The class $\mathcal{U}(X)$ consisting of
isomorphism classes of direct summands of $X$ in
$\mathcal{G}(\mathbb{A})$ is a set. In fact, there is a surjection
from the set of idempotent elements of $\End_{\mathcal{G}(\mathbb
A)}(X)$ to $\mathcal{U}(X)$. $(b)$ The class $\mathcal{V}(X)$
consisting of all dg $\mathbb{A}$-modules with $X$ as the underlying
graded $\mathbb{A}$-module is also a set since $\mathcal{V}(X)$ is
contained into the set $\{(X, d^i)_{i\in\mathbb{Z}}\mid
d^i\in\Hom_k(X^i, X^{i+1})\}$, which is a countable union of sets.

Furthermore, $\mathcal{X}(\mathbb A)$ is closed under extensions,
shifts, direct summands and countable direct sums in $\C{\mathbb
A}$.

Let $\C{\mathbb A, \aleph_0}$ be the smallest full subcategory of
$\mathcal{X}(\mathbb A)$ such that it

$(1)$ contains all relatively countable projective
$\mathbb{A}$-modules;

$(2)$ is closed under extensions and shifts;

$(3)$ is closed under countable direct sums.

\noindent Then $\C{\mathbb A, \aleph_0}$ is essentially small,
inherits an exact structure from $\C{\mathbb A}$, and becomes a
fully exact subcategory of $\C{\mathbb A}$. Even more, $\C{\mathbb
A, \aleph_0}$ is a Frobenius subcategory of $\C{\mathbb A}$, in
which projective-injective objects are the ones of $\C{\mathbb A}$
belonging to $\C{\mathbb A, \aleph_0}.$ This can be concluded from
the following fact: For each $M\in\C{\mathbb A}$, there is a
canonical conflation $M\to C(M)\to M[1]$ in $\C{\mathbb A}$ such
that $C(M)$ is a projective-injective object of $\C{\mathbb A}$ (see
\cite[Section 2.2]{keller1}). Hence $\C{\mathbb A, \aleph_0}$
provides a natural Frobenius pair $(\C{\mathbb A, \aleph_0},
\C{\mathbb A, \aleph_0}\mbox{-proj})$, and the inclusion $\C{\mathbb
A, \aleph_0}\subseteq \C{\mathbb A}$ induces a fully faithful
inclusion from the derived category $\DF{\C{\mathbb A, \aleph_0}}$
of $\C{\mathbb A, \aleph_0}$ to $\K{\mathbb A}$.

We denote by $\K{\mathbb A, \aleph_0}$ the full subcategory of
$\K{\mathbb A}$ consisting of those complexes which are isomorphic
in $\K{\mathbb A}$ to objects of $\C{\mathbb A, \aleph_0}$. Then
$\K{\mathbb A, \aleph_0}$ is a triangulated subcategory of
$\K{\mathbb A}$ by the condition $(2)$, and the inclusion
$\DF{\C{\mathbb A, \aleph_0}}\subseteq\K{\mathbb A, \aleph_0}$ is a
triangle equivalence. Since the full subcategory of
$\mathcal{X}(\mathbb A)$ consisting of all dg $\mathbb{A}$-modules
with the property $(P)$ satisfies the above conditions $(1)$-$(3)$,
we deduce that each object of $\C{\mathbb A, \aleph_0}$ has the
property $(P)$. This implies that $\K{\mathbb A,
\aleph_0}\subseteq\K{\mathbb A}_p$. Furthermore, by definition,
$\C{\mathbb A, \aleph_0}$ is closed under countable direct sums in
$\C{\mathbb A}$, and therefore $\K{\mathbb A, \aleph_0}$ is closed
under countable direct sums in $\K{\mathbb A}_p$. It follows from
Lemma \ref{literature} (1) that $\K{\mathbb A, \aleph_0}$ is closed
under direct summands in $\K{\mathbb A}_p$.

Now, let $\mathscr{X}(\mathbb A)$ be the full subcategory of
$\D{\mathbb A}$ consisting of all those objects which are isomorphic
in $\D{\mathbb A}$ to the images of objects of $\K{\mathbb A,
\aleph_0}$ under the equivalence $\widetilde{q}: \K{\mathbb
A}_p\lraf{\simeq} \D{\mathbb A}$. Then $\mathscr{X}(\mathbb A)$ is a
triangulated subcategory of $\D{\mathbb A}$ closed under direct
summands, and $\widetilde{q}$ induces a triangle equivalence from
$\K{\mathbb A, \aleph_0}$ to $\mathscr{X}(\mathbb A)$. In all, we
have
$$
\DF{\C{\mathbb A, \aleph_0}}\subseteq\K{\mathbb A,
\aleph_0}\subseteq \K{\mathbb A}_p,\;\; \mathscr{X}(\mathbb
A)\subseteq \D{\mathbb{A}}
$$
and
$$\DF{\C{\mathbb A, \aleph_0}}\lraf{\simeq}\K{\mathbb A,
\aleph_0}\lraf{\simeq}\mathscr{X}(\mathbb A)$$ as triangulated
categories.

Recall that a dg $\mathbb{A}$-module $M$ is called a \emph{finite
cell module} if there is a finite filtration $$ 0=M_0\subseteq
M_1\subseteq M_2\subseteq\cdots \subseteq M_n=M
$$
of dg $\mathbb A$-modules such that, for each $0\leq i\leq
n-1\in\mathbb{N}$, the quotient module $M_{i+1}/M_i$ is isomorphic
to $\mathbb{A}[n_i]$ for some $n_i\in\mathbb{Z}$ (see \cite[Part
III]{km}). Clearly, each finite cell $\mathbb{A}$-module belongs to
$\C{\mathbb A, \aleph_0}$. Moreover, the category of finite cell
$\mathbb{A}$-modules is closed under extensions in $\C{\mathbb A,
\aleph_0}$. Actually, this category is a Frobenius subcategory of
$\C{\mathbb A, \aleph_0}$, in which projective-injective objects are
the ones of $\C{\mathbb A, \aleph_0}$ belonging to this subcategory.

An object $M\in \D{\mathbb A}$ is said to be \emph{compact} if
$\Hom_{\D{\mathbb A}}(M,-)$ commutes with direct sums in $\D{\mathbb
A}$. Let $\mathscr{D}^c (\mathbb A)$ be the full subcategory of
$\D{\mathbb A}$ consisting of all compact objects. Then
$\mathscr{D}^c (\mathbb A)$ is the smallest full triangulated
subcategory of $\D{\mathbb A}$ containing $\mathbb A$ and being
closed under direct summands of its objects. In fact, each compact
object of $\D{\mathbb A}$ is a direct summand of a finite cell
module in $\D{\mathbb A}$ (see \cite[Section 5]{keller1}). This
implies the following chain of full subcategories: $\mathscr{D}^c
(\mathbb A)\subseteq \mathscr{X}(\mathbb A)\subseteq \D{\mathbb A}.$

Now, we define $\mathcal{W}_{\mathbb A}$ to be the full subcategory
of $\C{\mathbb A, \aleph_0}$ consisting of all those objects in
$\mathscr{C}(\mathbb{A},\aleph_0)$ such that they are isomorphic in
$\D{\mathbb A}$ to compact objects of $\D{\mathbb A}$. Clearly,
$\mathcal{W}_{\mathbb A}$ is essentially small. Moreover, by
applying Lemma \ref{thickness} to the Frobenius pair $\C{\mathbb A,
\aleph_0}$ and the equivalence $\DF{\mathscr{C}(\mathbb{ A},
\aleph_0)}\lraf{\simeq} \mathscr{X}(\mathbb A)$ with the
triangulated subcategory $\mathscr{D}^c(\mathbb{A})$ of
$\mathscr{X}(\mathbb A)$, we deduce that $\mathcal{W}_{\mathbb A}$
is a Frobenius subcategory of $\C{\mathbb A, \aleph_0}$ with the
same projective objects, and that the following diagram of
triangulated categories commutes:
$$(\star)\quad \xymatrix{\DF{\mathcal{W}_{\mathbb A}}\ar@{^{(}->}[r] \ar[d]^-{\simeq} &
\DF{\mathscr{C}(\mathbb{A},
\aleph_0)}\ar[d]^-{\simeq}\ar@{^{(}->}[r] &\K{\mathbb
A}_p\ar[d]^-{\simeq}_-{\widetilde{q}} \ar@{^{(}->}[r]& \K{\mathbb
A}\ar[dl]_-{q}\\
\mathscr{D}^c(\mathbb A)\ar@{^{(}->}[r] &\mathscr{X}(\mathbb
A)\ar@{^{(}->}[r] &\D{\mathbb A} &}
$$

From now on, we regard $\mathcal{W}_{\mathbb A}$ as a Waldhausen
category in the sense of Subsection \ref{fpairs}, namely, it arises
exactly from the Frobenius pair $(\mathcal{W}_{\mathbb A},
\mathcal{W}_{\mathbb A}\mbox{-proj})$.

\begin{Def}
The algebraic $K$-theory space of the dg $k$-algebra $\mathbb A$ is defined
to be the space $K(\mathcal{W}_{\mathbb A})$ of the Waldenhausen category $\mathcal{W}_{\mathbb A}$, denoted by $K(\mathbb{A})$.
For each $n\in\mathbb{N}$, the $n$-th algebraic $K$-group of $\mathbb A$ is defined to be the $n$-th homotopy group of $K(\mathbb{A})$ and denoted by $K_n(\mathbb A)$.
\end{Def}

Note that $K_0(\mathbb A)$ is isomorphic to $K_0(\DF{\mathcal{W}_{\mathbb
A}})$, the Grothendieck group of the (essentially small)
triangulated category $\DF{\mathcal{W}_{\mathbb{A}}}$ of the
Frobenius pair $(\mathcal{W}_{\mathbb A}, \mathcal{W}_{\mathbb
A}\mbox{-proj})$(see Subsection \ref{fpairs}).
As a result, we have the following fact.

\begin{Lem} \label{dgk}
The Verdier localization functor $\K{\mathbb A}\to \D{\mathbb A}$
induces a triangle equivalence: $\DF{\mathcal{W}_{\mathbb
A}}\lraf{\simeq} \mathscr{D}^c (\mathbb A)$. In particular,
$K_0(\mathcal{W}_{\mathbb A})$ is isomorphic to the Grothendieck
group $K_0(\mathscr{D}^c(\mathbb A))$ of $\mathscr{D}^c(\mathbb A)$.
\end{Lem}

Our definition of $K$-theory spaces of dg algebras has the following property.

\begin{Lem}\label{fin-cell}
Let $\mathcal{F}_{\mathbb A}$ be the full subcategory of
$\mathcal{W}_{\mathbb A}$ consisting of all finite cell
$\mathbb{A}$-modules. Then the inclusion $\mathcal{F}_{\mathbb A}\to
\mathcal{W}_{\mathbb A}$ induces an injection
$K_0(\mathcal{F}_{\mathbb A})\to K_0(\mathcal{W}_{\mathbb A})$ and
an isomorphism $K_n(\mathcal{F}_{\mathbb A})\lraf{\sim}
K_n(\mathcal{W}_{\mathbb A})$ for each $n>0$.
\end{Lem}

{\it Proof.} Note that $\mathcal{F}_{\mathbb A}$ is a Frobenius
subcategory of $\mathcal{W}_\mathbb{A}$ and that the inclusions
$\mathcal{F}_\mathbb{A}\subseteq\mathcal{W}_\mathbb{A}\subseteq\mathscr{C}(\mathbb{A})$
induce fully faithful inclusions $\DF{\mathcal{F}_{\mathbb
A}}\subseteq \DF{\mathcal{W}_\mathbb{A}}\subseteq \K {\mathbb A}_p$
(see Subsection \ref{fpairs}).

To show that the inclusion $\DF{\mathcal{F}_{\mathbb A}}\to
\DF{\mathcal{W}_{\mathbb A}}$ is an equivalence up to factors, we
shall compare the images of these two categories under the
equivalence $\tilde{q}: \K{\mathbb A}_p\to\D{\mathbb A}$ in the
above diagram ($\star$). In fact, by Lemma \ref{dgk}, the
restriction of the functor $\tilde{q}$  to $\DF{{\mathcal W}_\mathbb
A}$ gives rise to a triangle equivalence $\DF{{\mathcal W}_\mathbb
A} \lraf{\simeq} \mathscr{D}^c (\mathbb A)$. Let $\mathscr{Y}$ be
the smallest full triangulated subcategory of $\Dc{\mathbb A}$
containing $\mathbb A$. Since the objects of $\DF{{\mathcal
F}_\mathbb A}$ are the same as the ones of $\mathcal{F}_\mathbb{A}$,
the image of the restriction of the functor $\tilde{q}$ to
$\DF{{\mathcal F}_\mathbb A}$ is contained in $\mathscr{Y}$, and
therefore is equal to $\mathscr{Y}$. Thus $\tilde{q}$ induces a
triangle equivalence $\DF{\mathcal{F}_\mathbb{A}}\lraf{\simeq}
\mathscr{Y}$. Since $\Dc{\mathbb A}=\tria(\mathbb{A})$ and
$\mathbb{A}\in\mathscr{Y}\subseteq\Dc{\mathbb{A}}$, we have
$\tria(\mathscr{Y})=\Dc{\mathbb{A}}$. So the inclusion
$\mathscr{Y}\to\Dc{\mathbb A}$ is an equivalence up to factors.
Consequently, the inclusion $\DF{\mathcal{F}_{\mathbb A}}\to
\DF{\mathcal{W}_{\mathbb A}}$ induced from $\mathcal{F}_{\mathbb
A}\subseteq \mathcal{W}_{\mathbb A}$ is also an equivalence up to
factors. Now, Lemma \ref{fin-cell} follows from Lemma \ref{loc} (3).
$\square$.

\begin{Rem}\label{remark0}
In \cite[Section 12.3]{Sch1}, a $K$-theory spectrum
$\mathbb{K}(\mathcal{F}_\mathbb A)$ is defined for the category
$\mathcal{F}_\mathbb{A}$. Moreover, it is known in \cite[Theorem
8]{Sch1} that, for each $n\in\mathbb{N}$, the $n$-th homology group
of $\mathbb{K}(\mathcal{F}_\mathbb A)$ is given by
$$
\pi_n\big(\mathbb{K}(\mathcal{F}_\mathbb A)\big)=\left\{\begin{array}{ll} K_n(\mathcal{F}_\mathbb{A}) & \mbox{if } n>0,\\
K_0(\mathscr{D}^c(\mathbb A)) & \mbox{if } n=0.\end{array}\right.
$$ Thus Lemmas \ref{fin-cell} and \ref{dgk} show that
$\pi_n\big(\mathbb{K}(\mathcal{F}_\mathbb A)\big)\simeq K_n(\mathbb
A)$ for all $n\in\mathbb{N}$, and therefore, at the level of
homotopy groups, our definition of $K$-theory for dg algebras is
isomorphic to the one defined by Schlichting in \cite{Sch1}.
\end{Rem}

The following result, together with  Lemma \ref{fin-cell}, may
explain the advantage of defining $K$-theory of arbitrary dg
algebras by using the category $\mathcal{W}_\mathbb{A}$ rather than
$\mathcal{F}_\mathbb{A}$.

\begin{Lem}\label{same}
Let $A$ be an algebra with identity, and let $\mathbb{A}$ be the dg
algebra $A$ concentrated in degree $0$. Then $K(A)$ and
$K(\mathbb A)$ are homotopy equivalent as K-theory spaces.
\end{Lem}

{\it Proof.} Clearly, $\C{\mathbb A}=\C{A}$, $\K{\mathbb A}=\K{A}$
and $\D{\mathbb A}=\D{A}$. In particular, $\mathscr{D}^c(\mathbb
A)=\mathscr{D}^c(A)$. By the construction of $\mathcal{W}_{\mathbb
A}$, we see that $\Cb{\pmodcat A}\subseteq \mathcal{W}_{\mathbb A}$ and
$\mathscr{C}^b(\pmodcat A)\mbox{-proj}=
\mathscr{C}^b_{ac}(\pmodcat A)\subseteq\mathcal{W}_{\mathbb A}\mbox{-proj}$.
Thus the inclusion $j:\Cb{\pmodcat A}\to\mathcal{W}_{\mathbb A}$
is a fully faithful map of Frobenius pairs. In other words, $\Cb{\pmodcat A}$
is a Frobenius subcategory of $\mathcal{W}_{\mathbb A}$. This implies that
the triangle functor $\DF{j}:\DF{\Cb{\pmodcat A}}\to\DF{\mathcal{W}_{\mathbb A}}$ is fully faithful
(see Subsection \ref{fpairs}). Now we show that $\DF{j}$ is an equivalence.
On the one hand, the localization functor $q:\K{A}\to \D{A}$
induces an equivalence $q_1:\DF{\mathcal{W}_{\mathbb A}}\to
\mathscr{D}^c(A)$ by Lemma \ref{dgk}. On the other hand, the composite
of the following functors:
$$\Kb{\pmodcat R}=\DF{\Cb{\pmodcat
A}}\lraf{\DF{j}\;}\DF{\mathcal{W}_{\mathbb A}}\lraf{q_1}
\mathscr{D}^c(A)$$ is also an equivalence induced by $q$. Thus
$\DF{j}$ is a triangle equivalence. By Lemma \ref{loc} (2), we know that
$K(A)\lraf{\sim} K(\mathcal{W}_{\mathbb A})=:K(\mathbb A)$ as
$K$-theory spaces. $\square$

\medskip
\subsection{Homotopy equivalences of $K$-theory spaces from perfect dg modules}
In this subsection, we introduce the definition of perfect dg modules over dg rings, and discuss
homotopy equivalences of $K$-theory spaces of dg algebras linked by perfect dg modules.

Let $\mathbb{A}$ be a dg algebra. A dg $\mathbb{A}$-module is said to be
{\emph {perfect}} if it belongs to $\mathcal{W}_\mathbb{A}$. Recall that
each perfect dg $\mathbb{A}$-module has the property $(P)$ and is compact in
$\D{\mathbb{A}}$. Conversely, each compact dg $\mathbb{A}$-module is isomorphic in $\D{\mathbb A}$ to a perfect one, but itself may not have the property $(P)$. Moreover, if $\mathbb{A}$ is an ordinary ring concentrated in degree $0$, then each bounded complex of finitely generated projective $\mathbb{A}$-modules is perfect.

First of all, we point out the following result, which may illustrate the importance of perfect dg modules. For a proof, we refer to \cite[Section 3.1]{keller1}.

\begin{Lem}\label{dg-0}
Let $M$ be a dg $\mathbb{A}$-module and let $\mathbb{S}:=\cpx{\End}_{\mathbb{A}}(M)$. If ${_\mathbb{A}}M$ is perfect, then the left-derived functor $M\otimesL_{\mathbb{S}}-:\D{\mathbb{S}}\to\Tria({_\mathbb{A}}M)$
is a triangle equivalence.
\end{Lem}

In the following lemma, we can see that perfect dg modules always provide us with maps of Frobenius pairs which define algebraic $K$-theory spaces of dg algebras.

\begin{Lem}\label{functors}
Let $\mathbb{B}$ be a dg algebra and let $M$ be a dg $\mathbb{A}$-$\mathbb{B}$-bimodule. If $_{\mathbb{A}}M$ is perfect, then the tensor functor ${_\mathbb{A}}M\cpx{\otimes}_{\mathbb{B}}-:\mathcal{W}_\mathbb{B}\to\mathcal{W}_\mathbb{A}$
is a map of Frobenius pairs.
\end{Lem}

{\it Proof.} For simplicity, we denote by $G$ the tensor functor $M\cpx{\otimes}_{\mathbb B}-: \C{\mathbb{B}}\to \C{\mathbb{A}}$. In the following, we show that $G(\mathcal{W}_\mathbb{B})\subseteq \mathcal{W}_\mathbb{A}$.

Let $\mathcal{X}(\mathbb B)$ and $\mathcal{X}(\mathbb{A})$ be the full
subcategories of $\C{\mathbb B}$ and $\C{\mathbb{A}}$ consisting of
all countable projective modules, respectively. Recall that $\mathcal{X}(\mathbb{A})$ is
closed under shifts, direct summands and countable direct sums in $\C{\mathbb{A}}$. Then, it follows from $G(\mathbb{B})=M\cpx{\otimes}_{\mathbb B}\mathbb{B}\simeq M\in\mathcal{W}_\mathbb{A}\subseteq \mathcal{X}(\mathbb{A})$ that the functor $G:\mathcal{X}(\mathbb B)\to \mathcal{X}(\mathbb A)$ is well defined. Since $G$ always preserves conflations and commutes with both shifts and countable direct sums, the following full subcategory
$$G^{-1}(\C{\mathbb{A}, \aleph_0}):=\{N\in \mathcal{X}({\mathbb B})\mid
G(N)\in\C{A, \aleph_0} \}$$ of $\mathcal{X}(\mathbb{B})$ contains
all relatively countable projective $\mathbb{B}$-modules, and is
closed under extensions, shifts and countable direct sums.
Given that $\C{\mathbb B, \aleph_0}$ is the smallest subcategory of
$\mathcal{X}(\mathbb B)$ admitting these properties, we have
$\C{\mathbb B, \aleph_0}\subseteq G^{-1}(\C{\mathbb{A}, \aleph_0})$. Thus
$G(\C{\mathbb B, \aleph_0})\subseteq \C{\mathbb{A}, \aleph_0}$ and
$G:\C{\mathbb B, \aleph_0}\to\C{\mathbb{A}, \aleph_0}$ is a well-defined
functor.

Furthermore, since each object $N\in\C{\mathbb B, \aleph_0}$ always has the
property $(P)$, we see that $G(N)=M\otimesL_{\mathbb B}N$ in
$\D{\mathbb{A}}$. So, to show that $G(\mathcal{W}_\mathbb{B})\subseteq
\mathcal{W}_\mathbb{A}$, it suffices to prove that if $N\in \mathcal{W}_\mathbb{B}$, then $M\otimesL_{\mathbb B}N \in \Dc{\mathbb{A}}$. For checking this, we take an object $N\in \mathcal{W}_\mathbb{B}$. Then
$N\in\Dc{\mathbb{B}}$. Since each perfect dg $\mathbb{A}$-module is compact in $\D{\mathbb{A}}$, we have $M\otimesL_{\mathbb B}\mathbb{B}=M\otimes_{\mathbb B}\mathbb{B}\simeq M\in\Dc{\mathbb{A}}$. This implies that the functor
$M\otimesL_{\mathbb B}-:\D{\mathbb B}\to \D{\mathbb{A}}$ preserves
compact objects. Thus $M\otimesL_{\mathbb B}N \in\Dc{\mathbb{A}}$ and $G(\mathcal{W}_\mathbb{B})\subseteq \mathcal{W}_\mathbb{A}$.

Recall that, for an arbitrary dg algebra $\mathbb{S}$, the category $\mathcal{W}_\mathbb{S}\mbox{-proj}$ consists of all those objects which are homotopy equivalent to the zero object in $\C{\mathbb{S}}$. As $G$ always preserves conflations and homotopy equivalences, we see that $G$ sends projective objects of
$\mathcal{W}_\mathbb{B}$ to the ones of $\mathcal{W}_\mathbb{A}$. Thus $G:\mathcal{W}_\mathbb{B}\to \mathcal{W}_\mathbb{A}$ is a map of Frobenius pairs. $\square$

\medskip
Next, we show that perfect dg modules can offer homotopy equivalences of algebraic $K$-theory spaces.

\begin{Lem}\label{dg-1}
Let $\mathbb{B}$ be a dg algebra and let $M$ be a dg $\mathbb{A}$-$\mathbb{B}$-bimodule such that $_{\mathbb{A}}M$ is perfect. Let $\mathcal{P}$ be the full
subcategory of $\mathcal{W}_{\mathbb A}$ consisting of all those dg $\mathbb{A}$-modules,
which, regarded as objects of $\D{\mathbb{A}}$, belong to $\Tria{({_\mathbb{A}}M)}$. Then the followings hold true:

$(1)$ The category $\mathcal{P}$ is a Frobenius subcategory of $\mathcal{W}_\mathbb{A}$ and the map ${_\mathbb{A}}M\cpx{\otimes}_{\mathbb{B}}-:\mathcal{W}_\mathbb{B}\to\mathcal{W}_\mathbb{A}$  factorizes through the inclusion $\mathcal{P}\hookrightarrow \mathcal{W}_\mathbb{A}$.

$(2)$ If the left-derived functor
${_\mathbb{A}}M\otimesL_{\mathbb{B}}-:\D{\mathbb{B}}\to\Tria({_\mathbb{A}}M)$
is an equivalence, then ${_\mathbb{A}}M\cpx{\otimes}_{\mathbb{B}}-:\mathcal{W}_\mathbb{B}\to\mathcal{P}$ induces a homotopy equivalence $K(\mathbb{B})\lraf{\sim}K(\mathcal{P})$ of $K$-theory spaces.
If, in addition,  $\D{\mathbb{A}}=\Tria{({_\mathbb{A}}M)}$, then
$K(\mathbb{B})\lraf{\sim} K(\mathbb{A})$ as $K$-theory spaces.
\end{Lem}

{\it Proof.} $(1)$ Let  $\mathscr{X}:=\Tria{({_\mathbb{A}}M)}\cap {\mathscr{D}^c(\mathbb{A})}$.
Then $\mathscr{X}$ is a full triangulated subcategory of $\mathscr{D}^c(\mathbb{A})$. Since the
localization functor $q:\K{\mathbb{A}}\to \D{\mathbb{A}}$ induces a triangle equivalence $\DF{\mathcal{W}_{\mathbb A}}\lraf{\simeq} \mathscr{D}^c(\mathbb{A})$ by Lemma \ref{dgk}, we see that $\mathcal{P}$ is exactly the full subcategory of $\mathcal{W}_{\mathbb A}$ consisting of all those dg $\mathbb{A}$-modules, which are isomorphic in $\mathscr{D}^c(\mathbb{A})$ to objects of $\mathscr{X}$. Hence, by Lemma \ref{thickness}, $\mathcal{P}$ is a Frobenius subcategory
of $\mathcal{W}_{\mathbb A}$ and $q$ induces a triangle equivalence
$q_1:\DF{\mathcal{P}}\lraf{\simeq} \mathscr{X}$.

By Lemma \ref{functors}, the functor $G:\mathcal{W}_\mathbb{B}\to \mathcal{W}_\mathbb{A}$ is a map of Frobenius pairs. Note that $\D{\mathbb B}=\Tria(\mathbb{B})$ and $M\otimesL_{\mathbb B}-$ commutes with arbitrary direct sums. By Lemma \ref{literature} (2), we have $M\otimesL_{\mathbb B}N \in\Tria{({_\mathbb{A}}M)}$. It follows that $M\otimesL_{\mathbb B}N\in\Tria{({_\mathbb{A}}M)}\cap \Dc{\mathbb{A}}=\mathscr{X}$, and therefore $G(\mathcal{W}_\mathbb{B})\subseteq\mathcal{P}$. This implies that  $M\cpx{\otimes}_{\mathbb{B}}-:\mathcal{W}_\mathbb{B}\to\mathcal{W}_\mathbb{A}$ factorizes through the inclusion $\mathcal{P}\hookrightarrow \mathcal{W}_\mathbb{A}$.

$(2)$ Since $\D{\mathbb{A}}=\Tria(\mathbb{A})$ and ${_\mathbb{A}}M\in\Dc{\mathbb{A}}$, we know from \cite[Theorem 4.4.9]{ne3} that $\mathscr{X}$ coincides with the full subcategory of $\Tria({_\mathbb{A}}M)$ consisting of all compact objects in $\Tria({_\mathbb{A}}M)$. Now, suppose that the functor $M\otimesL_{\mathbb{B}}-:\D{\mathbb{B}}\to\Tria({_\mathbb{A}}M)$
is an equivalence. Then this functor restricts to a triangle equivalence $\Dc{\mathbb{B}}\lraf{\simeq}\mathscr{X}$. Moreover, by Lemma \ref{dgk}, the localization functor $\K{\mathbb B}\to\D{\mathbb B}$ induces an equivalence $\tilde{q}: \DF{\mathcal{W}_{\mathbb B}}\lraf{\simeq} \mathscr{D}^c(\mathbb B)$.
From the following commutative diagram:
$$
\xymatrix{
\DF{\mathcal{W}_\mathbb{B}}\ar[d]^-{\tilde{q}}_-{\simeq}\ar[rr]^-{\DF{M\cpx{\otimes}_{\mathbb{B}}-}}&& \DF{\mathcal{P}}\ar[d]_-{q}^-{\simeq}\\
\mathscr{D}^c({\mathbb B})\ar[rr]^-{M\otimesL_{\mathbb B}-}_-{\simeq}&& \mathscr{X}}
$$
we infer that $\DF{M\cpx{\otimes}_{\mathbb{B}}-}:\DF{\mathcal{W}_\mathbb{B}}\lraf{\simeq}
\DF{\mathcal{P}}$ is a triangle equivalence. It follows from Lemma \ref{loc} (2) that
$K(\mathbb B)\lraf{\sim} K(\mathcal{P})$ as $K$-theory spaces. Clearly, if  $\D{\mathbb{A}}=\Tria{({_\mathbb{A}}M)}$, then $\mathcal{P}=\mathcal{W}_\mathbb{A}$. Thus $(2)$ follows.  $\square$

\medskip
As a consequence of Lemma \ref{dg-1}, we re-obtain the following result in  \cite[Proposition 6.7 and Corollary 3.10]{DS} where its proof uses knowledge on model categories.

\begin{Koro}\label{kell}
Let $\lambda: \mathbb{B} \to \mathbb{A}$ be a homomorphism of dg
algebras which is a quasi-isomorphism. Then the functor
$\mathbb{A}\cpx{\otimes}_{\mathbb B}-: \C{\mathbb{B}}\ra
\C{\mathbb{A}}$ induces a homotopy equivalence $K(\mathbb
B)\lraf{\sim}K(\mathbb A)$ of $K$-theory spaces. In particular, if
$H^i(\mathbb{A})=0$ for all $i\neq 0$, then $K(\mathbb A)\lraf{\sim}
K(H^0(\mathbb A))$.
\end{Koro}

{\it Proof.} In Lemma \ref{dg-1}, we take $M=\mathbb{A}$. Then $M$ is a dg $\mathbb{A}$-$\mathbb{B}$-bimodule via $\lambda: \mathbb{B} \to \mathbb{A}$ such that it is perfect as a dg $\mathbb{A}$-module, and that $\Tria({_\mathbb{A}}M)=\D{\mathbb{A}}$. Since $\lambda$ is a quasi-isomorphism of dg algebras, it follows from \cite[Proposition 4.2]{km} that the functor $\mathbb{A}\cpx{\otimes}_{\mathbb B}-$ induces a triangle
equivalence $\D{\mathbb B}\lraf{\sim} \D{\mathbb A}$ (see also \cite[Section 3.1]{keller1}).
Now, the first part of Corollary \ref{kell} follows from Lemma \ref{dg-1} (2).

Suppose that $\mathbb{A}:=(A^i, d^i)_{i\in\mathbb{Z}}$ with $H^i(\mathbb{A})=0$ for all $i\neq 0$. We define $\tau^{\leq 0}(\mathbb A)$ to be the following dg algebra:
$$
\cdots \lra A^{-3}\lraf{d^{-3}} A^{-2}\lraf{d^{-2}}
A^{-1}\lraf{d^{-1}} \Ker(d^0)\lra 0\lra \cdots.
$$
Then there exist two canonical quasi-isomorphisms $\tau^{\leq 0}(\mathbb A)\to \mathbb A$ and $\tau^{\leq 0}(\mathbb A)\to
H^0(\mathbb A)$ of dg algebras. It follows from the first part of
Corollary \ref{kell} that $$K(\tau^{\leq 0}(\mathbb
A))\lraf{\sim}K(\mathbb A)\;\,\mbox{and}\;\, K(\tau^{\leq 0}(\mathbb
A))\lraf{\sim} K(H^0(\mathbb A)).$$ Combining these homotopy
equivalences with Lemma \ref{same}, we see that $K(\mathbb A)\lraf{\sim} K(H^0(\mathbb A))$ as $K$-theory spaces. $\square$

\medskip
Combining Lemma \ref{dg-1} with Lemma \ref{dg-0}, we have the following applicable result.

\begin{Koro}\label{dg-2}
Let $M$ be a perfect dg $\mathbb{A}$-module and let $\mathbb{S}:=\cpx{\End}_{\mathbb{A}}(M)$. Define $\mathcal{P}$ to be the full subcategory of $\mathcal{W}_{\mathbb A}$ consisting of all those dg $\mathbb{A}$-modules, which, regarded as objects in $\D{\mathbb{A}}$, belong to $\Tria{({_\mathbb{A}}M)}$. Then  $K(\mathbb{S})\lraf{\sim}K(\mathcal{P})$ as $K$-theory spaces.
Moreover, if $\D{\mathbb{A}}=\Tria{({_\mathbb{A}}M)}$, then
$K(\mathbb{S})\lraf{\sim} K(\mathbb{A})$ as $K$-theory spaces.

\end{Koro}

As a consequence of Corollary \ref{dg-2}, we obtain the following fact.

\begin{Koro}\label{dg-3}
Let $M$ and $N$ be two perfect dg $\mathbb{A}$-modules. If $\Tria(M)=\Tria(N)\subseteq \D{\mathbb{A}}$, then $$K\big(\cpx{\End}_{\mathbb{A}}(M)\big)\lraf{\sim}K\big(\cpx{\End}_{\mathbb{A}}(N)\big).$$
\end{Koro}

\smallskip
The following result conveys that, for ordinary rings, we can choose smaller subcategories of perfect complexes to realize the homtopy equivalence in Corollary \ref{dg-2}.

\begin{Koro}\label{Qe}
Let $A$ be an algebra and $\cpx{P}\in\Cb{\pmodcat A}$. Define
$\mathbb{S}:=\cpx{\End}_A(\cpx{P})$ and $\mathcal{P}$ to be the full
subcategory of $\Cb{\pmodcat A}$ consisting of all those complexes
which, regarded as objects in $\D{A}$, belong to $\Tria{(\cpx{P})}$.
Then $K(\mathbb S)\lraf{\sim} K(\mathcal P)$ as $K$-theory spaces.
\end{Koro}

{\it Proof.} We remark that $\mathcal{P}$ is a Frobenius subcategory of $\Cb{\pmodcat A}$ such that its derived category $\DF{\mathcal{P}}$ is equivalent to $\mathscr{X}:=\Tria{(\cpx{P})}\cap {\mathscr{D}^c(A)}$ via the localization functor $q:\K{A}\to \D{A}$.

Actually, since $\mathscr{X}$ is a full triangulated subcategory of
$\mathscr{D}^c(A)$ and $\DF{\Cb{\pmodcat A}}=\Kb{\pmodcat
A}\lraf{\simeq}\mathscr{D}^c(A)$, we see that $\mathcal{P}$ is
exactly the full subcategory of $\Cb{\pmodcat A}$, in which the
objects are complexes in $\Cb{\pmodcat A}$ such that they are
isomorphic in $\mathscr{D}^c(A)$ to objects of $\mathscr{X}$. Hence,
by Lemma \ref{thickness}, $\mathcal{P}$ is a Frobenius subcategory
of $\Cb{\pmodcat A}$ and the functor $q$ induces an equivalence
$q_1:\DF{\mathcal{P}}\lraf{\simeq} \mathscr{X}$.

Now we view $A$ as a dg algebra concentrated in degree $0$, and let $\mathcal{X}$ be the full subcategory of $\mathcal{W}_A$ consisting of those objects that are isomorphic in $\D{A}$ to objects of $\Tria{({_A}\cpx{P})}$. Since each object of $\mathcal{W}_A$ is compact in $\D{A}$, we clearly have $\mathcal{X}=\{X\in\mathcal{W}_A\mid X\in\mathscr{X}\}$. Note that $\cpx{P}$ is a dg $A$-$\mathbb{S}$-bimodule such that it is perfect as a dg $A$-module. So, from the proof of Lemma \ref{dg-1} (1),  we know that $\mathcal{X}$ is a Frobenius subcategory of $\mathcal{W}_A$, the functor $\cpx{P}\cpx{\otimes}_{\mathbb S}-:\mathcal{W}_\mathbb{S}\to \mathcal{X}$ is a map of Frobenius pairs and $q$ induces a triangle equivalence $q_2:\DF{\mathcal{X}}\lraf{\simeq} \mathscr{X}$.

In the following, we first show that
$K(\mathcal{P})\lraf{\sim}K(\mathcal{X})$, and then that $K(\mathbb
S)\lraf{\sim} K(\mathcal{X})$ as $K$-theory spaces. With these two
homotopy equivalences in mind, we will obviously have  $K(\mathbb
S)\lraf{\sim} K(\mathcal P)$, as desired.

Let us check that $K(\mathcal{P})\lraf{\sim}K(\mathcal{X})$.
Actually, it follows from $\Cb{\pmodcat A}\subseteq \mathcal{W}_A$
that $\mathcal{P}\subseteq \mathcal{X}$. Since
$\mathcal{P}\mbox{-proj} = \mathscr{C}^b_{ac}(\pmodcat A)\subseteq
\mathcal{W}_A\mbox{-proj}=\mathcal{X}\mbox{-proj}$, the inclusion
$\mu: \mathcal{P}\to\mathcal{X}$ of Frobenius categories induces a
fully faithful functor $\DF{\mu}: \DF{\mathcal{P}}\ra
\DF{\mathcal{X}}$. Since $q_1=q_2\DF{\mu}$, we see that $\DF{\mu}$
is an equivalence. Thus the map $K(\mu): K(\mathcal{P})\to K(\mathcal{X})$
is a homotopy equivalence by Lemma \ref{loc} (2).

It remains to show that the map $\cpx{P}\cpx{\otimes}_{\mathbb S}-:\mathcal{W}_\mathbb{S}\to \mathcal{X}$ induces a homotopy equivalence $K(\mathbb S)\lraf{\sim} K(\mathcal{X})$. In fact, since each object of $\Cb{\pmodcat A}$ is perfect, it follows from Lemma \ref{dg-0} that the functor $\cpx{P}\otimesL_{\mathbb S}-: \D{\mathbb S}\to\Tria{({_A}\cpx{P})}$ is a triangle equivalence. Thus $K(\mathbb S)\lraf{\sim} K(\mathcal{X})$ by Lemma \ref{dg-1} (2). $\square$

\subsection{Decomposition of higher algebraic $K$-groups}
In this subsection, we shall establish reduction formulas for calculation of algebraic $K$-groups of dg algebras. The main result of this subsection is Proposition \ref{dg-5}, which will be applied in the next subsection to show Theorem \ref{new-theorem}.

First, we extend a result of Berrick and Keating (see \cite{bk}) on algebraic $K$-groups of triangular matrix rings to the ones of dg triangular matrix rings.

\begin{Lem}\label{dg-4}
Let $\mathbb{R}=\left(\begin{array}{lc} \mathbb{S} &\mathbb{M} \\
0 & \mathbb{T}\end{array}\right)$ be the dg triangular matrix algebra defined by dg algebras $\mathbb{S}$, $\mathbb{T}$ and a dg $\mathbb{S}$-$\mathbb{T}$-bimodule $\mathbb{M}$.
Then
$$K_n(\mathbb R)\simeq K_n(\mathbb S)\oplus K_n(\mathbb T)\quad \mbox{for all}\;\; n\in
\mathbb{N}.$$
\end{Lem}

{\it Proof.} Let $e:=\left(\begin{array}{ll} 0 &0\\
0 & 1\end{array}\right)\in \mathbb{R}$,  $f:=\left(\begin{array}{ll} 1 &0\\
0 & 0\end{array}\right)\in\mathbb{R}$ and $\mathbb{J}:=\mathbb{R}e\mathbb{R}$.
Then $e^2=e$, $f^2=f$, $e\mathbb{R}e=\mathbb{T}$ and $\mathbb{R}/\mathbb{J}=\mathbb{S}$. On the one hand, for each $n\in\mathbb{Z}$, we have
$$\Hom_{\D{\mathbb R}}(\mathbb{R}e, \mathbb{R}f[n])\simeq\Hom_{\K{\mathbb R}}(\mathbb{R}e, \mathbb{R}f[n])\simeq H^n\big(\cpx{\Hom}_\mathbb{R}(\mathbb{R}e, \mathbb{R}f)\big)\simeq  H^n\big(e\cpx{\Hom}_\mathbb{R}(\mathbb{R}, \mathbb{R})f\big)\simeq H^n\big(e\mathbb{R}f\big)=0.$$
On the other hand, both $\mathbb{R}e$ and $\mathbb{R}f$ are compact in $\D{\mathbb{R}}$ and $\Tria(\mathbb{R}e\oplus\mathbb{R}f)=\Tria(\mathbb{R})=\D{\mathbb{R}}$. Then, by \cite[Theorem 3.3]{Jor}, there exists a recollement of derived categories of dg algebras: $$\xymatrix{\D{\mathbb{S}}\ar^-{D(\lambda_*)}[r]&\D{\mathbb{R}}\ar^-{e\mathbb{R}\cpx{\otimes}_\mathbb{R}-}[r]
\ar^-{\rHom_\mathbb{R}(\mathbb{S},-)\;}@/^1.2pc/[l]\ar_-{\mathbb{S}\otimesL_\mathbb{R}-}@/_1.6pc/[l] &\D{\mathbb{T}}\ar_-{\mathbb{R}e\otimesL_\mathbb{T}-}@/_1.6pc/[l]\ar^-{\rHom_{\mathbb{T}}(e\mathbb{R},-)}@/^1.2pc/[l] }$$
where $D(\lambda_*)$ is the restriction functor induced from the canonical surjection $:\mathbb{R}\to \mathbb{S}$. Note that the functors
$\mathbb{S}\otimesL_\mathbb{R}-$ and $\mathbb{R}e\otimesL_{\mathbb T}-$
preserve compact objects, and that ${_\mathbb{R}}\mathbb{S}=\mathbb{R}f\in\Dc{\mathbb{R}}$ and $e\mathbb{R}=\mathbb{T}\in\Dc{\mathbb{T}}$. Thus, from the above recollement
we can get the following ``half recollement" for the subcategories
of compacts objects:
$$(*)\quad
\xymatrix{\Dc{\mathbb{S}}\ar^-{D(\lambda_*)}[r]&\Dc{\mathbb{R}}\ar^-{e\mathbb{R}\cpx{\otimes}_\mathbb{R}-}[r]
\ar_-{\mathbb{S}\otimesL_\mathbb{R}-}@/_1.6pc/[l] &\Dc{\mathbb{T}}\ar_-{\mathbb{R}e\otimesL_\mathbb{T}-}@/_1.6pc/[l]}$$
This implies that the following sequence of triangulated categories
$$
 \xymatrix{\Dc{\mathbb{S}}&& \Dc{\mathbb{R}}\ar[ll]_-{\mathbb{S}\otimesL_\mathbb{R}-}&& \Dc{\mathbb{T}}\ar[ll]_-{\mathbb{R}e\otimesL_\mathbb{T}-}.}
$$
is exact. Since ${_\mathbb{R}}\mathbb{R}e\in\mathcal{W}_\mathbb{R}$ and ${_\mathbb{S}}\mathbb{S}\in\mathcal{W}_\mathbb{S}$, we see from Lemma \ref{functors} that the following functors
$$\mathbb{R}e\cpx{\otimes}_{\mathbb{T}}-:
\mathcal{W}_\mathbb{T}\longrightarrow \mathcal{W}_\mathbb{R}, \quad\quad\mbox{and}\quad\quad
\mathbb{S}\cpx{\otimes}_\mathbb{R}-: \mathcal{W}_\mathbb{R}\longrightarrow \mathcal{W}_\mathbb{S}
$$
are well-defined maps of Frobenius pairs. Moreover, by Lemma \ref{dgk}, we can construct the  following commutative diagram:
$$\xymatrix@C=1.2cm{\Dc{\mathbb{S}}
&\Dc{\mathbb{R}}\ar[l]_-{\mathbb{S}\otimesL_{\mathbb{R}}-} & \Dc{\mathbb{T}}\ar[l]_-{\mathbb{R}e\otimesL_{\mathbb{T}}-} \\
\DF{\mathcal{W}_\mathbb{S}}\ar[u]_{\simeq}
&\DF{\mathcal{W}_\mathbb{R}}\ar[l]_-{\DF{\mathbb{S}\cpx{\otimes}_{\mathbb{R}}-}}\ar[u]_{\simeq}&
\DF{\mathcal{W}_\mathbb{T}}\ar[l]_-{\DF{\mathbb{R}e\cpx{\otimes}_{\mathbb{T}}-}}\ar[u]_{\simeq}
}\vspace{0.3cm}$$
This implies that the second row is an exact sequence of triangulated categories:
$$
\xymatrix{\DF{\mathcal{W}_\mathbb{S}}&& \DF{\mathcal{W}_\mathbb{R}}\ar[ll]_-{\DF{\mathbb{S}\cpx{\otimes}_{\mathbb{R}}-}}&& \DF{\mathcal{W}_\mathbb{T}}\ar[ll]_-{\DF{\mathbb{R}e\cpx{\otimes}_{\mathbb{T}}-}}.}
$$
By Lemma \ref{loc} (1), the following sequence of maps among $K$-theory spaces
$$
\xymatrix{K(\mathbb{S})&& K(\mathbb{R})\ar[ll]_-{K(\mathbb{S}\cpx{\otimes}_{\mathbb{R}}-)}&& K(\mathbb{T})\ar[ll]_-{K(\mathbb{R}e\cpx{\otimes}_{\mathbb{T}}-)}.}
$$ is a homotopy fibration, and therefore there is a long exact
sequence of $K$-groups:
$$\cdots \lra K_{n+1}(\mathbb{S})\lra
K_n(\mathbb{T})\lraf{K_n(\mathbb{R}e\cpx{\otimes}_{\mathbb{T}}-)} K_n(\mathbb{R})\lraf{K_n(\mathbb{S}\cpx{\otimes}_{\mathbb{R}}-)}
K_n(\mathbb{S})\lra K_{n-1}(\mathbb{T})\lra$$
$$\cdots\lra K_0(\mathbb{T})\lra K_0(\mathbb{R})\lra
K_0(\mathbb{S})\lra 0$$ for all $n\in \mathbb{N}.$
It remains to show that this sequence breaks up into a series of split short exact sequences.

Actually, since $e\mathbb{R}=\mathbb{T}\in\mathcal{W}_\mathbb{T}$, the functor $e\mathbb{R}\cpx{\otimes}_\mathbb{R}-:\mathcal{W}_\mathbb{R}\to \mathcal{W}_\mathbb{T}$ is a map of Frobenius pairs due to Lemma \ref{functors}. Note that $$\big(e\mathbb{R}\cpx{\otimes}_\mathbb{R}-\big)(\mathbb{R}e\cpx{\otimes}_{\mathbb{T}}-)\simeq (e\mathbb{R}e)\cpx{\otimes}_{\mathbb{T}}-\simeq Id_{\mathcal{W}_{\mathbb{T}}}: \mathcal{W}_{\mathbb{T}}\lra \mathcal{W}_{\mathbb{T}}.$$
Thus the composite of the map $K(\mathbb{R}e\cpx{\otimes}_{\mathbb{T}}-):K(\mathbb{T})\to K(\mathbb{R})$ with the map $K(e\mathbb{R}\cpx{\otimes}_\mathbb{R}-):K(\mathbb{R})\to K(\mathbb{T})$ is homotopic to the identity map on $K(\mathbb{T})$. In view of $K_n$-groups, we have $$K_n(\mathbb{R}e\cpx{\otimes}_{\mathbb{T}}-)K_n(e\mathbb{R}\cpx{\otimes}_\mathbb{R}-)=Id_{K_n(\mathbb{T})}:
K_n(\mathbb{T})\lra K_n(\mathbb{T})$$
for each $n\in\mathbb{N}$. This implies that $K_n(\mathbb{R}e\cpx{\otimes}_{\mathbb{T}}-):K_n(\mathbb{T})\to K_n(\mathbb{R})$ is a split-injection. Combining this with the above long exact sequence,
we see that $K_n(\mathbb R)\simeq K_n(\mathbb S)\oplus K_n(\mathbb T)$ for each $n\in\mathbb{N}$. $\square$

\medskip
Now, we give the main result of this subsection.

\begin{Prop}\label{dg-5}
Let $\mathbb{A}$ be a dg algebra, and let $M$ and $N$ be two perfect dg $\mathbb{A}$-modules. Suppose that $\Hom_{\D{\mathbb A}}(M,N[i])=0$ for all $i\in\mathbb{Z}$. Then $$K_n\big(\cpx{\End}_{\mathbb A}(M\oplus N)\big)\simeq
K_n\big(\cpx{\End}_{\mathbb A}(M)\big)\oplus K_n\big(\cpx{\End}_{\mathbb A}(N)\big)
\quad \mbox{for all}\;\; n\in
\mathbb{N}.$$
If, in addition, $\D{\mathbb{A}}=\Tria(M\oplus N)$, then
$$K_n(\mathbb{A})\simeq K_n\big(\cpx{\End}_{\mathbb A}(M)\big)\oplus K_n\big(\cpx{\End}_{\mathbb A}(N)\big)\quad \mbox{for all}\;\; n\in
\mathbb{N}.$$
\end{Prop}

{\it Proof.} We define $\mathbb{B}:=\cpx{\End}_{\mathbb A}(M\oplus N)$. Then $$\mathbb{B}=\left(\begin{array}{ll} \cpx{\End}_{\mathbb A}(M) &\cpx{\Hom}_{\mathbb A}(M,N)\\
\cpx{\Hom}_{\mathbb A}(N,M)& \cpx{\End}_{\mathbb A}(N)\end{array}\right).$$
Since $M$ is perfect, it always has the property $(P)$. This implies that
$\Hom_{\K{\mathbb A}}(M,N[i]) \simeq \Hom_{\D{\mathbb A}}(M,N[i])$ for each $i\in\mathbb{Z}$.
Consequently, we have $$H^i\big(\cpx{\Hom}_{\mathbb A}(M,N)\big)=\Hom_{\K{\mathbb A}}(M,N[i])
\simeq \Hom_{\D{\mathbb A}}(M,N[i])=0$$
and therefore the following canonical inclusion:
$$
\mathbb{C}:=\left(\begin{array}{ll} \cpx{\End}_{\mathbb A}(M) & 0\\
\cpx{\Hom}_{\mathbb A}(N,M)& \cpx{\End}_{\mathbb A}(N)\end{array}\right)
\hookrightarrow
\left(\begin{array}{ll} \cpx{\End}_{\mathbb A}(M) &\cpx{\Hom}_{\mathbb A}(M,N)\\
\cpx{\Hom}_{\mathbb A}(N,M)& \cpx{\End}_{\mathbb A}(N)\end{array}\right)=\mathbb{B}
$$
is a quasi-isomorphism of dg algebras. It follows from Lemma \ref{kell} that
$K(\mathbb{B})\lraf{\sim} K(\mathbb{C})$ as $K$-theory spaces, and therefore
$K_n(\mathbb{B})\simeq K_n(\mathbb{C})$ for each $n\in\mathbb{N}$.
Further, due to
Lemma \ref{dg-4}, we have $$K_n(\mathbb{C})\simeq K_n(\cpx{\End}_{\mathbb A}(M))
\oplus K_n(\cpx{\End}_{\mathbb A}(N)).$$
Thus the first part of Proposition \ref{dg-5} follows.

To show the second part of Proposition \ref{dg-5}, we note that $K_n(\mathbb{A})\simeq K_n(\mathbb{B})$ by Corollary \ref{dg-2} because the dg $\mathbb{A}$-module $M\oplus N$ is
perfect and $\D{\mathbb{A}}=\Tria(M\oplus N)$. $\square$.

\medskip
Following \cite[Section 4]{NS}, we say that a homomorphism $\lambda:\mathbb{R}\to\mathbb{S}$ of dg algebras is a {\emph{homological epimorphism}} if the restriction functor $D(\lambda_*):\D{\mathbb{S}}\to\D{\mathbb{R}}$ is fully faithful. This is also equivalent to that the canonical homomorphism $\mathbb{S}\otimesL_{\mathbb{R}}\mathbb{S}\to \mathbb{S}$ is an isomorphism in $\D{\mathbb{S}}$. Clearly, each homological ring epimorphism is a homological epimorphism of dg algebras concentrated in degree $0$.

\begin{Koro}\label{dg-7}
Let $\lambda: \mathbb{R}\to \mathbb{S}$ be a homological epimorphism of dg algebras.
If the dg $\mathbb{R}$-module $\mathbb{S}$ is compact in $\D{\mathbb{R}}$, then there exists a dg algebra $\mathbb{T}$ determined by $\lambda$ such that $$K_n(\mathbb{R})\simeq K_n(\mathbb{S})\oplus K_n(\mathbb{T})\quad \mbox{for all}\;\; n\in
\mathbb{N}.$$
\end{Koro}

{\it Proof.} Since $\lambda$ is a homological epimorphism of dg algebras, it follows from \cite[Section 4]{NS} that there is a recollement of triangulated categories:
$$
\xymatrix@C=1.2cm{\D{\mathbb{S}}\ar[r]^-{i_*}&\D{\mathbb{R}}\ar[r]^-{j^!}
\ar@/^1.2pc/[l]\ar@/_1.2pc/[l]_-{i^*}
&{\rm{Tria}}({_\mathbb{R}}Q)\ar@/^1.2pc/[l]\ar@/_1.2pc/[l]_-{j_!}}\vspace{0.2cm}$$
where $Q$ is a dg $\mathbb{R}$-$\mathbb{R}$-bimodule such that $\mathbb{R}\lraf{\lambda} \mathbb{S}\lra Q\lra \mathbb{R}[1]$ is a distinguished triangle in $\K{\mathbb{R}\otimes_k\mathbb{R}\opp}$, and where
$j_!$ is the canonical embedding and
$j^!=Q\otimesL_{\mathbb{R}}-,\, i^*=\mathbb{S}\otimesL_{\mathbb{R}}-, \,i_*=D(\lambda_*).$ This implies that $\Hom_{\D{\mathbb{R}}}(Q, \mathbb{S}[m])=0$ for any $m\in\mathbb{Z}$, and that $\D{\mathbb{R}}=\Tria(Q\oplus \mathbb{S})$.

Assume that ${_\mathbb{R}}\mathbb{S}\in\Dc{\mathbb{R}}$. Then $Q\in\Dc{\mathbb{R}}$. As each compact object of $\D{\mathbb{R}}$ is isomorphic in $\D{\mathbb{R}}$ to a
perfect dg $\mathbb{R}$-module, there are two perfect dg $\mathbb{R}$-modules $N$
and $M$ such that $N\simeq \mathbb{S}$ and $M\simeq Q$ in $\D{\mathbb{R}}$. It follows that
$\Hom_{\D{\mathbb{R}}}(M, N[m])=0$ for any $m\in\mathbb{Z}$, and that $\D{\mathbb{R}}=\Tria(M\oplus N)$. By Proposition \ref{dg-5},
$$K_n(\mathbb{R})\simeq K_n\big(\cpx{\End}_{\mathbb R}(M)\big)\oplus
K_n\big(\cpx{\End}_{\mathbb R}(N)\big) \quad \mbox{for all}\;\; n\in
\mathbb{N}.$$
Now, we define $\mathbb{T}:=\cpx{\End}_{\mathbb R}(M)$ and $\mathbb{B}:=\cpx{\End}_{\mathbb R}(N)$. To finish the proof of Corollary \ref{dg-7}, it suffices to show that  $K_n(\mathbb{B})\simeq K_n(\mathbb{S})$ as $K$-groups for each $n\in\mathbb{N}$.

On the one hand, since ${_\mathbb{R}}N$ is perfect, the left-derived functor $N\otimesL_{\mathbb{B}}-:\D{\mathbb{B}}\to\Tria({_\mathbb{R}}N)$
is a triangle equivalence by Lemma \ref{dg-0}. On the other hand, since the functor $i_*$ is fully faithful, the adjoint pair $(i^*,i_*)$ implies that $i^*$ restricts to a triangle equivalence $\Tria({_\mathbb{R}}\mathbb{S})\lraf{\simeq}\D{\mathbb{S}}$.
Moreover, due to $N\simeq \mathbb{S}$ in $\D{\mathbb{R}}$, we have $\Tria({_\mathbb{R}}N)=\Tria({_\mathbb{R}}\mathbb{S})$. Thus the composite
$
(\mathbb{S}\otimesL_{\mathbb{R}}-)(N\otimesL_{\mathbb{B}}-):\D{\mathbb{B}}\to\D{\mathbb{S}}
$
of the functors $N\otimesL_{\mathbb{B}}-$ and $i^*$ is a triangle equivalence. Since ${_\mathbb{R}}N$ is perfect, we see that ${_\mathbb{R}}N$ has the property $(P)$ and that $\mathbb{S}\cpx{\otimes}_{\mathbb{B}}N$ is a perfect dg $\mathbb{S}$-module by Lemma \ref{functors}. As the functor $N\cpx{\otimes}_{\mathbb{B}}-:\K{\mathbb{B}}\to\K{\mathbb{R}}$ preserves dg modules with the property $(P)$, we clearly have
$$(\mathbb{S}\otimesL_{\mathbb{R}}-)(N\otimesL_{\mathbb{B}}-)
\simeq(\mathbb{S}\cpx{\otimes}_{\mathbb{R}}N)\otimesL_\mathbb{B}-:\D{\mathbb{B}}\lraf{\simeq}\D{\mathbb{S}}.$$
It follows from Lemma \ref{dg-1} (2) that $K(\mathbb{B})\lraf{\sim} K(\mathbb{S})$ as $K$-theory spaces. This gives rise to $K_n(\mathbb{B})\simeq K_n(\mathbb{S})$.   $\square$

\medskip
Applying Corollary \ref{dg-7} to homological epimorphisms of ordinary rings, we obtain the following result.

\begin{Koro}\label{prop1}
Let $\lambda: R\to S$ be a homological ring epimorphism such that $_RS$ has a finite-type resolution. Denote by $\cpx{Q}$ the two-term complex $0\ra R\lraf{\lambda} S\ra 0$ with $R$ and $S$ in degrees $0$ and $1$, respectively. Let $\cpx{P}\in \Cb{\pmodcat R}$ such that $\Tria(\cpx{P})=\Tria({_R}\cpx{Q})\subseteq\D{R}$. Then $$K_n(R)\simeq K_n(S)\oplus K_n\big(\cpx{\End}_R(\cpx{P})\big) \;\,\mbox{for all}\;\; n\in
\mathbb{N}.$$
\end{Koro}

{\it Proof.} Since $_RS$ has a finite-type resolution, we can choose a complex $\cpx{N}$ in $\mathscr{C}^b(\pmodcat{R})$ such that $_RS$ is isomorphic to $\cpx{N}$ in $\D R$. So we get a chain map from $_RR$ to $\cpx{N}$ such that its mapping cone $\cpx{M}$ is isomorphic in $\D{R}$ to $\cpx{Q}$. It follows that $\cpx{M}\in \Cb{\pmodcat R}$ and $\Tria(\cpx{M})=\Tria({_R}\cpx{Q})\subseteq\D{R}$.

Next, we regard $R$ and $S$ as dg $\mathbb{Z}$-algebras concentrated in degree $0$. Then $\lambda:R\to S$ is a homological epimorphism of dg algebras. Moreover, both $\cpx{N}$ and $\cpx{M}$ are perfect dg $R$-modules. By Lemma \ref{same} and the proof of Corollary \ref{dg-7}, we see that $K_n(R)\simeq K_n(S)\oplus K_n(\cpx{\End}_R(\cpx{M}))$ for all $n\in\mathbb{N}$.

Note that both $\cpx{M}$ and $\cpx{P}$ are perfect, and that $\Tria(\cpx{P})=\Tria({_R}\cpx{Q})=\Tria(\cpx{M})$. By Corollary \ref{dg-3}, we have $K_n\big(\cpx{\End}_{R}(\cpx{M})\big)\simeq K_n\big(\cpx{\End}_{R}(\cpx{P})\big)$. Thus $K_n(R)\simeq K_n(S)\oplus K_n(\cpx{\End}_R(\cpx{P}))$. $\square$

\begin{Rem}\label{biWaldhausen}
Let us give a comment on the relationship between Corollary \ref{prop1} and
\cite[Theorem 1.1]{xc4}. Recall that, in \cite[Theorem 1.1]{xc4}, we describe the difference between $K_n(R)$ and $K_n(S)$ by the $n$-th algebraic $K$-group of a complicial biWaldhausen category $\mathbf{W}(R,\lambda)$ defined in \cite[Theorem 14.9]{kr}.

Concretely, $\mathbf{W}(R,\lambda)$ is the full subcategory of $\Cb{\pmodcat R}$ consisting of all those complexes $\cpx{X}$ such that $S\otimes_R\cpx{X}$ is acyclic. As a Waldhausen category, it has injective chain maps which are degreewise split as cofibrations, and has homotopy equivalences as weak equivalences. In this sense, the cofibrations and weak equivalences of $\mathbf{W}(R,\lambda)$ are induced from the Frobenius pair $\big(\mathbf{W}(R,\lambda),\mathscr{C}^b_{ac}(R\pmodcat)\big)$.

In \cite[Theorem 1.1]{xc4}, it was shown that $K_n(R)\simeq K_n(S)\oplus K_n(R,\lambda)$ for all $n\in\mathbb{N}$, where $K_n(R,\lambda):=K_n(\mathbf{W}(R,\lambda))$.  Now, we point out $K_n\big(\cpx{\End}_R(\cpx{P})\big)\simeq K_n(R,\lambda)$ for the complex $\cpx{P}$ in Corollary \ref{prop1}.

In fact, since $\lambda:R\to S$ is homological, we see that $\Ker(S\otimesL_R-)=\Tria{(_R\cpx{Q})}\subseteq \D{R}$ by Lemma \ref{homological}. Note that $\mathbf{W}(R,\lambda)$ consists of all those complexes $\cpx{X}\in\Cb{\pmodcat R}$ such that $S\otimesL_R\cpx{X}=0$ in $\D{S}$. Thus $\mathbf{W}(R,\lambda)$ is the same as the full subcategory of $\Cb{\pmodcat R}$ consisting of all those complexes which, regarded as objects in $\D{R}$, belong to $\Tria{(_{R}\cpx{Q})}$. Since $\cpx{P}\in\Cb{\pmodcat R}$ and $\Tria({_R}\cpx{Q})=\Tria({_R}\cpx{P})$, we know from Corollary \ref{Qe} that $K_n\big(\cpx{\End}_R(\cpx{P})\big)\simeq K_n\big(\mathbf{W}(R,\lambda)\big)=K_n(R,\lambda)$. $\square$
\end{Rem}

Finally, we give an example to illustrate that the dg algebra in Corollary \ref{prop1} cannot be substituted by its underlying ring (just forgetting the differential). Note that, in this example, the ring homomorphism $\lambda:R\to S$ has already been considered in \cite{xc4} to illustrate \cite[Theorem 1.1]{xc4}.

Let $R$ be the following quiver algebra over a field $k$ with
relations
$$\xymatrix{1\,\bullet\ar^-{\bf{\alpha}}@/^0.8pc/[r]
&\bullet\, 2\ar^-{\bf{\beta}}@/^0.8pc/[l]},
\quad\alpha\beta=\beta\alpha=0.$$

Further, let $e_i$ be the idempotent element of $R$ corresponding to the
vertex $i$ for $i=1,2$, and let $\lambda:R\to S$ be the noncommutative localization
of $R$ at the homomorphism $\varphi: Re_2\to Re_1$ induced by $\alpha$.
Then $S$ is equal to the following quiver algebra over $k$ with relations:
$$\xymatrix{1\,\bullet\ar^-{\bf{\alpha}}@/^0.8pc/[r]
&\bullet\, 2\ar^-{\bf{\alpha^{-1}}}@/^0.8pc/[l]},
\quad\alpha\alpha^{-1}=e_1\;\;\mbox{and}\;\;\alpha^{-1}\alpha=e_2,$$
and $\lambda:R\to S$ is given
explicitly by
$$
e_1\mapsto e_1,\, e_2\mapsto e_2,\, \alpha\mapsto \alpha,
\,\beta\mapsto 0.
$$
For an explanation, we refer the reader to \cite{xc4}. Thus $S$ is isomorphic to the usual $2\times 2$ matrix ring $M_2(k)$ over $k$.
Note that $Se_2\simeq Se_1\simeq Re_1$ and that $S\simeq
Se_1\oplus Se_2\simeq Re_1\oplus Re_1$ as $R$-modules. Hence $\lambda$ is a homological ring epimorphism such that ${_R}S$ is finitely generated and projective.

In \cite{xc4}, we show that $K_n(R)\simeq K_n(S)\oplus K_n(R,\lambda)$ for all $n\in\mathbb{N}$, where $\mathbf{W}(R,\lambda)$ coincides with the full subcategory of $\Cb{\pmodcat R}$ consisting of those complexes $\cpx{X}$ such that $H^i(\cpx{X})\in\add(S_1)$ for all $i\in\mathbb{Z}$.
Here $S_1$ denotes the simple $R$-module corresponding to the vertex $1$.

Now, we follow Corollary \ref{prop1} and Remark \ref{biWaldhausen} to describe $K_n(R,\lambda)$ as the $K_n$-group of a dg algebra.

Let $$\cpx{Q}:= 0\lra Re_2\lraf{\varphi} Re_1\lra 0\quad
\mbox{and}\quad \cpx{P}:=0\lra R\lraf{\lambda} S\lra 0$$ where
$Re_2$ and $R$ are of degree $0$. Clearly, $\cpx{Q}\in\Cb{\pmodcat
R}$ and $\cpx{P}[1]$ is the mapping cone of $\lambda$. Since
$Se_2\simeq Se_1\simeq Re_1$ as $R$-modules, we infer that
$\cpx{Q}\simeq \cpx{P}$ in $\C{R}$ and
$\Tria(\cpx{Q})=\Tria({_R}\cpx{P})\subseteq \D{R}$. Thus all the
assumptions of Corollary \ref{prop1} are satisfied. It follows
from Corollary \ref{prop1} that
$$K_n(R)\simeq K_n(S)\oplus K_n(\mathbb{T})\;\;\mbox{for all}\;\;n\in\mathbb{N},$$
where $\mathbb{T}:=\cpx{\End}_R(\cpx{Q})$ is the dg endomorphism
algebra of the complex $\cpx{Q}$ (see Subsection \ref{subsection2.1} for
definition). By Remark \ref{biWaldhausen}, we also have $K_n(R,\lambda)\simeq K_n(\mathbb{T})$ for all $n\in\mathbb{N}$.

It is easy to check that the dg algebra
$\mathbb{T}:=(T^i)_{i\in\mathbb{Z}}$ is given by the following data:
$$
T^{-1}=k,\,T^0=k\oplus k, \,T^1=k, T^i=0\;\,\mbox{for}\;\, i\neq -1,
0, 1,
$$
with the differential:
$$ 0\lra T^{-1}\lraf{0}T^0\lraf{\left({1 \atop{-1}}\right)}T^1 \lra 0
$$
and the multiplication $\circ: \mathbb{T}\times \mathbb{T}\to
\mathbb{T}$ (see Subsection \ref{subsection2.1}):
$$
T^{-1}\circ T^{-1}=T^1\circ T^1=0=T^{-1}\circ T^1=T^1\circ T^{-1},\,
$$
$$
(a,b)\circ (c,d)=(ac,bd),\, f\circ (a, b)=fa,\; (a,b)\circ f=bf,\;
g\circ (a,b)=gb,\;(a,b)\circ g=ag,$$ where $(a,b), (c,d)\in T^0$,
$f\in T^{-1}$ and $g\in T^1$.

Since $H^1(\mathbb T)=0$, the dg algebra $\mathbb{T}$ is
quasi-isomorphic to the following dg algebra $\tau^{\leq 0}(\mathbb
T)$ over $k$:
$$
0\lra T^{-1}\lraf{0} \Ker(d^0)\lra 0
$$
where $d^0={\left({1\atop{-1}}\right)}: T^0\to T^1$. Clearly, the
latter algebra is isomorphic to the dg algebra
$$\mathbb{A}:=0\lra k \lraf{0} k\lra 0$$
where the first $k$ is of degree $-1$ and has a $k$-$k$-bimodule
structure via multiplication. Thus the algebra structure of
$\mathbb{A}$ (by forgetting its differential) is precisely the
trivial extension $k\ltimes k$ of $k$ by the bimodule $k$. Now, by
Lemma \ref{kell}, we know that
$$K(\mathbb{T})\lraf{\sim}K(\tau^{\leq 0}(\mathbb T))\lraf{\sim}
K(\mathbb A)$$ as $K$-theory spaces. This implies that $K_n(\mathbb{T})\simeq K_n(\mathbb{A})$,
and therefore $K_n(R,\lambda)\simeq K_n(\mathbb{A})$.
Thus $$K_n(R)\simeq K_n(S)\oplus K_n(\mathbb{T})\simeq K_n(S)\oplus K_n(\mathbb A)\;\;\mbox{for all}\;\;n\in\mathbb{N}.$$  It is worth noting that we cannot replace the dg
algebra $\mathbb{A}$ in the above isomorphism by the trivial
extension $k\ltimes k$ since the algebraic $K$-theory of dg algebras
is different from that of usual rings. In fact, in this example,
$K_1(R)=K_1(k)\oplus K_1(k)=k^{\times}\oplus k^{\times}$,
$K_1(S)=k^{\times}$ and $K_1(\mathbb{A})=k^{\times}$, but
$K_1(k\ltimes k)=k\oplus k^{\times}$. So $K_1(R)\not\simeq
K_1(S)\oplus K_1(k\ltimes k)$.

\subsection{ Proofs of Theorem \ref{new-theorem} and Corollary \ref{corollary1} }\label{sub-section5.1}

With previous preparations, now we prove the first two results in the introduction.

{\bf Proof of Theorem \ref{new-theorem}.}

 We regard the ordinary ring $R$ as a dg algebra $\mathbb{R}$ concentrated in degree $0$. Then
 $\C{R}$ is exactly the category of dg $\mathbb{R}$-modules, and $\D{R}$ coincides with $\D{\mathbb{R}}$. Moreover, by Lemma \ref{same}, $K(R)$ is homotopy equivalent to $K(\mathbb{R})$ as $K$-theory spaces, and therefore $K_n(R)\simeq K_n(\mathbb{R})$ for all $n\in\mathbb{N}$.
 Note that a complex $X\in\D{R}$ is compact if and only if $X$ is quasi-isomorphic to a complex $Y\in\Cb{\pmodcat R}$.

Now, we assume that there exists a recollement
$$\xymatrix@C=1.2cm{\D{S}\ar[r]^-{i_*}&\D{\mathbb{R}}\ar[r]
\ar@/^1.2pc/[l]\ar@/_1.2pc/[l]_{i^*}
&\D{T}\ar@/^1.2pc/[l]\ar@/_1.2pc/[l]_-{j_!}}\vspace{0.3cm}$$ such that
$i_*(S)$ is compact in $\D{\mathbb{R}}$. On the one hand, the dg module $j_!(T)$ is always compact in $\D{\mathbb{R}}$ by Lemma \ref{prop-recell} (3). On the other hand, we see that $\D{\mathbb{R}}=\Tria\big(j_!(T)\oplus i_*(S)\big)$, and that
$$\Hom_{\D{\mathbb R}}(j_!(T), i_*(S)[m])\simeq \Hom_{\D{\mathbb R}}(i^*j_!(T), S[m])=0$$
for each $m\in\mathbb{Z}$ because $(i^*,i_*)$ is an adjoint pair and $i^*j_!=0$.
Recall that each compact object of $\D{\mathbb{R}}$ is isomorphic in $\D{\mathbb{R}}$ to a
perfect dg $\mathbb{R}$-module. So, there exist two perfect dg $\mathbb{R}$-modules $M$
and $N$ such that $M\simeq j_!(T)$ and $N\simeq i_*(S)$ in $\D{\mathbb{R}}$.
Consequently, we have $\Tria(M\oplus N)=\D{\mathbb{R}}$ and $\Hom_{\D{\mathbb R}}(M, N[m])=0$
for all $m\in\mathbb{Z}$. It follows from Proposition \ref{dg-5} that
$$K_n(\mathbb{R})\simeq K_n\big(\cpx{\End}_{\mathbb R}(M)\big)\oplus
K_n\big(\cpx{\End}_{\mathbb R}(N)\big) \quad \mbox{for all}\;\; n\in
\mathbb{N}.$$ Next, we claim that
$$(1)\quad H^m\big(\cpx{\End}_{\mathbb R}(M)\big)\simeq\left\{\begin{array}{ll} 0 & \mbox{if}\; m\neq 0,\\ T& \mbox{if }\; m=0,\end{array} \right.
\;\;\mbox{and that}\;\;\;(2)\quad
H^m\big(\cpx{\End}_{\mathbb R}(N)\big)\simeq\left\{\begin{array}{ll} 0 & \mbox{if}\; m\neq 0,\\
S& \mbox{if }\; m=0.\end{array} \right.$$
In fact, since ${_\mathbb{R}}M$ is perfect, it has the property $(P)$. This implies that
$$H^m\big(\cpx{\End}_{\mathbb R}(M)\big)=\Hom_{\K{\mathbb R}}(M, M[m])
\simeq\Hom_{\D{\mathbb R}}(M, M[m]).$$
As $j_!:\D{T}\to\D{\mathbb{R}}$ is fully faithful, we have
$$\Hom_{\D{\mathbb R}}(M, M[m])\simeq
\Hom_{\D{\mathbb R}}(j_!(T), j_!(T)[m])\simeq \Hom_{\D{T}}(T, T[m])
\simeq\left\{\begin{array}{ll} 0 & \mbox{if}\; m\neq 0,\\
T& \mbox{if }\; m=0.\end{array} \right.$$
Thus $(1)$ holds. Similarly, we can show that $(2)$ also holds since ${_\mathbb{R}}N$ is perfect and $i_*:\D{S}\to\D{\mathbb{R}}$ is fully faithful.

Now, by Corollary \ref{kell}, it follows from $(1)$ and $(2)$ that $K_n\big(\cpx{\End}_{\mathbb R}(M)\big)\simeq K_n(T)$ and  $K_n\big(\cpx{\End}_{\mathbb R}(N)\big)\simeq K_n(S)$, respectively.  Thus  $K_n(R)\simeq K_n(\mathbb{R})\simeq K_n(S)\oplus K_n(T)$ for each $n\in\mathbb{N}$.  $\square$

\smallskip
We remark that the proof of Theorem \ref{new-theorem} shows a little bit more: If  $\D{R}$ in Theorem \ref{new-theorem} is replaced by $\D{\mathbb R}$ with $\mathbb{R}$ a dg algebra, then $K_n(\mathbb{R})\simeq K_n(S)\oplus K_n(T)$.

\medskip
{\bf Proof of Corollary \ref{corollary1}.}

Note that $i_*(S)=S$ and that $_RS$ is quasi-isomorphic to a bounded complex of finitely generated projective $R$-modules if and only if $_RS$ admits a finite-type resolution. So, for the recollement in Corollary \ref{corollary1}, if $_RS$ has a finite-type resolution, then it follows from Theorem \ref{new-theorem} that $K_n(R)\simeq K_n(S)\oplus K_n(T)$ for each $n\in\mathbb{N}$.
Similarly, we can prove Corollary \ref{corollary1} for the case
that $S_R$ has a finite-type resolution. In fact, this can be
understood from Lemma \ref{dual} and the following fact: For any
ring $A$, there is a homotopy equivalence $K(A)\lraf{\sim}K(A\opp)$
(see \cite[Sections 1 (3) and 2 (5)]{Quillen}). Thus Corollary \ref{corollary1}
follows.  $\square$.

\section{Applications to algebraic $K$-theory of homological exact contexts}\label{Exact context}

In this section, we apply our results to algebraic $K$-theory of exact contexts (see  \cite{xc3}). We mainly concentrate on two classes of exact contexts, one is induced from noncommutative localizations, and the other is from the free products of groups.

\subsection{$K$-theory of noncommutative localizations}

First, we recall some results about noncommutative localizations in algebraic $K$-theory (see \cite{nr, ne2}).

Let $\Sigma$ a set of homomorphisms between finitely generated projective $R$-modules. By abuse of notation, we always identify each map $P_1\lraf{f} P_0$ in $\Sigma$ with the two-term complex $0\ra P_1\lraf{f} P_0\ra 0$ in $\Cb{\pmodcat R}$, where $P_i$ is in
the degree $-i$ for $i=0,1$. Moreover, let $\lambda_\Sigma: R\to R_\Sigma$
be the noncommutative localization  of $R$ at $\Sigma$. Note that the terminology ``noncommutative localization" was originally called ``universal localization" in the literature (for example, see \cite[Part I, 4]{schofield}). Moreover, $\lambda_\Sigma$ is a ring epimorphism with $\Tor^R_1(R_\Sigma,R_\Sigma)=0$.

Now, we recall from \cite[Definition 0.4]{nr}the definition of a small Waldhausen category $\mathbf{R}$. Precisely, the category $\mathbf{R}$ is the smallest full subcategory of $\Cb{\pmodcat R}$ which

$(i)$ \;\,contains all the complexes in $\Sigma$,

$(ii)$ \;contains all acyclic complexes,

$(iii)$ is closed under the formation of mapping cones and shifts,

$(iv)$ contains any direct summands of any of its objects.

\smallskip
We remark that, in $\mathbf{R}$, the cofibrations are injective chain maps which are degreewise split, and the weak equivalences are homotopy equivalences. So, the cofibrations and weak equivalences of $\mathbf{R}$ are exactly induced from the Frobenius pair $\big(\mathbf{R},\mathscr{C}^b_{ac}(R\pmodcat)\big)$. Following Remark \ref{biWaldhausen}, let $\mathbf{W}(R,\lambda_{\Sigma})$ be the full subcategory of $\Cb{\pmodcat R}$ consisting of all those complexes $\cpx{X}$ such that $R_\Sigma\otimes_R\cpx{X}$ is acyclic. Then $\mathbf{R}=\mathbf{W}(R,\lambda_{\Sigma})$ as Waldhausen categories by the proof of  \cite[Corollary 1.2]{xc4}.

\medskip
The following result follows from \cite[Theorem 0.5]{nr}.

\begin{Lem} \label{KT}
If $\lambda_\Sigma:R\to R_\Sigma$ is homological, then there is a weak homotopy fibration of $K$-theory spaces:
$$
\xymatrix{K(\mathbf{R})\ar[rr]^-{K(F)}&& K(R)\ar[rr]^-{K(\lambda_\Sigma)}&& K(R_\Sigma)}
$$
where $F:\mathcal{R}\to \Cb{\pmodcat R}$ is the inclusion. In this case, we have a long exact sequence of $K$-groups:
$$
\xymatrix{\cdots\lra K_{n+1}(R_\Sigma)\lra
K_n(\mathbf{R})\ar[r]^-{K_n(F)} & K_n(R) \ar[r]^-{K_n(\lambda_\Sigma)}&
K_n(R_\Sigma)\lra K_{n-1}(\mathbf{R})\lra }
$$
$$\cdots\lra K_0(\mathbf{R})\lra K_0(R)\lra
K_0(R_\Sigma)$$ for all $n\in \mathbb{N}.$
\end{Lem}

\medskip
One of the significant methods for calculating $K$-groups is to have a kind of long
Mayer-Vietoris sequences which link $K$-groups of rings together. In the following, we shall establish some long exact sequences of this type, which are induced from homological exact contexts introduced in \cite{xc3}.

We follow all the notations introduced in Section $1$.

Let $(\lambda,\mu,M,m)$ be an arbitrary but fixed exact context, where $\lambda:R\to S$ and $\mu:R\to T$ are ring homomorphisms, and where $M$ is an $S$-$T$-bimodule with an element $m$. Then there is map $\gamma: S\otimes_RT\ra M$, defined by $s\otimes t\mapsto smt$ for $s\in S$ and $t\in T$. Also, by the definition of exact context, there is a map $\beta: M\ra T\otimes_RS$ which makes the following diagram commutative:
$$\xymatrix{R\ar@{=}[d]\ar[r]^-{(\lambda,\,\mu)\,} &S\oplus T \ar@{=}[d]\ar[r]^-{\left({\cdot m\,\atop{-m \cdot}}\right)}
                & M \ar[d]^-{\beta}\\
 R\ar[r]^-{(\lambda,\,\mu)\,} & S\oplus T\ar[r]^-{\left({\rho\,\atop{-\phi}}\right)}& T\otimes_RS}$$
where
$\rho=\mu\otimes S:\; S\ra T\otimes_RS,\; s\mapsto 1\otimes s \quad \mbox{and}\quad \phi=T\otimes\lambda:\; T\ra T\otimes_RS,\; t\mapsto t\otimes 1.$ In \cite{xc5}, we use $\gamma\beta$ as a twisting to define a ring structure on $T\otimes_RS$, called the \emph{noncommutative tensor product} of $(\lambda,\mu,M,m)$ and denoted $T\boxtimes_RS$ (see \cite[Section 4.1]{xc5} for details).

Define $$B:=\left(\begin{array}{lc} S & M\\
0 & T\end{array}\right), \;e_1:= \left(\begin{array}{ll} 1 &0\\
0 & 0\end{array}\right)\quad \mbox{and}\quad  e_2:=\left(\begin{array}{ll} 0 &0\\
0 & 1\end{array}\right)\in B.$$
Then, by \cite[Lemma 5.1]{xc3}, the noncommutative localization of $B$ at the map:
$$\varphi: Be_1\lra Be_2:\;\left(\begin{array}{l} s \\
0 \end{array}\right)\mapsto \left(\begin{array}{c} sm\\
0 \end{array}\right)\;\;\mbox{for}\;\; s\in S,$$
is given by the following ring homomorphism:
$$\theta: \;B=\left(\begin{array}{lc} S & M\\
0 & T\end{array}\right)\lraf{\left(\begin{array}{cc} \rho & \beta\\
0 & \phi\end{array}\right)} \left(\begin{array}{lc} T\boxtimes_RS & T\boxtimes_RS\\
T\boxtimes_RS& T\boxtimes_RS\end{array}\right)=C.$$

Furthermore, let $\cpx{P}$ be the complex $$0\ra Be_1\lraf{\varphi} Be_2\ra 0$$
over $B$ with $Be_1$ and $Be_2$ in degrees $-1$ and $0$, respectively. Note that $Be_1$ and $Be_2$ are also right $R$-modules via $\lambda$ and $\mu$, respectively, and that the map $\cdot m: S\to M$ is a homomorphism of $S$-$R$-bimodules. Thus $\varphi$ is actually a homomorphism of $B$-$R$-bimodules, and therefore $\cpx{P}$ is a bounded complex over $B\otimes_\mathbb{Z}R\opp$. Since $_B\cpx{P}\in\Cb{\pmodcat B}$, it makes sense to discuss the tensor functor $\cpx{P}\cpx{\otimes}_R-:\Cb{\pmodcat R}\to \Cb{\pmodcat B}$.

\smallskip
Let $\mathbf{W}(B,\theta)$ be the full subcategory of $\Cb{\pmodcat B}$ consisting of those complexes $\cpx{X}$ such that $C\otimes_B\cpx{X}$ is acyclic. As a Waldhausen category, $\mathbf{W}(B,\theta)$ is exactly induced from the Frobenius pair $\big(\mathbf{W}(B,\theta),\mathscr{C}^b_{ac}(B\pmodcat)\big)$. Now, we regard  $\mathbf{W}(B,\theta)$ as a Frobenius subcategory of $\Cb{\pmodcat B}$, and define $K_n(B,\theta):=K_n\big(\mathbf{W}(B,\theta)\big)$.

\begin{Lem}\label{fraction}
The functor $\cpx{P}\cpx{\otimes}_R-:\Cb{\pmodcat R}\to \Cb{\pmodcat B}$ induces a homotopy equivalence from $K(R)$ to $K\big(\mathbf{W}(B,\theta)\big)$.
\end{Lem}

{\it Proof.} We first claim that the functor $\cpx{P}\cpx{\otimes}_R-:\C{\pmodcat R}\to \C{\pmodcat B}$ factorizes through the inclusion  $\mathbf{W}(B,\theta)\hookrightarrow \Cb{\pmodcat B}$.

In fact, since $\theta$ is the noncommutative localization of $B$ at $\varphi$, the map $C\otimes_B\varphi:C\otimes_BBe_1\to C\otimes_BBe_2$ is an isomorphism of $C$-modules. This implies that the complex $C\otimes\cpx{P}$ is acyclic. Thus $\cpx{P}\in\mathbf{W}(B,\theta)$ and $\cpx{P}\otimes_R-:\Cb{\pmodcat R}\lra \Cb{\pmodcat B}$ admits a factorisation as follows:
$$\Cb{\pmodcat R}\lraf{G}\mathbf{W}(B,\theta)\hookrightarrow\Cb{\pmodcat B}.$$
Note that $\mathbf{W}(B,\theta)$ is a Frobenius subcategory of $\Cb{\pmodcat R}$ such that its derived category $\DF{\mathbf{W}(B,\theta)}$ is equal to the full subcategory of $\Kb{\pmodcat B}$ consisting of all objects of $\mathbf{W}(B,\theta)$. In this sense, $G$ is a map of Frobenius pairs.

Next, we show that the map $K(G):K(R)\to K\big(\mathbf{W}(B,\theta)\big)$ induced from $G$ is a homotopy equivalence.

Actually, by \cite[Lemma 5.4]{xc3}, the left-derived functor $\cpx{P}\otimesL_R-:\D{R}\to \D{B}$ is fully faithful. This induces a triangle equivalence $$\D{R}\lraf{\simeq} \Tria(\cpx{P})$$
which restricts to an equivalence between full subcategories of compact objects: $$\Dc{R}\lraf{\simeq}\Tria(\cpx{P})^c.$$
Since $\D{B}=\Tria(B)$ and $\cpx{P}\in\Dc{B}$, we see from \cite[Theorem 4.4.9]{ne3} that $\Tria(\cpx{P})^c=\Tria(\cpx{P})\cap \Dc{B}$. Further, by \cite[Theorem 0.11]{ne2}, the category $\Tria(\cpx{P})^c$ coincides with the full subcategory of $\Dc{B}$ consisting of all those complexes $\cpx{X}$ such that $C\otimesL_B\cpx{X}=0$ in $\D{C}$. Now, we identify $\Kb{\pmodcat{R}}$ and $\Kb{\pmodcat{B}}$ with $\Dc{R}$ and $\Dc{B}$ up to triangle equivalences, respectively. Then $\cpx{P}\otimes_R-:\Kb{\pmodcat R}\to\Kb{\pmodcat B}$ induces a triangle equivalence
$$
\DF{G}: \DF{\Cb{\pmodcat R}}=\Kb{\pmodcat R}\lraf{\simeq}\DF{\mathbf{W}(B,\theta)}.
$$
By Lemma \ref{loc} (2), the map $K(G):K(R)\to K\big(\mathbf{W}(B,\theta)\big)$ is a homotopy equivalence. $\square$

\medskip
As a preparation for the proof of Theorem \ref{K-theory}, we need the following result (see \cite[Theorem 1.1]{xc3}).

\begin{Lem} \label{hexact}
Let  $(\lambda,\mu,M,m)$ be a homological exact context.
Then the ring homomorphism $\theta:B\to C$ is a homological noncommutative localization, and there is a recollement of derived module categories:
$$
\xymatrix@C=1.2cm{\D{C}\ar[r]^-{D(\theta_*)}
&\D{B}\ar[r]^-{\;j^!}\ar@/^1.2pc/[l]\ar_-{C\otimesL_B-}@/_1.2pc/[l]
&\D{R} \ar@/^1.2pc/[l]\ar@/_1.2pc/[l]_{j_!\;}}$$

\medskip
\noindent where $D(\theta_*)$ is the restriction functor induced by
$\theta$, and where
$$j_!={_B}\cpx{P}\otimesL_R-\;\; \mbox{and}\;\; j^!=\cpx{\Hom}_B(\cpx{P},-).$$
\end{Lem}

\medskip
For homological exact contexts, we can establish the following result, which links $K$-theory spaces of rings involved in exact contexts together.

\begin{Lem}\label{free-coproduct}
Let $(\lambda,\mu,M,m)$ be a homological exact context. Then the sequence of $K$-theory spaces:
$$
\xymatrix{K(R)\ar[rr]^-{\big(-K(\lambda),\, K(\mu)\big)} &&
K(S)\times K(T) \ar[rr]^-{\left({K(\rho)\; \atop{K(\phi)}}\right)}
&& K(T\boxtimes_RS)}
$$
is a weak homotopy fibration, where $-K(\lambda)$ denotes the map $K(S[1]\cpx{\otimes}_R-):K(R)\to K(S)$ induced from the functor $S[1]\cpx{\otimes}_R-:\Cb{\pmodcat R}\to \Cb{\pmodcat S}$.

\end{Lem}

{\it Proof.} By the proof of Lemma \ref{fraction}, we have a factorisation of the functor $\cpx{P}\otimes_R-:\Cb{\pmodcat R}\to \Cb{\pmodcat B}$ as
$$\Cb{\pmodcat R}\lraf{G}\mathbf{W}(B,\theta)\lraf{F}\Cb{\pmodcat B}$$
such that $K(G):K(R)\to K\big(\mathbf{W}(B,\theta)\big)$ is a homotopy equivalence, where $F:\mathbf{W}(B,\theta)\to\Cb{\pmodcat B}$ is the inclusion. Since $(\lambda,\mu,M,m)$ is a homological exact context, the ring homomorphism $\theta:B\to C$ is a homological noncommutative localization by Lemma \ref{hexact}. Then, it follows from Lemma \ref{KT} that there is a weak homotopy fibration: $$ K\big(\mathbf{W}(B,\theta)\big)\lraf{K(F)} K(B)\lraf{K(\theta)} K(C).$$

Next, we simplify $K(B)$ and $K(C)$ up to homotopy equivalence.

Indeed, let $i: S\times T\to B$ be the inclusion and let $j: B\to S\times T=B/M$ be the projection. Since $B$ is a triangular matrix ring with $S$ and $T$ in the diagonal, it is known that the ring homomorphisms $i$ and $j$ induce inverse homotopy equivalences: $$K(i):K(S)\times K(T)\lraf{\sim} K(B)\quad \mbox{and}\quad
K(j):K(B)\lraf{\sim} K(S)\times K(T).$$
For a proof of this result using Waldhausen $K$-theory, one may refer to the proof of
\cite[Proposition 5.7 (iv)]{RS}, where the additivity theorem for Quillen $K$-theory
(see \cite[Proposition 1.3.2 (4)]{wald}) was applied. Note that the isomorphisms $K_n(B)\simeq K_n(S)\oplus K_n(T)$ were first obtained by Berrick and Keating (see \cite{bk}).

\smallskip
Now, we define $e:=\left(\begin{array}{ll} 1 &0\\0 & 0\end{array}\right)\in C$ and $\Lambda:=T\boxtimes_RS$. Since the functor $eC\otimes_C-:\pmodcat C\to \pmodcat \Lambda$ is an equivalence of categories,  we see that $K(eC\otimes_C-):K(C)\lraf{\sim} K(\Lambda)$. Moreover, there are the following natural isomorphisms of exact functors:
$$
\big(B/M\otimes_B-\big)\big(\cpx{P}\otimes_R-\big)\lraf{\simeq} \big(S[1]\otimes_R-,\, T\otimes_R-\big):\;\Cb{\pmodcat R}\lra\Cb{\pmodcat S}\times \Cb{\pmodcat T},
$$
$$
\big(eC\otimes_B-\big)\big(C\otimes_B-\big)\big(B\otimes_{(S\times T)}-\big)\lraf{\simeq} \big(\Lambda\otimes_S-\big) \oplus \big(\Lambda\otimes_T-\big):\;\pmodcat {S} \times \pmodcat{T}\lra \pmodcat{\Lambda}.
$$

\smallskip
With the above preparations, we can construct the following commutative diagram of $K$-theory spaces (up to homotopy equivalence):
$$
\xymatrix{
K\big(\mathbf{W}(B,\theta)\big)\ar[rr]^-{K(F)} && K(B)\ar@/_1.0pc/[d]_-{K(j)}\ar[rr]^-{K(\theta)} && K(C)\ar[d]^-{\wr}_-{K(eC\otimes_C-)}\\
K(R)\ar[u]^-{K(G)}_-{\wr}\ar[rr]^-{\big(-K(\lambda),\, K(\mu)\big)} && K(S)\times K(T)\ar[rd]^-{K(\rho)\times K(\phi)}\ar@/_1.0pc/[u]_-{K(i)}\ar[rr]^-{\left({K(\rho)\; \atop{K(\phi)}}\right)} && K(\Lambda)\\
& & & K(\Lambda)\times K(\Lambda)\ar[ru]^-{K(\oplus)}  && }
$$
where $K(\oplus)$ is the map induced from the coproduct functor $\oplus:\; \pmodcat {\Lambda} \times \pmodcat{\Lambda} \to \pmodcat{\Lambda}$. Note that the first row is a weak homotopy fibration. This means that the second row is also a weak homotopy fibration. Thus the proof is completed. $\square$

\medskip
As a byproduct of Lemma \ref{free-coproduct}, we have the following corollary, which says that,
although the multiplication of the noncommutative tensor product $T\boxtimes_RS$ of a homological exact context $(\lambda,\mu,M,m)$ depends on the pair $(M,m)$, the loop space of the $K$-theory space of $T\boxtimes_RS$ is independent of the pair $(M,m)$, up to homotopy equivalence. For definitions of loop spaces and homotopy fibres, we refer the reader to Section \ref{subsection3.1} for details.

\begin{Koro}\label{HTP}
Let $(\lambda,\mu,M,m)$ be a homological exact context. Then the loop space $\Omega\big(K(T\boxtimes_RS)\big)$ of the $K$-theory space $K(T\boxtimes_RS)$ is homotopy equivalent to the homotopy fibre of the map $$\big(-K(\lambda),\, K(\mu)\,\big): K(R)\lra K(S)\times K(T).$$
\end{Koro}

{\bf Proof of Theorem \ref{K-theory}.}

$(1)$ Note that the long exact sequence of $K$-groups in $(1)$ is exactly the one of homotopy groups (see Section \ref{subsection3.1}) induced from the weak homotopy fibration in Lemma \ref{free-coproduct}.

$(2)$ Recall that we have a recollement of derived module categories described in Lemma \ref{hexact}. Due to Corollary \ref{corollary1}, in order to show $(2)$, it is sufficient to prove that $_BC$ (respectively, $C_B$) has a finite-type resolution if and only if so is $_RS$
(respectively, $T_R$).

In fact,  ${_B}C$ has a finite-type resolution if and only if $_BC\in\Dc{B}$. Applying Lemma \ref{prop-recell} (3) to the recollement in Lemma \ref{hexact}, we see that $D(\theta_*)(C)={_B}C\in\Dc{B}$) if and only if $j^!(B)\in\Dc{R}$. However,
$$j^!(B)=\cpx{\Hom}_B(\cpx{P}, B)\simeq\big(S\oplus \Cone(\lambda)\big)[-1]\in\D{R}$$ where $\Cone(\lambda)$ stands for the two-term complex $0\to R\lraf{\lambda} S\to 0$ with $R$ of degree $-1$. This implies that $j^!(B)\in\Dc{R}$ if and only if $_RS\in\Dc{R}$, while the latter is equivalent to saying that $_RS$ has a finite-type resolution. Thus $_BC$ has a finite-type resolution if and only if so is $_RS$.

Note that, for the ring homomorphisms $\mu\opp:R\opp\to T\opp$ and $\lambda\opp:R\opp\to S\opp$
of opposite rings, the quadruple $(\mu\opp, \lambda\opp, _{T\opp}M_{S\opp}, m)$ is also a homological exact context. In a similar way, we can show that $C_B$ has a finite-type resolution if and only if so is $T_R$ . $\square$

\medskip
{\bf Proof of Corollary \ref{milnorsquare}.}

Since $(i_1, i_2)$ is an exact pair, we know from \cite[Remark 5.2]{xc3} that $R_2\boxtimes_RR_1$ is isomorphic to the coproduct $R_1\sqcup_RR_2$. We can check, however, that $R_1\sqcup_RR_2$ is isomorphic to $R'$ (see also \cite[Lemma 2.3]{xc6}). Thus Corollary \ref{milnorsquare} follows immediately from Theorem \ref{K-theory} (1). $\square$.

\medskip
{\bf Proof of Corollary \ref{ring extension}.}

Let $\lambda:R\to S$ be the inclusion, $\pi:S\to S/R$ the
canonical surjection and $\lambda':R\to S'$ the induced map by right multiplication.
Recall from \cite[Section 3]{xc3} that the quadruple $\big(\lambda, \lambda', \Hom_R(S,S/R),\pi \big)$ is an exact context. So the noncommutative tensor product $S'\boxtimes_RS$ of this exact context is well defined. Assume that $_RS$ is finitely generated and projective. Then  $\Tor^R_i(S', S)=0$ for all $i\geq 1$. This means that this exact context is homological.
Now, Corollary \ref{ring extension} follows from Theorem \ref{K-theory} (2).  $\square$

\medskip
As a consequence of Theorem \ref{K-theory} (1), we reobtain the
following result of Karoubi \cite[Chapter V, Proposition 7.5
(2)]{weibelbook}.

\begin{Koro}\label{Karoubi} Let $A$ and $B$ be arbitrary rings, and
let $f:A\to B$ be a ring homomorphism and $\Phi$ a central
multiplicatively closed set of nonzerodivisors in $A$ such that the
image of $\Phi$ under $f$ is a central set of nonzerodivisors in
$B$. Assume that $f$ induces a ring isomorphism $A/sA\lraf{\simeq}
B/sB$ for each $s\in\Phi$. Then there is a Mayer-Vietoris sequence
$$\cdots\lra K_{n+1}(\Phi^{-1}B)\lra K_n(A)\lra K_n(\Phi^{-1}A)\oplus K_n(B) \lra K_n(\Phi^{-1}B)\lra
K_{n-1}(A)\lra
$$
$$\cdots\lra K_0(A)\lra K_0(\Phi^{-1}A)\oplus K_0(B)\lra
K_0(\Phi^{-1}B)$$ for all $n\in \mathbb{N},$ where $\Phi^{-1}A$
stands for the localization of $A$ at $\Phi$.
\end{Koro}

{\it Proof.} Define $R:=A$, $S:=\Phi^{-1}A$, $T:=B$ and $\mu:=f$.
Let $\lambda:R\to S$ be the canonical map of the localization. By
\cite[Lemma 6.2]{xc1}, we have $S\sqcup_R T=\Phi^{-1}B$, which is
defined by the canonical maps $\rho:\Phi^{-1}A\to\Phi^{-1}B$ and
$\phi:B\to \Phi^{-1}B$. Since  $\Phi$ and $(\Phi)f$ do not contain
zerodivisors, both $\lambda$ and $\phi$ are injective. As the
modules $_A\Phi^{-1}A$ and $_B\Phi^{-1}B$ are flat, both $\lambda$
and $\phi$ are homological ring epimorphisms.

Now, we claim that $(\lambda,\mu)$ is an exact pair. To show this,
we first prove that the following well-defined map
$$h:\Phi^{-1}A\otimes_AB\lra
\Phi^{-1}B,\;\; a/s\otimes b\mapsto ((af)b)/(sf)$$ for $a\in A$,
$s\in \Phi$ and $b\in B$, is an isomorphism of
$\Phi^{-1}A$-$B$-bimodules. In fact, since
$\Phi^{-1}A=\displaystyle\varinjlim_{s\in \Phi} s^{-1}A$, where
$s^{-1}A:=\{a/s\mid a\in A\}\subseteq \Phi^{-1}A$, we have
$$\Phi^{-1}A\otimes_AB=(\displaystyle\varinjlim_{s\in \Phi}
s^{-1}A)\otimes_AB\lraf{\simeq}\displaystyle\varinjlim_{s\in \Phi}
(s^{-1}A\otimes_AB)\lraf{\simeq}\displaystyle\varinjlim_{s\in \Phi}
(sf)^{-1}B=\Phi^{-1}B.$$

Next, we show that the cokernels of $\lambda$ and $\phi$ are
isomorphic as $A$-modules. Actually,
$$
\Phi^{-1}A/A=(\displaystyle\varinjlim_{s\in \Phi}
s^{-1}A)/A\lraf{\simeq}\displaystyle\varinjlim_{s\in \Phi}
(s^{-1}A/A)\lraf{\simeq}\displaystyle\varinjlim_{s\in \Phi}
(A/sA).$$ Similarly,
$\Phi^{-1}B/B\lraf{\simeq}\displaystyle\varinjlim_{s\in \Phi}
(B/sB)$. Since $A/sA\lraf{\simeq} B/sB$ for each $s\in\Phi$, the map
$f$ induces an isomorphism of $A$-modules:
$\Phi^{-1}A/A\lraf{\simeq}\Phi^{-1}B/B$, that is,
$\Coker(\lambda)\simeq\Coker(\phi)$.

Finally, we point out that the map $\lambda': B\to
\Phi^{-1}A\otimes_AB$, defined by $b\mapsto 1\otimes b$ for $b\in
B$, is injective and that $\Coker(\lambda)\simeq
\Coker(\lambda')$. This is due to the equality $\phi=\lambda' h$.

Thus $$ 0\lra A\lraf{(-\lambda,\mu)}\Phi^{-1}A\oplus
\;B\lraf{{\left({\mu\,'\atop{\lambda\,'}}\right)}}
\Phi^{-1}A\otimes_AB \lra 0
$$ is an exact sequence of $A$-modules, where $\mu': \Phi^{-1}A\to
\Phi^{-1}A\otimes_AB$ is defined by $x\mapsto x\otimes 1$ for $x\in
\Phi^{-1}A$. By definition, the pair $(\lambda,\mu)$ is exact.

Since $\Phi$ consists of central, nonzerodivisor elements in $A$,
the $A$-module $\Phi^{-1}A$ is flat. This implies that $\Tor^A_i(B,
\Phi^{-1}A)$ = $0$ for all $i>0$, and therefore Corollary \ref{Karoubi} follows
from  Theorem \ref{K-theory} (1) immediately. $\square$

\subsection{$K$-theory of free products of groups\label{groupring}}

Finally, we apply Theorem \ref{K-theory} to algebraic $K$-theory of group rings.
As a preparation, we first recall some definitions and results from \cite{wald1, xc3} about
pure extensions.

Recall that an extension $R\subseteq C$ of rings is called \emph{pure} if there exists a spliting $C=R\oplus X$ of $R$-$R$-bimodules. The actual splitting is not part of the data, just its existence. In general, the $R$-$R$-bimodule $X$ may not be unique. For example, for a group $G$, the canonical embedding $R\subseteq RG$ is pure. In this case, $X$ has two natural choices. One is the free $R$-submodule of $RG$ generated by the nonidentity elements of $G$. The other is  the kernel of the canonical surjective ring homomorphism  $$\delta_G:\; RG\lra R, \quad\sum_{g\in G} r_g g\mapsto \sum_gr_g$$ where $r_g\in R$. The latter motivates the following definition. A pure extension $R\subseteq C$ is said to be \emph{strictly pure} if the $R$-$R$-bimodule $X$ is even a $C$-$C$-bimodule. In other words, $X$ is an ideal of $C$ such that the composite of the inclusion $R\to C$ with the canonical surjection $C\to C/X$ is the identity map.

Pure extensions were originally used by Waldhausen to study algebraic $K$-theory of generalized free products in \cite{wald1}. Now we briefly recall some of the results there.

Let $\alpha:R\to C$ and $\beta: R\to D$ be two pure extensions of rings. We denote by  $C\sqcup_RD$ the coproduct of $\alpha$ and $\beta$ in the category of $R$-rings. In \cite{wald1}, coproducts of rings were called generalized free products. Note that $\alpha$ and $\beta$ give rise to a map  $$\big(K(\alpha), -K(\beta)\big): K(R) \lra K(C) \times K(D)$$ of $K$-theory spaces, where $-K(\beta)$ denotes the map $K(D[1]\cpx{\otimes}_R-):K(R)\to K(D)$ induced from the functor $D[1]\cpx{\otimes}_R-:\Cb{\pmodcat R}\to \Cb{\pmodcat D}$. Since $K([1]):K(R)\to K(R)$ is a homotopy equivalence and a homotopy inverse of $K(R)$ (see the statements at the end of Subsection \ref{fpairs}), we see that both $\big(K(\alpha), -K(\beta)\big)$ and $\big(-K(\alpha), K(\beta)\big)$ have the same homotopy fibre up to homotopy equivalence.

Further, we fix two split decompositions of $R$-$R$-bimodules:
$$
C=R\oplus X\;\;\mbox{and}\;\; D=R\oplus Y.
$$

In order to describe the $K$-theory space of the coproduct $C\sqcup_RD$, Waldhausen introduced a topological space $\widetilde{K}\mathfrak{Nil}(R,X,Y)$ in \cite{wald1}, of which the homotopy type depends only on the ring $R$ and the $R$-bimodules $X$ and $Y$. For the original definition of $\widetilde{K}\mathfrak{Nil}(R,X,Y)$, we refer to \cite[Page 217]{wald1}; for further explanations of this space, one may find in \cite[Section 2.4]{Ranicki} and \cite[Section 0.2]{DKR}. The $n$-th algebraic $K$-group of $\widetilde{K}\mathfrak{Nil}(R,X,Y)$, usually called the $n$-th reduced Nilgroup, will be simply denoted by $\widetilde{\mathfrak{Nil}}_n(R,X,Y)$.

\begin{Lem}{\rm \cite[Theorem 1 and Theorem 4]{wald1}}\label{wd}
Suppose that $X$ and $Y$ are free right $R$-modules. Then the loop space $\Omega\big(K(C\sqcup_RD)\big)$ of the $K$-theory space of the ring $C\sqcup_RD$ is the direct product, up to homotopy, of the space $\widetilde{K}\mathfrak{Nil}(R,X,Y)$ and of the homotopy fibre of the map $\big(K(\alpha), - K(\beta)\big)$. Moreover, if the ring $R$ is regular coherent, then the space $\widetilde{K}\mathfrak{Nil}(R,X,Y)$ is contractible.
\end{Lem}

From the first part of Lemma \ref{wd}, we obtain the following Mayer-Vietoris exact sequence of $K$-groups:
$$
\cdots\ra K_n(R)\oplus\widetilde{\mathfrak{Nil}}_n(R,X,Y)
\ra K_n(C)\oplus K_n(D) \ra K_n(C\sqcup_RD)\ra K_{n-1}(R)\oplus\widetilde{\mathfrak{Nil}}_{n-1}(R,X,Y) \ra
$$
$$\cdots\ra K_1(R)\oplus\widetilde{\mathfrak{Nil}}_1(R,X,Y) \ra K_1(C)\oplus K_1(D)\ra
K_1(C\sqcup_RD)\ra K_0(R)\oplus\widetilde{\mathfrak{Nil}}_0(R,X,Y)\ra K_0(C)\oplus K_0(D) $$
where $K_n(C\sqcup_R D)\to \widetilde{\mathfrak{Nil}}_{n-1}(R,X,Y)$ is a split surjection for
$n\geq 1$.

\medskip

Now, let us reveal how the homotopy fibre of the map $\big(K(\alpha), - K(\beta)\big)$ can be related to noncommutative tension products when $\alpha$ and $\beta$ are strictly pure. This is based
on a construction of exact contexts in \cite[Section 4.2.2]{xc3}.

Assume that $\alpha$ and $\beta$ are strictly pure. Then the pair $(\alpha,\beta)$ can be completed into an exact context in the following way:
Let $M=R\oplus X\oplus Y$, the direct sum of abelian groups. We endow $M$ with the following multiplication:
$$
(r_1+x_1+y_1)(r_2+x_2+y_2):=r_1r_2+(r_1x_2+x_1r_2+x_1x_2)+(r_1y_2+y_1r_2+y_1y_2)
$$
for $r_i\in R$, $x_i\in X$ and $y_i\in Y$ with $i=1,2$. Under this multiplication, $M$ is a ring with identity $1$, and contains both $C$ and $D$ as subrings. Moreover, the quadruple
$(\alpha, \beta, M, 1)$ is an exact context. Now, we identify the noncommutative tensor product $D\boxtimes_RC$ with $R\oplus X\oplus Y\oplus Y\otimes_RX$ as $R$-$R$-bimodules. Then the multiplication of $D\boxtimes_RC$ is given by
$$
\big(r_1+x_1+y_1+y_3\otimes x_3\big)\circ \big(r_2+x_2+y_2+y_4\otimes x_4\big)$$
{\small
$$
=r_1r_2+(r_1x_2+x_1r_2+x_1x_2)+(r_1y_2+y_1r_2+y_1y_2)+\big(y_1\otimes x_2+y_3\otimes
(x_3r_2)+(r_1y_4)\otimes x_4+(y_1y_4)\otimes x_4+y_3\otimes (x_3x_2)\big).
$$}
where $r_1,r_2\in R$, $x_i\in X$ and $y_i\in Y$ for $1\leq i\leq 4$.

Note that $(\alpha, \beta, M, 1)$ is homological if and only if $\Tor^R_i(Y,X)=0$ for all $i\geq 1$. In particular, $(\alpha, \beta, M, 1)$ is homological if $Y_R$ or $_RX$ is free. In this case, by Corollary \ref{HTP}, the homotopy fibre of the map $\big(-K(\alpha), K(\beta)\big): K(R)\to K(C) \times K(D)$ (and thus also the map $\big(K(\alpha), -K(\beta)\big)$) is homotopy equivalent to the loop space $\Omega\big(K(D\boxtimes_RC)\big)$.
So, the following result follows immediately from Lemma \ref{wd}, Corollary \ref{HTP} and the fact that $C\sqcup_RD\simeq D\sqcup_RC$ as rings.

\begin{Koro}\label{Regular}
Let $\alpha:R\to C$ and $\beta:R\to D$ be strictly pure such that $X$ and $Y$ are free right $R$-modules. Then the following hold:

$(1)$ There is a homotopy equivalence: $$\Omega\big(K(C\sqcup_RD)\big)\lraf{\sim} \widetilde{K}\mathfrak{Nil}(R,X,Y)\times \Omega\big(K(D\boxtimes_RC)\big).$$
In particular, $K_n(C\sqcup_RD)\simeq \widetilde{\mathfrak{Nil}}_{n-1}(R,X,Y)\oplus K_n(D\boxtimes_RC)$ for $n\geq 1$.

$(2)$ If $R$ is regular coherent, then there are homotopy equivalences: $$\Omega\big(K(C\sqcup_RD)\big)\lraf{\sim}\Omega\big(K(D\boxtimes_RC)\big)\lraf{\sim}\Omega\big(K(C\boxtimes_RD)\big).$$
In particular, $K_n(C\sqcup_RD)\simeq K_n(D\boxtimes_RC)\simeq K_n(C\boxtimes_RD)$ for $n\geq 1$.
\end{Koro}

Now, let us consider pure extensions from group rings.
Suppose that $H$ and $G$ are two groups. Let $RH$ and $RG$ be the group rings of $H$ and $G$ over $R$, respectively. We take $C:=RG$ and $D:=RH$. Let $\alpha$ and $\beta$ be the canonical inclusions.
Then $\alpha$ and $\beta$ are strictly pure. Let $X$ and $Y$ be the kernels of the ring homomorphisms $\delta_G:RG\to R$ and $\delta_H:RH\to R$, respectively. Then $RG=R\oplus X$ and $RH=R\oplus Y$ as $R$-$R$-bimodules. Note that $Y$ and $X$ are free $R$-modules with $R$-basis $\{h-1\mid h\in H\setminus \{e_H\}\}$ and $\{g-1\mid g\in G\setminus \{e_G\}\}$, respectively, where $e_G$ denotes the identity of the group $G$. In this case, the multiplication of the ring $RH\boxtimes_RRG$ is exactly the one defined in Section \ref{Section 1}. We leave checking the details to the reader.
Note that the construction of $RH\boxtimes_RRG$ still makes sense if both $H$ and $G$ are semigroups with identity element.

\medskip
{\bf Proof of Corollary \ref{GP}.}

Let $\mathscr{G}$ be the category of groups, and let $\mathscr{R}$ be the category of $R$-rings. Recall that an $R$-ring is a ring $U$ with identity and a ring homomorphism from $R$ to $U$ preserving identity. The group ring functor
$$R(-):\;\mathscr{G}\lra\mathscr{R},\; G\mapsto RG\quad \mbox{for}\quad G\in\mathscr{G}$$ is left adjoint to the functor which sends an $R$-ring to its group of invertible elements. So the functor $R(-)$ preserves coproducts. Since the group $H*G$ is the coproduct of $H$ and $G$ in the category $\mathscr{G}$, we see that $R(H*G)$ is the coproduct of $RH$ and $RG$ in the category $\mathscr{R}$. Thus $R(H*G)=RH\sqcup_RRG$. Now, Corollary \ref{GP} follows from Corollary \ref{Regular} (2). $\square$

\smallskip
Finally, we apply Corollary \ref{GP} to fundamental groups of topological spaces.

Let $U$ be a topological space which is the union of two open and path connected subspaces $U_1$ and $U_2$. Suppose that $V:=U_1\cap U_2$ is path connected and nonempty. Let $x\in V$ be a point.
We consider the fundamental groups $\pi_1(U),\pi_1(U_1),\pi_1(U_2)$ and $\pi_1(V)$ of
$U$, $U_1$, $U_2$ and $V$ at $x$. By the Seifert-van Kampen theorem, the diagram of fundamental groups
$$
\xymatrix{\pi_1(V)\ar[r]\ar[d] & \pi_1(U_1)\ar[d]\\
\pi_1(U_2)\ar[r] &\pi_1(U)
}
$$
is a pushout in the category of groups. Let $R(\pi_1(U))$ be the group ring of $\pi_1(U)$ over the ring $R$.

As a consequence of Corollary \ref{GP}, we have the following result.

\begin{Koro}\label{fgroups}
Let $R$ be a regular coherent ring (for example, the ring $\mathbb{Z}$ of integers). Suppose that $\pi_1(V)$ is trivial. Then $$ K_n\big(R(\pi_1(U))\big)\simeq K_n\big(\,R(\pi_1(U_2))\boxtimes_R R(\pi_1(U_1))\,\big)\simeq K_n\big(\,R(\pi_1(U_1))\boxtimes_R R(\pi_1(U_2))\,\big)\quad \mbox{for all}\;\; n\geq 1.$$
\end{Koro}

\begin{Rem} Corollaries \ref{GP}, \ref{Regular} and \ref{fgroups} can be applied to some cases considered in \cite{Ranicki,DKR}. For example, for an arbitrary ring $R$, we have $$K_n\big(R(D_\infty)\big)\simeq K_n(R\,\mathbb{Z}_2\boxtimes_R R\,\mathbb{Z}_2)\oplus
\widetilde{\mathfrak{Nil}}_{n-1}(R),$$where $\mathbb{Z}_2$ is the group of order $2$. In fact, we know that $D_{\infty}=\mathbb{Z}_2*\mathbb{Z}_2$ and $R(\mathbb{Z}_2*\mathbb{Z}_2)=R\mathbb{Z}_2\sqcup_RR\mathbb{Z}_2$. It follows from Corollary \ref{Regular}(1) that $K_n\big(R(D_\infty)\big)\simeq K_n(R\,\mathbb{Z}_2\boxtimes_R R\,\mathbb{Z}_2)\oplus \widetilde{\mathfrak{Nil}}_{n-1}(R,R,R)$ for $n\geq 1$. Note that $\widetilde{\mathfrak{Nil}}_{*}(R,R,R)\simeq \widetilde{\mathfrak{Nil}}_{*}(R)$ by \cite[Corollary 3.27 (1)]{DKR}, where the reduced Nilgroups $\widetilde{\mathfrak{Nil}}_{*}(R)$ appears in algebraic $K$-groups of the polynomial ring $R[x]$ with one variable $x$: $K_*(R[x])\simeq K_*(R)\oplus \widetilde{\mathfrak{Nil}}_{*-1}(R)$. Thus $K_n\big(R(D_\infty)\big)\simeq K_n(R\,\mathbb{Z}_2\boxtimes_R R\,\mathbb{Z}_2)\oplus
\widetilde{\mathfrak{Nil}}_{n-1}(R)$, as desired. This is different from the decomposition: $K_n\big(R(D_\infty)\big)\simeq \big(K_n(R\,\mathbb{Z}_2)\oplus K_n(R\,\mathbb{Z}_2)\big)/K_n(R)
\oplus\widetilde{\mathfrak{Nil}}_{n-1}(R)$, given in \cite[Corollary 3.27 (2)]{DKR}.
\end{Rem}

We have considered strictly pure extensions in this section by using exact contexts. Now we mention the following open question for arbitrary pure extensions.

\medskip
{\bf Question.} Let $\alpha:R\to C$ and $\beta:R\to D$ be \emph{pure} extensions such that $X$ and $Y$ are free right $R$-modules. How can one describe the homotopy fibre of the map $\big(K(\alpha), -K(\beta)\big): K(R)\to K(C) \times K(D)$ in terms of algebraic $K$-theory space of a ring?


\medskip
{\bf Acknowledgements.}  Both authors would like to thank Marco Schlichting for email discussions on the definition of $K$-theory spaces of dg algebras. They are also grateful to Andrew Ranicki and Stefan Schwede for conversations during a min-workshop in Oberwolfach, 2013. The author H.X.Chen thanks Dong Yang for a discussion on the proof of Theorem \ref{new-theorem} during a visit to Nanjing Normal University in April, 2013. The corresponding author C.C.Xi thanks the Natural Science Foundation for partial support.


\medskip

{\footnotesize

\begin{tabular}{lcl}
Hongxing Chen, & \qquad \qquad\qquad & Changchang Xi \\
School of
Mathematical Sciences, BCMIIS, & \qquad &
School of
Mathematical Sciences, BCMIIS,\\
Capital Normal University & \qquad & Capital Normal University\\
100048 Beijing & \qquad & 100048 Beijing \\
People's Republic of China & \qquad & People's Republic of China\\
{\tt Email: chx19830818@163.com} & \qquad & {\tt Email:
xicc@cnu.edu.cn}
\end{tabular}

\bigskip
First version: August 25, 2012. Revised: April 3, 2013}
\end{document}